\documentclass[12pt]{article}
\newcommand{\R}{\mbox{I$\!$R}}

\newcommand{\qed}{{\hfill {$\rlap{$\sqcap$}\sqcup$}}\\[0.2in]\hspace*{0.5in}}
\newcommand{\qedwh}{{\hfill {$\rlap{$\sqcap$}\sqcup$}}\\[0.2in]}
\newcommand{\bk}{\\[0.05in] \hspace*{0.5in} }
\newcommand{\btd}{\bigtriangledown}
\newcommand{\mfor}{\ \ \ \ {\mbox{for}} \ \ }

\begin{document}
\pagenumbering{gobble}


\begin{center}
{\LARGE {\bf Conformal Scalar Curvature Equation on $S^n$\,:       }} \medskip \medskip \medskip  \smallskip   \\
{\LARGE {\bf       Functions  With Two Close Critical Points } } \medskip \medskip  \medskip  \medskip   \\
{\LARGE {\bf     (\,Twin Pseudo\,-\,Peaks)} }\medskip \medskip \smallskip \\
\end{center}

\vspace{0.8in}

\centerline{\Large  {Man Chun  {\LARGE L}EUNG${\,}^{*}$\ \ \ \&\ \ \ Feng \,{\LARGE Z}HOU${\,}^{**}$ }}

\vspace*{0.2in}

\centerline{\large {National University of Singapore}}

\vspace{0.73in}


\begin{abstract}
\vspace{0.33in}
\noindent By using  the Lyapunov\,-\,Schmidt reduction method without perturbation\,,\, we consider  existence results for the conformal scalar curvature on $\,S^n\,$ ($\,n \,\ge\,3\,$) when the prescribed function (after being  projected to $\,\R^n\,$) has two {\it close}\, critical points\,,\, which have the same value (positive), equal ``\,flatness\,"\  (\,`\,{\it twin}\,'\,;\, flatness \ $< \ n\,-\,2\,$)\,,\, and exhibit maximal behavior in certain directions (\,`\,{\it pseudo\,-\,peaks}\,'\,)\,.\, The proof relies on a balance between the two main contributions to the reduced functional \,-\, one from the critical points and the other from the interaction of the two bubbles.

\end{abstract}

\vspace*{1.2in}

{\bf 2000 AMS MS Classification}\,:\ \ Primary 35J60\,;\ \ Secondary 53C21. \\[0.05in]
$\overline{\ \ \ \ \ \ \ \ \ \ \ \ \ \ \ \ \ \ \ \ \  \ \ \ \ \ \ \ \ \ \ \ \ \ \ \ \ \ \ \ \ \ \ \ }$\\[0.03in]
\noindent{$\!\!{\!}^{*}$ \,{\ \tt{matlmc@nus.edu.sg}}  \\[0.075in]
\noindent{$\!\!{\!}^{**}$ \,{\tt{zhoufeng@u.nus.edu.sg}}  \\[0.075in]

\newpage
\pagenumbering{arabic}

{\bf \Large {\bf  1. \ \ Introduction.}}\\[0.1in]
As a counterpart of the Yamabe problem \cite{Brendle-1} \cite{Brendle-2} \cite{Compact} \cite{Obata} \cite{Schoen} \cite{Trudinger} (cf. also \cite{Bahri-Yamabe}\,)\,,\,   the prescribed scalar curvature problem in $\,S^n\,$ ($\,n\ \ge \ 3\,$\,)  asks for a positive solution $\,U\,$ to the nonlinear partial differential equation
$$
\ \ \ \ \ \ \ \Delta_1\, U\ -\ {\tilde c}_n\,n\,(n-1)\,U\ +\ (\,{\tilde c}_n\,{\cal K})\,U^{{n\,+\,2}\over {n\,-\,2}}\ =\ 0\ \ \ \ \ \ \ \ \ \ \ {\mbox{in}}\ \ \ S^n\ \ \ \ \ \ \ (U \ > \ 0)\,,\leqno (1.1)
$$
where    $\,{\cal K}\,$  is a prescribed function on $\,S^n\,.\,$ Here $\,{\tilde c}_n \,=\,(n\,-\,2)/\,[\, 4\,(\,n\,-\,1\,)\,]\ .\,$ See \,{\bf \S\,1\,d}\, for the rather standard notations we use\,.\, Also known as the  Nirenberg/\,Kazdan\,-\,Warner problem \cite{Tristan}\,,\, it can be compared to the classical Minkowski problem on prescribing Gaussian curvature for convex compact surfaces in $\,\R^3\,.\,$  The hallmark of equation (1.1) is the critical Sobolev exponent: the injection $\,H^{1\,,\,2} (S^n)\ \hookrightarrow\ L^{{2n}\over {n\,-\ 2}} (S^n)\,$ is not compact, typified by   blow\,-\,up gathering at critical point(s) of $\,{\cal K}$ \,.\, Close to half a century (\,cf. an early work in 1972   by Dimitri Koutroufiotis \cite{1972}, whose thesis adviser is   Louis Nirenberg\,)\,,\,  equation (1.1) serves as a vehicle for sophisticated techniques in nonlinear partial differential equations to be deployed and developed\,.\, It can also be branched out to complete manifolds, CR manifolds, Q\,-\,curvature, as well as related to mean field equations. See  some recent works \,\cite{Progress-Book} \cite{Chen-Xu-Flow} \cite{Chen-Lin-2} \cite{Chen-Lin-3} \cite{Chtioui} \cite{Leung-1st}  \cite{Li-1} \cite{Li-2} \cite{Lin-Wang} \cite{Simple-Symmetry} \cite{Pak} \cite{Pino} \cite{Wei-Yan}\, on the topic\,,\,  and the references therein\,.\, In general\,,\, existence results involve   symmetry on $\,{\cal K} \,$,\, or   local conditions on the critical points of $\,{\cal K}\,$ together with index inequality(ies)\,.\, The following result provides a good picture \{\,see \cite{Progress-Book} \cite{Chang-Yang-Duke}\,,\, and in particular \cite{Chen-Xu-Flow} regarding \,{\bf (iv)} below\}\,.\,  Assume the following {\bf (i)}\,--\,{\bf (iv)}\,.\\[0.1in]
  {\bf (i) } \ \ $\,{\cal K}\,$ is a smooth Morse function [\,namely, all its critical points \,(\,collected in the set\\[0.075in]
  \hspace*{.4in}  denoted by Crt\,) are non\,-\,degenerate\,]\,.\, \\[0.05in]
  {\bf (ii) } \ \, $\Delta_{1} \,{\cal K} \,(x_c) \,\not=\,0\,$ for \,all \,$x_c \,\in\,  {\mbox{Crt}} \,$.\, \\[0.075in]
  {\bf (iii) } \
$ \displaystyle{
\sum_{x_c \, \in \  {\mbox{Crt}}_{<}  } (-\,1)^{\,{\small\mbox{Ind}} \, ( \,{\cal K}\,, \ x_c)  } \ \not= \  (\,-1)^{\,n} \ }$,\,\ \  here $\,{\mbox{Crt}}_{<}   = \ \{\, x_c  \,\in\,  {\mbox{Crt}} \ | \ \Delta_{1} \,{\cal K} (x_c) \ <  \ 0\, \}$\,.\,\\[0.075in]
  {\bf (iv) } \   $\,{\cal K}\,$ is ``\,sufficiently\," close to a positive constant. \\[0.1in]
  Then equation (1.1) has a positive solution (precisely, see Theorem 7.1 in \cite{Progress-Book}\,,\, pp. 103\,)\,.\,  Recall that the index of a non\,-\,degenerate critical point is the number of {\it negative\,} eigenvalues of the Hessian matrix at that point.\bk
We observe that in case $\,{\mbox{Crt}}_{<}\,$ contains only {\it one\,} point, say at the north pole $\,{\bf N}\,,\,$ then it must be the peak (maximal point),  and hence (together with the non\,-\,degenerate condition)
$$
{{\small\mbox{Ind}} \, ( \,{\cal K}\,, \ {\bf N}\, )  }  \ = \ n \ \ \Longrightarrow \ \ \sum_{x_c \, \in \,  {\mbox{Crt}}_{<}  } (-\,1)^{{\small\mbox{Ind}} \, ( \,{\cal K}\,, \ x_c) }  \  = \   (-\,1)^{{\small\mbox{Ind}} \, ( \,{\cal K}\,, \ {\bf N}\, )  }  \ = \ (-1)^{\,n}\ .
$$
Moreover, via the Kazdan\,-\,Warner (Pohozaev) identity, if $\,{\cal K}\,$ is strictly decreasing from $\,{\bf N}\,$ to $\,{\bf S}\,,\,$ measured via the geodesics, then equation (1.1) does not have any positive solution at all (cf. also \cite{Han-Li})\,.\,

\newpage

\hspace*{0.5in}Motivated by this, attention is given  to the situation where $\,{\mbox{Crt}}_{<}\,$ contains at least {\it two\,} points\,.\, Cf. \cite{Degree=} \cite{Bianchi-1} \cite{Yan+} \cite{Yan} (\,a discussion on the existence and non\,-\,existence  results can be found in {\bf \S\,1\,c})\,.\, Thus in this article we consider juxtaposed `\,twin\,' pseudo\,-\,peaks (described in {\bf{\S\,1\,a}}\,). We note that ``\,side\,-\,by\,-\,side\," is a kind of symmetry condition. To state the  local conditions on the Taylor expansions at the two pseudo\,-\,peaks\,,\, we introduce the stereographic projection $\,\hat{\cal P}\,$ [\,see (4.15)\,] which sends the north pole to infinity.  Equation (1.1) is transformed to
$$
\ \ \ \ \ \ \ \ \ \Delta\,v \ +\ (\,{\tilde c}_n\,K\,)\,v^{{n\,+\,2}\over{n\,-\,2}}\ = \ 0\ \ \ \ \ \ {\mbox{in}}\ \ \ \ \R^n\ \ \ \ \ \ (v \ > \ 0)\,. \leqno (1.2)
$$
See, for examples, \cite{I} \cite{Schoen-Yau-Book}\,. Here
\vspace*{-0.05in}
$$
\ \,v\,(y)\ :=\ U\,({\hat{\cal P}}^{-1}\,(y))\,\cdot \left({2\over {1\,+\,|\,y\,|^{2}}}\right)^{\!\!{{\,n\,-\,2\,}\over 2}}\ \  {\mbox{and}}\ \ \ K\,(y) \ = \ {\cal K}\,(x) \mfor y \ = \ \hat{\cal P}\,(x) \,\in\,\R^n\,. \leqno (1.3)
$$
(\,Note that max\ $K\,$ is not affected by whichever point we choose as the north pole.)

\vspace*{0.2in}

{\bf \S\,1\,a.}\ \ \ {\it Close\, `\,twin\,' pseudo\,-\,peaks and their key parameters.} \ \ Consider two critical points  $\,{\bf q}_1\,$ and $\,{\bf q}_{\,2}$\, of $\,K\,.\,$  Via a translation, one may assume without loss of generality that
$$
{\bf q}_1 \ = \ 0\,.\leqno (1.4)
$$
Let \,$\gamma$\, denote the distance (or gap) between the two critical points\,.\, Moreover, in this article, we always assume that $\,{\bf q}_1\,$ and $\,{\bf q}_{\,2}$\, are close, namely,
$$
\gamma \ := \ |\,{\bf q}_{\,2}\,| \ =  \ O\,(1)\,.\leqno (1.5)
$$
The two critical points are symmetric (or `\,twin\,') in the following sense [\,(1.6) \ \& \ (1.8)\,]\,:
$$
K\,(0) \ = \ K\,({\bf q}_2) \ >  \ 0 \leqno (1.6)
$$
[\,after a rescaling, we may accept without loss of generality that
$$
({\tilde c}_n\!\cdot K)\,(0) \ = \ ({\tilde c}_n\!\cdot K)\,({\bf q}_2) \ = \ n\,(n\,-\,2)\,]\,, \leqno (1.7)
$$
in conjunction with (similarity on the Taylor expansions)
$$
({\tilde c}_n\!\cdot K) \, (y) \ = \ n\, (n \,-\, 2)\   + \   {\bf P}_j^\ell \, (y\,-\, {\bf q}_j)  \
+ \ R_j^{\ell \ + 1}\,  (y\,-\, {\bf q}_j)       \ \ \ \ \ \ \ {\mbox{for}} \ \ \ y \,\in\, B_{{\bf q}_j}\, (\rho) \leqno (1.8)
$$
 and $\,j \ = \ 1\,,\, \ 2\,.\,$ Here
$$ \ \ \ \ \ \ \ \ \  \ \ \rho \ = \  \hbar \cdot \gamma \ \ \ \ \ \ \ \ \ \ \ \ (\, \hbar \ \ {\mbox{is \ \ a  \  \ fixed \ \ number \ \ less \ \ than \ \ half}}\,)\,, \leqno (1.9)_{\bf{(i)}}
$$
${\mbox{(1.9)}}_{ ({\bf{ii}}) }$ \ \ $\,{\bf P}_j^\ell\ $ is a {\it  homogeneous polynomial}\, of degree $\,\ell \ \ge \ 2\,$ (\,$\ell\ $-\, the flatness, is the same for  $\ j \ = \ 1\,,\, \ 2\,$)\,,\, and \\[0.1in]
${\mbox{(1.9)}}_{ ({\bf{iii}}) }$ \ \ $\, R_j^{\ell \ + 1} \,$ the remainder in the Taylor expansion\,,\, satisfying
\begin{eqnarray*}
(1.10) \ \ \ \ \ \ \ \ \  \ \ \ \   |\,R_1^{\ell \ + 1} \, (y)   \,| & \le &  C_{R_1} \cdot |\,y\,|^{\ell\,+\,1} \ \ \ \  \ \ \ \ \ \ \ \ \ {\mbox{for}} \ \ y \ \in \ B_o \,(\rho)\,, \\[0.15in]
|\,R_2^{\ell \ + 1} \, (y)   \,|  &  \le &  C_{R_2}  \cdot |\,y \ - \ {\bf q}_2\,|^{\ell\,+\,1} \ \ \ \  {\mbox{for}} \ \ y \ \in \ B_{{\bf q}_2} \,(\rho) \ .\ \ \ \ \ \ \ \ \ \ \ \ \ \ \ \ \ \ \ \ \ \ \ \ \ \ \ \ \ \ \ \ \ \ \ \
\end{eqnarray*}
Here $\,C_{R_1}\,$ and $\,C_{R_2}\,$ are positive constants.
(1.8) implies that
\begin{eqnarray*}
(1.11) \ \ \ \ \ \ \ \ \  \ \ \ \   |\,({\tilde c}_n\!\cdot K) \, (y)  \ - \ n\, (n \,-\, 2) \,| & \le &  C_{P_1} \cdot |\,y\,|^\ell \ \ \ \  \ \ \ \ \ \ \ \ \,{\mbox{for}} \ \ y \ \in \ B_o \,(\rho)\,, \\[0.15in]
 |\,({\tilde c}_n\!\cdot K) \, (y)  \ - \ n\, (n \,-\, 2) \,|  &  \le &  C_{P_2}  \cdot |\,y \ - \ {\bf q}_2\,|^\ell \ \ \ \  {\mbox{for}} \ \ y \ \in \ B_{{\bf q}_2} \,(\rho) \ .\ \ \ \ \ \ \ \ \ \ \ \ \ \ \ \ \ \ \ \ \ \ \ \ \ \ \ \ \ \ \ \ \ \ \ \
\end{eqnarray*}
Here the positive constant $\, C_{P_j}\,$ (\,$j\ = \ 1\,, \ 2\,$) is linked to the sum of the absolute values of the coefficients of ${\bf P}_j$\,.\, Assume that
$$
\ell \ \ \,{\mbox{is \ \ even\,,\, \ \ and \ \ let}} \ \ \,h_\ell \ = \ {1\over 2} \cdot \ell\ .
$$

\vspace*{-0.2in}

Hence
\begin{eqnarray*}
(1.12) \ \ \ \ \ \ \ \ \ \Delta_{\bar y}^{\!(h_\ell)} \, {\bf P}_j^\ell\, (\,\bar y\,) & = &  \Delta_{\bar y}\ (\,\cdot \cdot \cdot\, [\,\Delta_{\bar y} \,[\,\Delta_{\bar y}  \  {\bf P}^\ell_j\,(\,\bar y\,)\,]\,]\,) )  \ = \ {\varpi}_j \ \ \ \ \ \ \ \ \ \ \ \ \ (\,j \ = \ 1\,,\, \ 2\,) \ \ \ \ \ \ \ \ \ \\[0.1in]
& \ & \leftarrow \  \, h_\ell \ \,{\mbox{times}} \ \,\rightarrow
\end{eqnarray*}
is a number. Here $\,\bar y \ = \ y\,-\, {\bf q}_{\,j}\,.\,$
The key condition for the critical points $\,{\bf q}_1\,$ and $\,{\bf q}_{\,2}$\, to be called {\it pseudo\,-\,peaks}\, is the following\,:
$$
{\varpi}_j \  <  \ 0 \ \ \  \ {\mbox{for}} \ \ j \ =  \ 1\,, \ 2\,.\leqno (1.13)
$$
We add the following ``symmetry\," condition as well\,:
$$
{1\over {C_p}} \cdot |\,{\varpi}_1\,| \ \le \ |\,{\varpi}_2\,| \ \le \ C_p \cdot |\,{\varpi}_1\,| \,, \ \ \ \ \ \  \ \ \ \ \ \ \  \ \leqno (1.14)
$$

\vspace*{-0.15in}

and

\vspace*{-0.3in}

$$
{\varpi}_1\,  \ \ge \  -\,C_\omega\ .   \leqno (1.15)
$$
In (1.14) and (1.15)\,,\, $\,C_p\,$ and $\,C_\omega\,$ are positive constants.

\vspace*{0.15in}

{\bf Main Theorem 1.16.} \ \ {\it For $\,6 \ \le  \ n \ <  \ 10\,,\,$  let
$$
\,\ell \ \in [\,2\,, \ n\,-\,2\,)   \leqno (1.17)
$$
be an even integer $\,{\cal K} \, \in   C^{\ell \,+\,1}\,(S^n)\,,\,$ and $\,K\,$  the projection of  $\,{\cal K}\,$ to} $\,\R^n$\, {\it via\,} (1.3)\,.\, {\it Assume that\,}
 $$\,|\,K\,|\, \le \, {\bar C}_b\ \ \ \ \  {\it{in}} \ \ \R^n,\leqno (1.18)$$
  {\it and  $\,K\,$ has twin pseudo\,-\,peaks in the sense of } (1.7)\,,\, (1.8)\, {\it and} \ (1.13)\,,\,  {\it located at $\ {\bf q}_1\,=\,0$\, and} $\ {\bf q}_2\,\in\,\R^n.\,$  {\it Under the conditions in\,} (1.10)\,,\, (1.11)\,,\, (1.14) {\it and}\, (1.15)\,,\, {\it there is a positive constant $\,\gamma_o\,$ so that if }
$$
|\,{\bf q}_{\,2}\,| \ \le \ \gamma_o\,,\,
$$
{\it then equation\,} (1.1) {\it has a positive $\,C^2$-\,solution.} {\it Moreover,  $\,\gamma_o\,$ depends only on  $\,n\,$,\,   $\,\ell\,,\,$ $\,{\bar C}_b\,$,\, and the parameters of the twin pseudo\,-\,peaks}\, ({\it namely\,,\, $\hbar\,,\,$ ${C}_{R_1}$\,,\, ${C}_{R_2}$\,,\, ${C}_{P_1}$\,,\, ${C}_{P_2}$\,,\, $C_p\,$ {\it and\,} $\,C_\omega\,)$\,.\,}

\newpage

{\it Remarks.}\\[0.1in]
(1) \ \  To gain an idea on the dependence of $\,\gamma_o\,$ on   $\,C_\omega\,$ [\,appeared in (1.14)\,]\,,\,   we have
$$
\gamma_o \ \approx     {  {c_\mu  }\over { C_\omega^{1\over {\ell}}  }}   \ .
$$
Here the small positive number $\,c_\mu$\, depends on  the other parameters in Theorem 1.16\,.\, See {\bf \S\,5\,c\,.}

\vspace*{0.1in}

(2) \ \,With the help of Theorem 1.16, one can consider  multiple solutions for well-separated multiple twin pseudo\,-\,peaks\,.\,

\vspace*{0.1in}

(3) \ \  There is no condition on other critical points.

\vspace*{0.1in}

(4) \ \  Dimension restriction ($\,n \ =  \ 6\,,  \ 7\,, \ 8\, \ \& \ 9\,$) mainly due to the process when key information are extracted out of the reduced functional (\,refer to Proposition 4.1).

\vspace*{0.2in}

{\bf \S\,1\,b.}\ \ \ {\it  Lyapunov\,-\,Schmidt reduction method without perturbation. Organization.} \ \ As described in \cite{Progress-Book}\,,\,
the elegant Lyapunov\,-\,Schmidt reduction method is considered on those $\,K\,$ which is a perturbation of a positive constant, that is (after a rescaling),
$$
(\,{\tilde c}_n\,K\,) \ =  \ n\,(n \,-\,2) \ + \ \varepsilon \cdot (\,{\tilde c}_n\,H\,)\,.  \leqno (1.19)
$$
Here $\,\varepsilon\,$ is ``\,small enough\,". A new insight is introduced in \cite{Wei-Yan}, where Wei and Yan  bring home to the point that when a large number of standard bubbles are arranged near the critical points of $\,K\,$,\, one can still apply the Lyapunov\,-\,Schmidt reduction method, this time without the requirement on $\,\varepsilon\,$ being  close to zero (\,see also an earlier work of Yan \cite{Yan}\,)\,.\, Thus the number of bubbles replaces the parameter $\,\varepsilon\,.\,$ \bk
In this article, we show that by ``\,planting\," one bubble each near one of the twin pseudo\,-\,peaks, the   Lyapunov\,-\,Schmidt reduction method is also  applicable without the need for $\,K\,$ being close to a constant (\,{\bf \S\,2}   \,\&  {\bf \S\,3}\,)\,.\, In this case the ``\,gap\," $\,\gamma$\, take the place of the parameter  $\,\varepsilon\,.\,$ Moreover, we show that the reduced functional has two main contributions (\,Proposition 4.1\,;\, cf. also \cite{Yan+}\,)\,,\, one from the critical point (\,{\bf \S\,4}\,)\,,\, and the other one from the interaction with the other bubble (\,{\bf \S\,2\,b}\,). By properly balancing these two effects, we show that equation (1.2) has a solution if the peaks are close enough (\,{\bf \S\,5}\,)\,.\, This solution can be transferred back to $\,S^n\,$ via (1.3) as a solution of (1.1)\,.\,   Moreover, as the two bubbles are highly concentrated near the twin pseudo\,-\,peaks, other critical points (\,if any) do not contribute to the consideration\,.\,  This is   in harmony with a theme in \cite{Leung-Supported} (\,cf. also \cite{Leung-blow-up}\,) that concentration can be put to good use to find solutions of equation (1.1)\,.\,\bk

\newpage

{\bf \S\,1\,c.}\ \ \ {\it  Comparison with some related existence and non\,-\,existence results.} \ \
Our result should be compared to \cite{Yan+}\,,\, in which the authors use  a version of
Lyapunov\,-\,Schmidt reduction method for $\,\varepsilon\,$  small enough\,,\,   when $\,H\,$ in (1.19) has two critical points [\,among other possible critical point(s)\,]\,,\, say at $\,{\bf q}_{\,1}' \ = \ 0\ $ and $\ {\bf q}_{\,2}' \ = \ ({\bf q}_{\,2|_1}' \,,\ \cdot \cdot \cdot\,, \ {\bf q}_{\,2|_n}'\,)\,\in\,\R^n\,$ (not necessarily close)\,,\, which satisfy\\[0.1in]
(1.20)
\begin{eqnarray*}
H\,(y) & = & H\,(0) \ + \ \left(\,a_1 \,|\,y_1|^{\,\beta_1} \ + \ \cdot \cdot \cdot \ + \ a_n \,|\,y_n|^{\,\beta_1}\,\right) \ + \ O\left(|\,y_i|^{\,\beta_1 \,+\,\sigma_1} \right) \ \ \  \ \  \mfor y \,\in\,B_o(\rho)\,,\\[0.1in]
H\,(y) & = & H\,({\bf q}_{\,2}') \ + \    \left(\,b_1 \,|\,y_1 \ - \ {\bf q}_{\,2|_1}'|^{\,\beta_2}\ + \ \cdot \cdot \cdot \ + \ b_n \,|\,y_n \ - \ {\bf q}_{\,2|_n}'|^{\,\beta_2}\,\right) \\[0.1in]
 & \ &  \ \ \ \ \ \  \ \ \ \  \ \ \ \ \ \ \ \ \ + \ O\left(|\,y_i \ - \ {\bf q}_{\,2|_i}'  |^{\,\beta_2 \ +\ \sigma_2} \right) \ \ \ \ \       \mfor y \,\in\,B_{{\bf q}_{\,2}'}(\rho) \ \ \ \ \  \ \left[\,\rho \ < \  {{ |\,{\bf q}_{\,2}'\, |}\over 2}\,  \right]\,,\\[-0.3in]
\end{eqnarray*}
where $\,\beta_1\,, \ \beta_2 \ \in \ (\,0\,, \ n\,-\,2\,)\,,\,$
$$
a_i\  \not= \ 0\,, \ \ b_i\  \not= \ 0\,, \ \  \sum_{i \,=\,1}^n a_i \ <  \ 0\ \ \ \  \ {\mbox{and}} \ \ \ \  \sum_{i \,=\,1}^n\, b_i \ <  \ 0\ \ \ \ \  \ \ (\,i \ = \ 1\,, \ 2\,, \ \cdot \cdot \cdot\,, \ n)\,,
$$
then for $\,\varepsilon\,$ in (1.19) small enough, equation (1.2) has a (two peaks) solution (see Theorem 1.1 in \cite{Yan+} for the precise description). In the above, $\,\sigma_1\,, \ \sigma_2 \ \in \ (0, \ 1\,)$  are fixed numbers\,.\, Besides the requirement on $\,\varepsilon\,$  being small enough\,,\, we note that in (1.20), there is no cross over terms like $\  y_1 \times y_2 \cdot \cdot \cdot\,$,\, which is allowed in our Main Theorem 1.16. \bk
In \cite{Yan}\,,\, a counterpart to the situation above is considered. There Yan studies the case when $\,K\,$ has
a pair of strictly local maximum points  at $\,{\bf m}_1\,$ and $\,{\bf m}_2\,,\,$ whose distance $\,|\,{\bf m}_1\,-\,{\bf m}_2\,|\,$ is very {\it large} [\,flatness of  these two local maxima  is in the range $(n\,-\,2\,, \ n)$\,]\,.\, See Theorem 1.1 in \cite{Yan} for the complete statements.\bk
On the other hand, a non\,-\,existence result obtained by Bianchi in \cite{Bianchi-Non} suggests that for certain ``\,very sharp\," twin peaks with flatness lesser than or equal to $\,n\,-\,2\,$,\,\, equation (1.1) has no positive solution. For details, see \cite{Bianchi-1} \cite{Bianchi-Non}\,.\, Cf. also \cite{Simple-Symmetry}\,.\,  Thus the smallness of $\,\gamma\,$ in the Main Theorem cannot be totally removed.

\vspace*{0.2in}

{\bf \S\,1\,d.}\ \ \ {\it General conditions, assumptions and conventions.}\ \ Throughout this work,\,
$$
\ \ S^n\ =\ \{ \,x\,=\,(\,x_1\,,\ \cdots\ ,\ x_{n +1})\ \in\ \R^{n+1}\,\ |\ \,x_1^2\,+\,\cdots\,+\,x_{n+1}^2\ =\ 1\,\} \ \ \ \ \ (\,n\,\ge\,3)\,, \leqno (1.21)
$$
with the induced metric $\,g_1\,$\,. $\,\Delta_1\,$ is the Laplace\,-\,Beltrami operator associated with $\,g_1\,$ on $\,S^n\,$. Likewise, $\,\Delta\,$ is the Laplace\,-\,Beltrami operator associated with Euclidean metric $\,g_o\,$ on $\,\R^n\,$,\, with coordinates $\,y \ = \  (y_1\,,\ \cdots\ ,\, y_{n})\,\in \,\R^n\,.\,$ Moreover, the norm $\,\Vert\ \ \Vert\,$ and the inner product $\,\langle\ \,,\ \rangle\,$ are defined via Euclidean metric $\,g_o\,$ on $\,\R^n\,$.\\[0.1in]

\newpage

$\bullet_1$\ \, As mentioned earlier, $\,{\tilde c}_n\,=\,{{n\,-\,2}\over{4\,(n\,-\,1)}}$\,. We  observe the practice on  using `$\,C\,$'\,,\, possibly with sub\,-\,indices, to denote various positive constants, which may be rendered {\it differently\,} from line to line according to contents. {\it Whilst we use `$\,{\bar c}$' \,or\, `$\,{\bar C}$'\,,\, possibly with sub\,-\,indices, to denote a\,} fixed\, {\it positive constant which always keeps the same value as it is first defined\,}.\,\\[0.05in]
$\bullet_2$\ \, Denote by $\,B_y\,(r)\,$ the open ball in $\,(\R^n\,,\,g_o)\,$ with center at $\,y\,$ and radius $\,r\,>\,0\,$, and $\,\partial B_y\,(r)\,$ its boundary\,.\, Whenever there is no risk of misunderstanding, we suppress $\,dy\,$ from the integral expressions on domains in $\,\R^n\,$.

\vspace*{0.2in}

{\bf \S\,1\,e.}\ \ \ \textbf{e}\,-{\it Appendix.}\ \ \ Some of the preparatory estimates are situational modifications of well\,-\,established arguments. We gather those details in the \textbf{e}\,-Appendix, which is presented from pp. 36 onward.

\vspace*{0.35in}

{\bf \Large {\bf  \S\,2.    }} {\bf  \Large{  The Lyapunov\,-\,Schmidt reduction scheme sans }}\\[0.1in]
 \hspace*{0.66in} {\Large{\bf perturbation\,:\, the case of  two bubbles.}} \\[0.2in]
Equation (1.2) is naturally associated  with  the Hilbert space
$$
{\cal D}^{1,\,\,2}\,=\,{\cal D}^{1,\,\,2}\,(\R^n)\,:=\,\left\{\,f\in L^{{2n}\over {n\,-\,2}}\,(\R^n)\:\bigcap\: W^{1,\, 2}_{\small\mbox{loc}}\, (\R^n)\ \,\bigg\vert\,\ \int_{\R^n}\langle\,\btd\,f\,,\,\btd\,f\,\rangle\,<\,\infty\,\right\}\,.\leqno (2.1)
$$
The inner product is defined by
$$
\langle\,f\,,\ \psi\,\rangle_{\,\btd}\ :=\ \int_{\R^n}\langle\,\btd\,f\,,\,\btd\,\psi\,\rangle\ \ \ \ \mfor\ f\,,\ \psi\ \in\ {\cal D}^{1,\,\,2}\,,  \ \ \ {\mbox{and}} \ \ \  \Vert\,f \,\Vert_{\btd}^2 \ := \ \langle\,f\,,\ f\,\rangle_{\,\btd} \ .  \leqno (2.2)
$$
The functional corresponding to (1.2) is given by
$$
\ {\bf I}\,(f)\ = \ {1\over 2}\,\int_{\R^n}\,\langle\,\btd\,f\,,\,\btd\,f\,\rangle\ -\ \left({{n\,-\,2}\over {2n}}\,\right)\,\cdot\,\int_{\R^n}\,(\,{\tilde c}_n\!\cdot K\,)\,f_+^{{2n}\over {n\,-\,2}}\ \ \ \mfor f\,\in\,{\cal D}^{1,\,\,2}\,.\leqno (2.3)
$$
Here $\,f_+\,$ denotes the positive part of $\,f\,$.\, See Part I \cite{I} on the regularity of the critical points of (2.3). Cf. also \cite{Brezis-Kato} in relation to equation (1.1)\,.\,  Let
$$
(\,{\tilde c}_n\!\cdot K\,) \ = \  n\,(n\,-\,2) \ + \ (\,{\tilde c}_n\!\cdot H\,)  \ \ \  \ \Longleftrightarrow \ \ \ (\,{\tilde c}_n\!\cdot H\,) \ = \ (\,{\tilde c}_n\!\cdot K\,)  \ - \  n\,(n\,-\,2)\,.\leqno (2.4)
$$
Accordingly\,,\, $\,{\bf I}\,$ can be split into two parts
$$
{\bf I} \,(f)\ = \ {\bf I}_o\,(f)\ +\ {\bf G}\,(f)\,,\leqno (2.5)
$$

\vspace*{-0.2in}

\begin{eqnarray*}
(2.6)\ \ \ \ \ \ \ {\mbox{where}} \ \  \ \ {\bf I}_o\,(f) & = & {1\over 2}\,\int_{\R^n} \langle\,\btd\,f\,,\,\btd\,f\,\rangle\ -\ n\,(n\,-\,2) \cdot \left({{n\,-\,2}\over {2n}}\, \right) \int_{\R^n} f^{{2n}\over {n\,-\,2}}\ ,\ \ \ \ \\[0.1in]
(2.7)\ \ \ \ \ \ \ \ \ {\mbox{and}} \ \ \ \ \   {\bf G}\,(f) & = &  -\ {{ n\,-\,2 }\over {2\,n}}\cdot\int_{\R^n} (\,{\tilde c}_n\!\cdot H\,) \,f^{{2n}\over {n\,-\,2}}\hspace*{1.2in} \mfor\ \ f\ \in \ {\cal D}^{1,\,\,2}.
\end{eqnarray*}

\newpage

Here we pay special attention on the negative sign in $\,{\bf G}\,(f)\,.\,$ One of the key themes in this article is to expound the interaction between   $\,{\bf I}_o'\ $ and $\,{\bf G}'$\,.\, \bk
Let us present the following flow chart to guide our discussion.

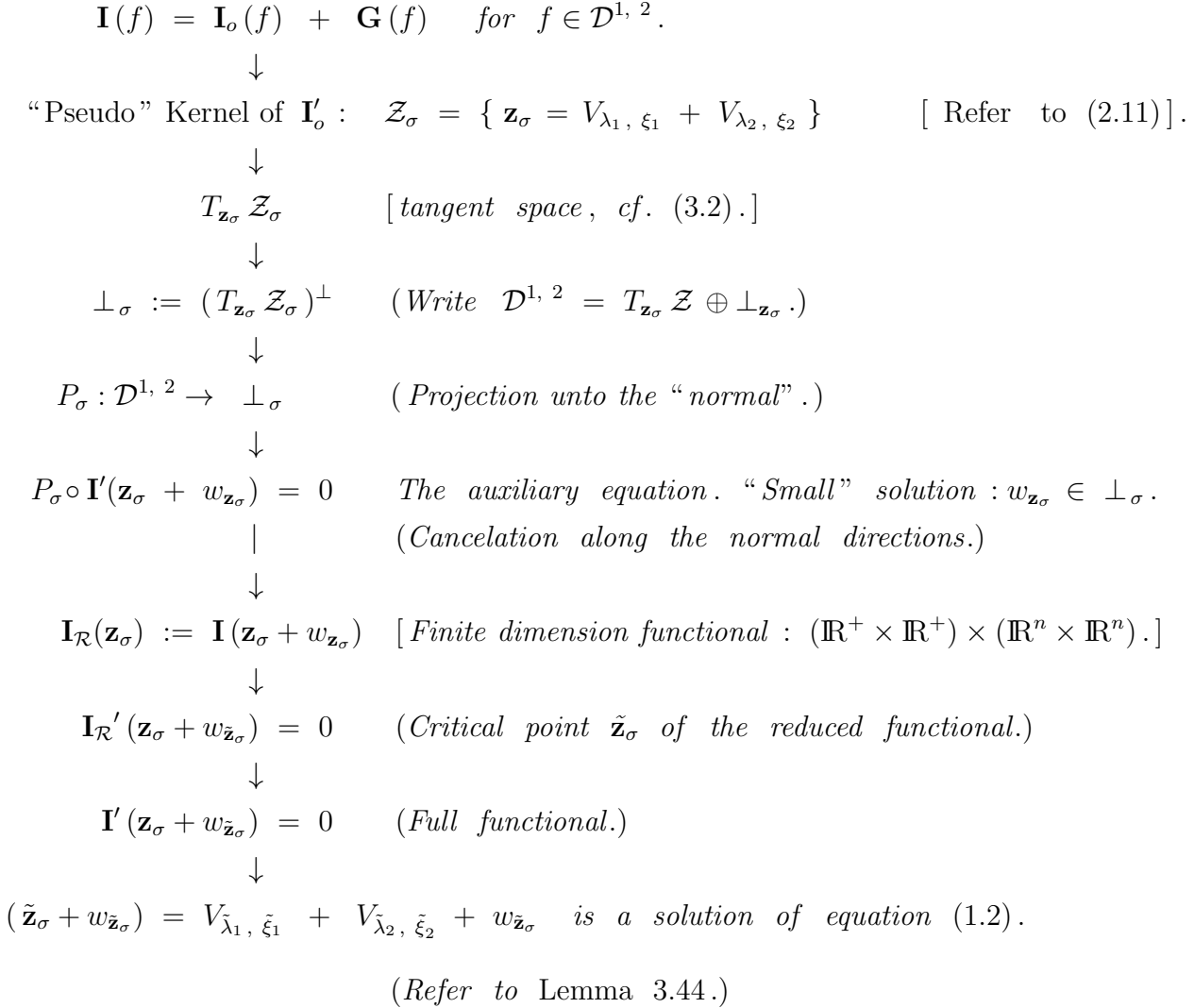
\vspace*{-0.2in}

\begin{eqnarray*}
{\bf I} \,(f)\ = \ {\bf I}_o \,(f) &+&   {\bf G}\,(f)\ \ \  \ \ {\it for} \ \ f \in {\cal D}^{1, \  2}\,.\\
\downarrow  \ \ &\,&\\
 {\mbox{``\,Pseudo\," Kernel of}} &\,&  \!\!\!\!\!\!\!\!\!\!\!{\bf I}_o'\,: \ \ \ {\cal Z}_\sigma \ = \ \{\ {\bf z}_\sigma \, = \, V_{\lambda_1\,,\ \xi_1} \ + \ V_{\lambda_2\,,\ \xi_2}\ \}\ \  \ \ \ \ \ \  \ \  [\,{\mbox{ Refer \ \ to} \ \ } (2.11)\,]\,.\\
\downarrow  \ \ &\,&\\
T_{{\bf z}_\sigma }\  \!{\cal Z}_\sigma   &\,& \ \ \ [{\it{\,tangent \ \ space\,,\,\  cf. }}\ \ (3.2)\,.\,] \\
\downarrow  \ \ &\,&\\
 {{\perp}_{\,\sigma }} \ := \ (\,T_{{\bf z}_\sigma }\  \!{\cal Z}_\sigma\,)^\perp\!\!\!\!\!\!\!\!\!\! &\,&\ \ \ \,( {\it{Write}} \ \ \,{\cal D}^{1, \  2}\ = \ T_{{\bf z}_\sigma}\  \!{\cal Z} \,  \oplus  {{\perp}_{{\bf z}_\sigma }}\,. ) \\
\downarrow  \ \ &\,&\\
P_{ \sigma }: {\cal D}^{1, \  2}   \to \ \ {{\perp}_{\,\sigma }}  & \ & \ \ \ \,({\it{\,Projection \ unto \ the \ ``\,normal"\,.\,}})\\
\downarrow  \ \ &\,&\\
P_{ \sigma }\!\circ {\bf I}' ({\bf z}_\sigma\  + \ w_{{\bf z}_\sigma })\ = \ 0  \!\!\!\!\!\!\!\!\!\! &\,& \ \ \ \  {\it The \ \ auxiliary \   \ equation\,. \ \ ``\,Small\," \ \ solution}\,: w_{{\bf z}_\sigma } \, \in \ {{\perp}_{\,\sigma }} \,.\\  |\ \,\,\,&\,&  \ \ \ \ ({\it Cancelation \ \ along \ \ the \ \ normal \ \ directions.)} \\
\downarrow  \ \ &\,&\\
{\cal {\bf I}_{R}} ({{\bf z}_\sigma })\ := \ {\bf I}\, ({{\bf z}_\sigma } + w_{{\bf z}_\sigma }) \!\!\!\!\!\!\!\!\!\!\!\!\!\!\!\!&\,& \ \   \ \ [\,{\it Finite \ dimension \   functional} \,: \   (\R^+ \times \R^+) \times (\R^n \times \R^n)\,.\,]\\
\downarrow  \ \ &\,&\\
{\cal {\bf I}_{R}}' \,({{\bf z}_\sigma } + w_{{\tilde{\bf z}}_\sigma }) \ = \ 0  \!\!\!\!\!\!\!\!\!\! &\,& \ \ \ \ ({\it Critical \ \ point \ \ {\tilde{\bf z}}_\sigma  \ \ of \ \ the \ \ reduced \ \ functional}.)  \\
\downarrow  \ \ &\,&\\
{\bf I}' \,({{\bf z}_\sigma } + w_{{\tilde{\bf z}}_\sigma })\ = \ 0    \!\!\!\!\!\!\!\!\!\! &\,& \ \ \ \ ({\it Full \ \ functional}.)  \\
\downarrow  \ \ &\,&\\
(\,{{\tilde{\bf z}}_\sigma } + w_{{\tilde{\bf z}}_\sigma }) \ = \ V_{{\tilde{\lambda}}_1\,,\ {\tilde{\xi}}_1} & + & V_{{\tilde{\lambda}_2\,,\  \tilde{\xi}}_2} \ + \   w_{{\tilde{\bf z}}_\sigma }   \ \  \ \, {\it is  \ \ a   \ \ solution \ \ of \ \ equation} \ \ {\mbox{(1.2)}}\,. \\[0.15in]
 & \,&  \ \ \ \,({\it Refer \ \ to} \ \ {\mbox{Lemma\  \ 3.44}}\,.)
\end{eqnarray*}

\centerline{--\ {\it{ Flow Chart of the Lyapunov\,-\,Schmidt reduction scheme without perturbation}.}\ --}

\vspace*{0.3in}

{\bf \S\,2\,a.}\ \ {\it First order property \,-\, interaction between two `\,well\,-\,separated' bubbles.}\ \ For $\,f\,\in\,{\cal D}^{1,\,\,2} \,$, a calculation using (2.6) shows that the Fr\'echet derivative of $\,{\bf I}_o\,$ at $\,f\,$ is given by
\vspace*{-0.02in}
$$
\ \ \ \ {\bf I}_o'\,(f)\,[\,h\,]\ = \ \int_{\R^n} \!\left[\,\langle\,\btd\,f\,,\,\btd\,h\,\rangle\,-\,n\,(n\,-\,2)\,f^{{\,n\,+\,2\,}\over{n\,-\,2}}_+ \cdot h\,\right]\ \ \ \ \ \ \ \ \mfor\ \ h\ \in\ {\cal D}^{1,\,\,2}\,. \leqno (2.8)
$$
The kernel of $\,{\bf I}_o'\,$ consists of functions of the type (see \cite{Caffarelli-Gidas-Spruck}\,)
$$
V_{\lambda\,,\ \xi}(y) \ = \ \left({\lambda\over {\lambda^2\,\,+\,\,|\,y\,-\,\xi\,|^2}}\right)^{\!\!{{n\,-\,2}\over 2} } \ \ \ \mfor\ \ \ \ (\lambda\,, \ \xi) \  \in\ \R^+ \,\times \,\R^n\,, \leqno (2.9)
$$
which satisfies the equation
$$
\Delta V_{\lambda\,,\ \xi}\,(y) \ + \ n\,(n\,-\,2) [\, V_{\lambda\,,\ \xi}\,(y)\,]^{{n\,+\,2}\over {n\,-\,2}} \ = \ 0 \ \ \ \ \ {\mbox{in}} \  \ \ \ \R^n.  \leqno (2.10)
$$
We consider juxtaposition of two bubbles
$$
{\bf z}_\sigma \ = \ V_{\lambda_1\,,\ \xi_1} \ +  \ V_{\lambda_2\,,\ \xi_2}\ \   \mfor\ \ \
(\,\lambda_1\,, \  \lambda_2\,; \ \ \xi_1\,, \ \,\xi_2\,) \  \in\ (\,\R^+  \times  \R^+) \,\times\,(\R^n \,\times \,\R^n)\,. \leqno (2.11)
$$

\vspace*{0.1in}

{\bf \S\,2\,b.}\ \  {\it Unit and restrictions.} \ \ In the following we assume that
$$
{\bar C }^{-1}  \cdot \lambda_2 \ < \ \lambda_1 \ <\  {\bar C} \cdot \lambda_2\,, \ \ \ \ |\,\xi_1\,| \ <  \   {\bar c} \,\cdot \lambda \ \  \ \ \  {\mbox{and}} \  \ \ \ \ \ |\,\xi_2 \ - \ {\bf q}_{\,2}\,| \ <  \  {\bar c} \,\cdot \lambda\,. \leqno (2.12)
$$
Here $\,{\bar C }\,$ ($\,>\,1$)\, and $\,{\bar c }\,$ ($\,\approx \,0^+\,$)\, are positive constants (to be more precisely described  in {\bf \S\,5\,}\,)\,.\,
With (2.11)\,,\, we define
$$
\lambda \  = \   \sqrt{\lambda_1 \cdot \lambda_2  \,}\ \ \ \ {\mbox{and}} \ \ \ \
\bf {D}  \  = \ {\gamma\over \lambda} \ \left( = {{|\,{\bf q}_2\,| }\over {\sqrt{\lambda_1 \cdot \lambda_2  \,}  }}\right) \,. \leqno(2.13)
$$
These imply  \\[0.1in]
(2.14)

\vspace*{-0.45in}

\begin{eqnarray*}
& \ &  {1\over {\sqrt{{\bar C}\,}}} \cdot \lambda_{\,j} \ \le \ \lambda  \ \le \ \sqrt{{\bar C}\,} \cdot \lambda_{\,j} \ \ \ \ \ \ \ {\mbox{for}} \ \ \ j \ = \ 1\,, \ 2\,,\\[0.15in]
  \left[\ {\bf {D }} \ - \ 2 \,{\bar c } \, \right] \   \le \  {\bf d}  & := &  {{|\,\xi_1\ -\ \xi_2\,|}\over{\,\sqrt{\,\lambda_1 \,\cdot\,\lambda_2\,}\,}}\    \le \   \left[\ {\bf {D }} \ + \ 2\,{\bar c }\,  \right]\ \ \
 {\mbox{and}} \ \ \ \ \   \displaystyle{1\over { {\bf d} }} \  = \  {1\over {\ {\bf {D }}   }}  \cdot \left[\ 1 \ + \ O\,\left( \,{{{\bar c }}\over { \bf {D } }} \right) \, \right]\ .
 \end{eqnarray*}

 {\bf \S\,2\,c.}\ \   {\it Weak interaction.} \ \  We know that
 $$
 {\bf I}_o'\,(V_{\lambda\,,\,\xi}) \ \equiv \ 0\,, \ \ \ \ {\mbox{but}} \ \ \  {\bf I}_o'\,(V_{\lambda_1\,,\,\xi_1}\ + \ V_{\lambda_2\,,\,\xi_2}) \ \not\equiv \ 0\,.
 $$
In this section we investigate the ``interaction\," in more detail\,.\, From (2.8) and (2.11) we have \\[0.1in]
(2.15)
 $$
 {\bf I}_o'\,({\bf z}_\sigma)\,[\,h\,]\ = \ n\,(n\,-\,2)  \int_{\R^n}  \left\{ \left( [\, V_{\lambda_1\,,\ \xi_1} \,]^{{n\,+\,2}\over {n\,-\,2}} \ + \ [\, V_{\lambda_2\,,\ \xi_2} \,]^{{n\,+\,2}\over {n\,-\,2}} \right) \ - \ \left[\  V_{\lambda_1\,,\ \xi_1}  \ + \ V_{\lambda_2\,,\ \xi_2} \,\right]^{\!{{n\,+\,2}\over {n\,-\,2}}} \,\right\} \cdot h \ .
$$
for $\,h\ \in\ {\cal D}^{1,\,\,2}\,,\,$

\newpage

{\bf Lemma 2.16 (Weak Interaction Lemma)\,.} \ \   {\it Assume that \,$\,n\,\ge\,6\,$,\, with the notations and conditions in\,} (2.12) {\it and}\, (2.13)\,,\,
{\it there exists a positive constant \,$\,\bar{D}_1 \ > \ 1\,$ such that}
$$
{\bf D}\  = \ {\gamma\over \lambda } \ \ge\ \bar{D}_1\ , \leqno (2.17)\ \ \ \ \ {\it if}
$$

\vspace*{-0.25in}

$$
\!\!\!\!\!\!\!\!\!\!\!\!\!\!\!\!\!\!\!\!\!\!\!\!\!\!\!\!\!\!\Vert \, {\bf I}_o'\,(\,{\bf z}_\sigma\,)\,\Vert\ \le\ {\bar C}_1 \cdot {{\ln\,{\bf D}   }\over{\ {\bf D}^{ {{n\,+\,2}\over 2} \ }}}\ .\leqno (2.18)\ \ \ \ \ {\it then}
$$
{\it In\,} (2.17) {\it and\,} (2.18)\,, {\it the positive constants $\,\bar{D}_1\,$ and \,$\,{\bar C}_1\,$} {\it can be precisely determined by $\,{\bar C}\,,\,$ $\,{\bar c}\,$} [\,{\it appeared in\,} (2.12)\,] {\it and $\,n\,$,\, and} {\it they are independent on} \,$\,(\,\lambda_1\,,\ \lambda_2\,;\ \xi_1\,,\ \xi_2\,)\,$ {\it as long as}\, (2.12)   {\it is satisfied\,.}

\vspace*{0.2in}

\hspace*{0.5in}The proof can be seen from the proof of Lemma 2.1 in   \cite{III}\,, together with Lemma A.5 in the Appendix.

\vspace*{0.2in}

{\bf \S\,2\,b.}\ \ {\it Interaction terms.}\ \  In the following we describe the interaction between two bubbles via (2.15). We first observe that  in a small neighborhood of $\,\xi_1\,,\,$  $\,V_{\lambda_2\,,\ \xi_2}\,$ is small when compared to $\,V_{\lambda_1\,,\ \xi_1}$\,.\,  Precisely,   we let
$$
\rho_\mu \ =  \ \mu \cdot |\,\xi_1 \ - \ \xi_2\,| \ . \leqno (2.19)
$$
Here $\,\mu\,$ is a chosen small positive number so that
$$
\mu \ \to \ 0^+\ \ ({\mbox{slowly}}) \ \ {\mbox{and}} \ \ \ \ \mu^M \cdot {\bf D} \ \to \ \infty \ \ \ \ {\mbox{when}} \ \ \ \ {\bf D} \ \to \ \infty\,. \leqno (2.20)
$$
Here $\,M\,$ is a (fixed) large integer\,.\,  For most particular purpose one can take
$$
\mu \ = \ {1\over 2} \cdot {1\over { {\bf D}^{\,\epsilon} }} \ \ \ \ \ \ {\mbox{for}} \ \ \ \ {\bf D} \ \gg \ 1\,,
$$
where $\,\epsilon\ < \ 1\,$ is any fixed small positive number\,.\, Under the conditions in (2.12)\,,\, we have

\vspace*{-0.3in}

\begin{eqnarray*}
(2.21) \ \ \ \ \ \ \ & \ & \!\!\!\!\!\!\!\!\!\!\!\!\!\!\!\!V_{\lambda_2\,,\ \xi_2} (\,y)   \ = \ \left( {{\lambda_2}\over { \lambda^2_2 \ + \ |\,y\ - \ \xi_2\,|^2  }}\right)^{\!\!{{n\,-\,2}\over 2}  } \ =  \ \left(\ {{{1\over\lambda_1}}\over { \left( {{\lambda_2}\over {\lambda_1}}  \right)\ + \ {{|\,y\ - \ \xi_2\,|^2}\over
{\lambda_1 \cdot \lambda_2}}  }}\right)^{\!\!{{n\,-\,2}\over 2}  }  \\[0.15in]
& = &  {1\over {\lambda_1^{{n\,-\,2}\over 2} }} \cdot \left(\,{{{1}}
\over { \left( {{\lambda_2}\over {\lambda_1}}  \,\right)\ + \ {{|\,(\,y\ - \ \xi_1)\ + \ (\,\xi_1\ - \ \xi_2)\,|^2}\over
{\lambda_1 \cdot \,\lambda_2}}  }}\ \right)^{\!\!{{{n\,-\,2}\over 2}  }} \\[0.15in]
& = &  {1\over {\lambda_1^{{n\,-\,2}\over 2} }} \cdot \left[\ {{{1}}
\over { \left( {{\lambda_2}\over {\lambda_1}}  \right) \ + \ {{|\, \,\xi_1\ - \ \xi_2\,|^2}\over
{\lambda_1 \cdot \,\lambda_2}}   \ + \ {{|\,\,y\ - \ \xi_1\,|^2}\over
{\lambda_1 \cdot \,\lambda_2}}   \ + \ {{\,2\,(\,y\ - \ \xi_1)\,*\,(\,\xi_1\ - \ \xi_2)\, }\over
{\lambda_1 \cdot \,\lambda_2}}}}\ \right]^{{{n\,-\,2}\over 2}  }\\[0.1in]
& \ & \hspace*{1.5in} \{\  \uparrow \ = \ {\bf d}^2 \ [\,\gg \ 1\,; \ \ {\mbox{cf.}} \ \ (2.14)\,]\,,\ {\mbox{dominating \ \ term}}\,\}\\[0.01in]
& = &  {1\over {\lambda_1^{{n\,-\,2}\over 2} }} \cdot {1\over {{\bf d}^{n\,-\,2} }} \cdot \left[\ {{{1}}
\over { 1 \ + \ {1\over {{\bf d}^{2} }} \cdot \left( {{\lambda_2}\over {\lambda_1}}  \right) \ + \ {1\over {{\bf d}^{2} }} \cdot {{|\,\,y\ - \ \xi_1\,|^2}\over
{\lambda_1 \cdot \,\lambda_2}}   \ + \ {1\over {{\bf d}^{2} }} \cdot {{\,2\,(\,y\ - \ \xi_1)\,*\,(\,\xi_1\ - \ \xi_2)\, }\over
{\lambda_1 \cdot \,\lambda_2}}}}\ \right]^{{{n\,-\,2}\over 2}  }  \\[0.2in]
& =  & {1\over {\lambda_1^{{n\,-\,2}\over 2} }} \cdot {1\over {{\bf D}^{n\,-\,2} }} \,\cdot  \, [\,1 \ + \ o\,(1)\,] \ = \ {{ \lambda_2^{{n\,-\,2}\over 2 } }\over {\gamma^{n\,-\,2}   }} \,\cdot  \, [\,1 \ + \ o\,(1)\,]  \mfor \ \ y \ \in \ B_{\xi_1}\,(\rho_\mu)\ . \ \ \ \ \ \ \ \ \ \
\end{eqnarray*}
In the above we apply (2.13), (2.14) and (2.20)\,.\,
Moreover, $\,o\,(1)\ \to \ 0\,$ as $\,{\bf D} \ \to \ \infty\,$.\,
Compare with
$$
\, V_{\lambda_1\,,\ \xi_1} \ \approx \ {1\over {\lambda_1^{{n\,-\,2}\over 2} }} \cdot {1\over {{\bf D}^{n\,-\,2} }} \cdot  {1\over {\mu^{n\,-\,2} }}\ \ \  \ \ \ \  \ \   {\mbox{on}} \ \ \  \ \partial B_{\xi_1}\,(\rho_\mu)\,.
$$

\vspace*{-0.25in}

Let
$$
{\cal I} \ := \ \left( V_{\lambda_1\,,\ \xi_1}^{{n\,+\,2}\over{n\,-\,2}}\ +\ V_{\lambda_2\,,\ \xi_2}^{{n\,+\,2}\over{n\,-\,2}} \right) \ - \  (V_{\lambda_1\,,\ \xi_1} \ + \ V_{\lambda_2\,,\ \xi_2})^{{n\,+\,2}\over{n\,-\,2}} \ \ (\ < \ 0\,)\ . \leqno (2.22)
$$
By using the partition
$$
\R^n \ = \ B_{\xi_1} ({\rho_\mu}) \ \cup \ B_{\xi_2} ({\rho_\mu}) \ \cup \ \{\,\R^n \,\setminus\,  [\,B_{\xi_1} ({\rho_\mu}) \,\cup\,B_{\xi_2} ({\rho_\mu})\,]\, \}\ ,
$$
and the inequality ($\,A\,$ and $\,B\,$ are positive numbers\,)
$$
\left( A^{{n\,+\,2}\over {n \,-\,2}} \, + \, B^{{n\,+\,2}\over {n \,-\,2}} \right) \ - \ (A \, + \, B)^{{n\,+\,2}\over {n \,-\,2}}
\ = \  -\,{{n\,+\,2}\over {n \,-\,2}}\,\cdot\,
A^{4\over {n \,-\,2}} \cdot B \ + \ O\,(1) \cdot B^{{n\,+\,2}\over {n \,-\,2}} \ \ \ \  \   \mbox{for}  \ \ \ \ {B\over A} \ \
{\mbox{small}}\,,
$$
we obtain (\,see {\bf \S\,A.1} in the {\bf e}\,-\,Appendix for more detail)
\begin{eqnarray*}
(2.23) \ \ \ \ \ \ \ \ \ \ \ \ {\bf I}_o'({\bf z}_\sigma)\,[\,\partial_{\lambda_1} V_{\lambda_1\,,\ \xi_1}\,]   &   = &  n\,(n\,-\,2) \int_{\R^n} {\cal I} \,\cdot\, [\,\partial_{\lambda_1} V_{\lambda_1\,,\ \xi_1}\,]\\[0.15in]
&  = &    \,-\ {\bar C}^+_1 (n)  \cdot {1\over {\lambda_1}} \cdot {{\lambda_1^{{n\,-\ 2}\over 2}  \cdot \lambda_2^{{n\,-\ 2}\over 2}}\over {\gamma^{n\,-\,2} }}
 \cdot [\,1\ + \ o\,(1)\,] \ . \ \ \ \ \ \ \ \ \ \ \ \ \ \ \ \ \ \ \ \ \ \ \ \
 \end{eqnarray*}

 \vspace*{-0.25in}

Here
$$\,{\bar C}^+_1 (n) \, = \,  n\,(n\,-\,2) \cdot  \omega_n  \cdot    {{n \,-\  2}\over {2\,n}}\,.\,\leqno (2.24)$$
Likewise\,,\, we extract  the leading term in $\xi$\,-\,derivatives  (\,refer to {\bf \S\,A.1} in the {\bf e}\,-\,Appendix for more of the calculations)\,:
\begin{eqnarray*}
(2.25) \ \ \ \ \ \  \ \ \ \ \ \  {\bf I}_o'({\bf z}_\sigma)\,[\,\partial_{\xi_{1_1}} V_{\lambda_1\,,\ \xi_1}\,]  & = &   n\,(n\,-\,2)  \int_{\R^n} {\cal I} \cdot [\,\partial_{\xi_{1_1}} \,V_{\lambda_1\,,\ \xi_1}\,]\\[0.15in]
 &= &  {\bar C}^+_2(n) \cdot  {1\over {\lambda_1}} \cdot  {1\over { \lambda_2}} \cdot (\xi_{1_1} - \,\xi_{2_1} ) \cdot  {1\over {{\bf D}^{n} }}
   \cdot [\,1\ + \ o\,(1)\,]\ .\ \ \ \ \ \ \ \ \ \ \ \
 \end{eqnarray*}
 Here $\,{\bar C}^+_2 (n)\,$  is a positive constant depending on $\,n\,$ only\,.\,

 \newpage

\hspace*{0.5in}Similarly for the expressions on $\,\partial_{\lambda_2}  V_{\lambda_2\,,\ \xi_2}\,$ [\,with the same constant  \,$\,{\bar C}^+_1 (n)\,$] and $\,\partial_{\xi_{2_1}} \, V_{\lambda_2\,,\ \xi_2}\,$ [\,with the same constant  \,$\,{\bar C}^+_2 (n)\,$]\,.\,
Here conditions (2.12) and (2.17) apply\,,\, and $\,o\,(1) \ \to \ 0\,$ when $\,{\bf D}  \ \to \ \infty\,.\,$ We present the rather standard calculations in  {\bf \S\,A.1} in the {\bf e}\,-\,Appendix, paying special attention of the sign ($+\,/\,-$)\,.\, Compare also with  Lemma B.2 and Lemma B.4 in \cite{Yan+} and formulas  2.119 and 2.206 in \cite{Bahri}\,.\,

 \vspace*{0.5in}

{\Large{\bf \S\,3.}\ \ {\bf Second order property \,-\, solving the equation in the}}\\[0.1in]
\hspace*{0.35in} {\Large\bf{ perpendicular directions.}}\\[0.2in]
As is often the case in mathematics,  simplicity is linked with orthogonality, such as the Lagrange multiplier method. Likewise,  the Lyapunov\,-\,Schmidt reduction method consists of solving the equation

\vspace*{-0.35in}

$$
{\bf I}' \,(u)\ \equiv\ 0 \ \ \mfor \ \ {\mbox{an \ \ unknown}} \ \ u \ \in \ {\cal D}^{1,\,\,2} \leqno (3.1)
$$
in two steps, first in the `perpendicular' direction, and then in the `horizontal' direction.  Here we first consider the  `perpendicular' direction. Given $\,{\bf z}_\sigma\,$ as in (2.11), let
$$
\!\!\!\!\!\!\!\!\!\!\!\!\!\!\!\!\!\!\perp_{\,\sigma}  \ \,= \ \bigg\{\  h\ \in \ {\cal D}^{1,\,\,2} \ \  | \ \
 \langle \,h\,,\, \partial_{\lambda_j} V_{\lambda_j\,, \, \xi_j} \,\rangle_\btd \ = \
 \langle \,h\,,\, \partial_{\xi_{j_{|_k}}} V_{\lambda_j\,, \, \xi_j}\,\rangle_\btd \ = \ 0\ \leqno(3.2)
$$

\vspace*{-0.35in}

$$
\hspace*{4in}{\mbox{for}} \ \ \   j \ = \ 1\,,\ 2\ ; \ \ \ k \ = \ 1\,,\ 2\,,\ \cdots\ ,\ n\,\bigg\}\ ,
$$

\vspace*{-0.3in}

and

\vspace*{-0.3in}

$$
{\cal P}_\sigma \,: \ {\cal D}^{1,\,\,2} \ \to \ \ \perp_{\,\sigma}\leqno(3.3)
$$
be the orthogonal projection\,. Here $\,\displaystyle \xi_{j}\,=\,\left(\,\xi_{j_{\,|_1}}\,,\ \cdots\ ,\ \xi_{j_{\,|_n}}\,\right)\,\in\,\R^n\ $ for $\,j\ = \ 1\,, \ 2\,.\,$  Fixed a $\,{\bf z}_\sigma\,$, to solve the {\it auxiliary equation} (\,`perpendicular' direction) is to find an unknown $\,w_{{\bf z}_\sigma}\,\in\ \,\perp_{\,\sigma}\,$ in the equation

\vspace*{-0.3in}

$$
\ \ \   {\cal P}_\sigma \!\circ\,{\bf I}'_\varepsilon \,\,(\, {\bf z}_\sigma\, + \, w_{{\bf z}_\sigma}) \ = \ 0\ \ \ \ \ \  (\,w_{{\bf z}_\sigma}\,\in\ \,\perp_{\,\sigma}\ \ {\mbox{is \   called \   informally \   }} ``\,\perp\,-\,{\mbox{solution}}\,"\,)\ .\leqno (3.4)
$$
In order to apply the implicit function theorem to solve equation (3.4), we proceed to the second Fr\'echet derivative of $\,{\bf I}_o\,$, in particular, the ``\,diagonal element\,"\,:
$$
(\,{\bf I}_o''\,({\bf z}_\sigma)\, [\,f\,] \,f\,)\ = \ \int_{\R^n} \!\left[\ \langle\,\btd\,f\,,\,\btd\,f\,\rangle\  -\  n\,(n\,+\,2)\,\left(V_{\lambda_1\,, \ \xi_1}\, +\, V_{\lambda_2\,,\ \xi_2}\right)^{{4}\over {n\,-\,2}}  \cdot f^2\ \right] \leqno (3.5)
$$
for $\,f \,\in\,{\cal D}^{1,\,\,2}\,$. For a proof of the following lemma, see \cite{Bahri}\,,\, and Lemma 2.5 in \cite{III}.

\vspace*{0.2in}

{\bf Lemma 3.6\ \ (Non\,-\,degeneracy Lemma)\,.}\ \ \ {\it Assume that \,$\,n\,\ge\,6\,$. Under the conditions in } (2.12)\,,\, {\it there exists a positive constant \,$\,\bar{D}_2\,$ such that }
$$
\ \ \ \ \ \ \ \ \ \ \ \ \ {\bf D} \ \ge \ \bar{D}_2\ , \leqno (3.7)\ \ \ \ \ \ \ \ {\it if}
 $$

\vspace*{-0.25in}

$$
|\, (\,{\bf I}_o'' \,({\bf z}_\sigma)\ [\,f\,]\ f\,)\,|\ \,\ge\,\ {\bar{c}}_{\,\sigma}^{\,2}\,\,\Vert\,f\,\Vert_\btd^2\ \ \ \ \ \ \ \ \ \ \ {\it{for \ \ all}} \ \ \  f \ \in\ \,\perp_{\,\sigma}\,. \leqno (3.8)\ \ \ \ \ \ \ \ {\it then}
$$

\newpage

{\it Here the constant \,$\,{\bar{c}}_{\,\sigma}\,$ is independent on \,$\,(\,\lambda_1\,,\, \lambda_2\,;\ \xi_1\,,\ \xi_2\,)\,$ as long as\,} (2.12) {\it and the condition in\,} (3.7) {\it are fulfilled\,.}\\[0.15in]
\hspace*{0.5in}As $\,\perp_{\,\sigma}\ \subset\ {\cal D}^{1,\,\,2}\,$ is itself a complete Hilbert space, we consider the restriction
$$
(\,{\bf I}_o''\,({\bf z}_\sigma)\,[\,f\,]\,\,h\,)\ = \ \int_{\R^n} \!\left[\,\,\langle\,\btd\,f\,,\,\btd\,h\,\rangle\ -\ n\,(n\,+\,2)\,{\bf z}_\sigma^{{4}\over {n\,-\,2}}\cdot f\cdot h\,\,\right]\ \ \mfor f\,,\ h\ \in\ \,\perp_{\,\sigma}\,. \leqno (3.9)
$$
Via the Riesz Representation Theorem, we obtain a linear map
\begin{eqnarray*}
\ \ \ \ \ \ \ \ \ {\bf I}_o''\,({\bf z}_\sigma ) \,:\ \ \perp_{\,\sigma}  & \to &  \perp_{\,\sigma}\ \ \ \ \ \ \ \ \ \ \ \ \ \ \ \ \ \ \ \ \ \ \ \ \ \ \ \ (\,h\,\in\ \perp_{\,\sigma}\!\setminus\,\{ 0\}\,,\ \ \  h\,\perp\,{\mbox{Ker}}\, {\bf I}_o''\,({\bf z}_\sigma )\,[f]\,) \ \ \ \ \ \ \ \ \ \ \ \ \ \ \ \ \ \ \ \ \ \ \ \ \ \ \\[0.15in]
 (3.10)\ \ \ \ \ \ \ \ \ \ \ \ \ \ \ \ \ \ \  \ \ \ \ \ f  & \mapsto  &  \left\{\,{{ ({\bf I}_o''\,({\bf z}_\sigma )\,[f]\,\,h)}\over {\Vert\,h\,\Vert_{\,\btd}^{\,2}}}\,\right\}\cdot h\ . \ \
\end{eqnarray*}
The following result can be seen as a direct consequence of Lemma 3.6. We refer to \cite{Progress-Book}\,,\, or  \cite{III}\,.

\vspace*{0.17in}

{\bf Lemma 3.11.}\ \ {\it Under the conditions in Lemma\,} 3.6\,,\, {\it the map  given in\,} (3.10) {\it is an isomorphism with uniformly bounded inverse\,.  }\\[0.15in]
\hspace*{0.5in}With the help of Lemma 3.12 below (a proof can be found in the Appendix), we now show that when the two bubbles are concentrated around the two critical points, one can solve the ``perpendicular\," direction, just like in the perturbation case \cite{Progress-Book}.

\vspace*{0.2in}

{\bf Lemma 3.12.} \ \ {\it Assume that $\,n \ \ge \ 6\,,\,$  $\,|\,H\,| \ \le \  {\bar C}_b\,$,\, and under the conditions in} (1.8)\,,\, (1.9)\,,\, (1.10)\,,\, (1.11)  {\it and\,} (2.12)\,.\,
{\it Given a number $\,m\,$ so that  $\,m \cdot \ell \ > \ 2\,,\,$ there is a positive number \,$\bar{D}_3\,$ such that\,}
$$\,\ \ \ \ \ \ \ {\bf  D}\ \ge\ \bar{D}_3\,,\, \hspace*{1.5in} \leqno{\it{if}}$$
$$
{\it{then}} \ \ \ \ \  \ \ \ \ \  \ \ \ \ \    \int_{\R^n} |\,H|^{\,m} \,\cdot\, {\bf z}_\sigma^{{2n}\over {n\,-\,2}} \ \leq\ \left\{
\begin{array}{l@{\ \ \ \ }l}
\displaystyle{C_4 \cdot \lambda^{m\,\ell}\,   \ \ \ \ \ \,+ \ {\bf R}  }& \displaystyle{\it{if}\ \ } \ \ m \cdot \ell \ <\ n\ ,\\[0.15in]
\displaystyle{C_5 \cdot\lambda^{\,n\,-\,o\,(1)}\,    \, + \ {\bf R}  } & \displaystyle{\it {if}\ \ } \ \ m \cdot \ell \ = \ n\ , \ \ \ \ \  \ \ \ \ \   \ \ \ \ \  \ \ \ \ \    \\[0.15in]
\displaystyle{C_6 \cdot \lambda^n\,  \ \ \ \ \ \ \ \,+ \ {\bf R}  }  & \displaystyle{\it {if}\ \ } \ \ m \cdot \ell \  > \ n\ .\\[0.15in]
\end{array}
\right.
$$
{\it Here}
$$
  {\bf R}  \ = \   \lambda^2 \cdot  {{C_7}\over { {\bf D}^{\,n\,-\ 2} }}   \ + \    {{C_8}\over { {\bf D}^{\,n\,-\ o\ (1)} }} \ \,,
$$
{\it and $\,o\,(1) \,\to\,0^+$\, as \,$\lambda\,\to\,0\,.\,$}
{\it Moreover, the constants $\,{\bar D}_3$\,,\, $\,C_4\,,\,$ $\,C_5\,,$ $\,C_6\,,\,$ $\,C_7\,$ and $\,C_8\,$  are independent on \,$\,(\,\lambda_1\,,\ \lambda_2\,;\ \xi_1\,,\ \xi_2\,)\,$ as long as\,} (2.12) {\it and the conditions in\,} ${\mbox{(1.9)}}_{i}\,, \ {\mbox{(1.9)}}_{ii}\,$ $\,{\mbox{(1.9)}}_{iii}\,$ (1.10) {\it and}\, (1.11)\, {\it are fulfilled\,.}\\[0.15in]

\newpage

{\bf Theorem 3.13 (\,Existence of small $\perp\!-\,$solution.\,)} \ \
{\it Assume that $\,n \ \ge \ 6\,,\,$ $\ell \  \ge \ 2\,,\,$   and the conditions in Lemma} 3.12\,.\,
{\it Then there exist   positive numbers $\,{\bar D}_4\,$} ({\it relatively }``{\it \,large}\,") {\it and} $\, \tilde \epsilon\,$ (``{\it \,small}\,") {\it such that for each $\,{\bf z}_\sigma\,$ with }

 \vspace*{-0.15in}

$$
{\bf D} \ \ge \ {\bar D}_4\,,\,
$$

 \vspace*{-0.15in}

{\it   the auxiliary equation}
$$\ \ \ \ \ \  \
{\cal P}_\sigma  \circ {\bf I}'\,\, (\,{\bf z}_\sigma \,+\,w) \ = \ 0 \ \ \ \ \ \  \ \ \ \ \ \ \  \ \ \ \ \ \ \  \ (\,w \ \in \ \perp_\sigma\,) \leqno (3.14)
$$
{\it has a unique   ``\,{\it small}\ "\,   solution $\,w_{{\bf z}_\sigma }\,\in\  \perp_{\,\sigma}\,,\,$ precisely,
 }
$$
 \Vert \, w_{{\bf z}_\sigma }\Vert_\btd \ \le \ \tilde \epsilon\ .\leqno (3.15)
 $$

 \vspace*{-0.1in}

{\it Moreover, one can take }
$$\,\tilde \epsilon \ \to \ 0^+ \ \ \ \ {\mbox{as}} \ \ \ \ {\bf D} \ \to \ \infty\,.$$
{\it The constant  $\,{\bar D}_4\,$ is independent on \,$\,(\,\lambda_1\,,\ \lambda_2\,;\ \xi_1\,,\ \xi_2\,)\,$ as long as\,} (2.12)\, {\it and the conditions in\,} ${\mbox{(1.9)}}_{i}\,, \ {\mbox{(1.9)}}_{ii}\,$ $\,{\mbox{(1.9)}}_{iii}\,$ (1.10) {\it and}\, (1.11)\, {\it are fulfilled\,.}\, {\it In addition, $\,w_{{\bf z}_\sigma}$
depends on the parameters $\,(\lambda_1\,,\,\lambda_2\,;\ \xi_1\,,\,\xi_2\, )\,$ of \ ${\bf z}_\sigma$ $\left( \ = \ V_{\lambda_1\,,\ \xi_1}\ +\ V_{\lambda_2\,,\ \xi_2}\right)$\, in a $\,C^1\,$ manner\,.
}

 \vspace*{0.2in}

{\bf Proof.} \ \ From (2.5)\,--\,(2.7)\,,\, we have
 \begin{eqnarray*}
{\bf I}'\,\, ({\bf z}_\sigma \,+\,w\,) [\,\bullet\,] & = &  \int_{\R^n} \left[\,\langle\, \btd \,({{\bf z}_\sigma}\,+\, w)\,,  \ \btd\, [\bullet]\,\rangle
 \ - \ n \, (n\,-\,2) ({{\bf z}_\sigma}+w)^{{n\,+\,2}\over {n\,-\,2}}\cdot [\bullet]\,\right]\\[0.2in]
& \ & \ \ \ \ \  \ \ \     - \    \int_{\R^n}
({\tilde c}_n \cdot H)\,({\bf z}_\sigma \,  +\, w\,)_+^{{n\,+\,2}\over {n\,-\,2}} \cdot [\,\bullet\,]  \ \ \mfor \  \bullet \ \in \ {\cal D}^{1\,,\,2}  \ \ \ \   (\,w \ \in \ \perp_\sigma\,)  \ .
\end{eqnarray*}

\vspace*{-0.1in}

Write\\[0.1in]
(3.16)

\vspace*{-0.3in}

$$
 {\cal P}_\sigma \circ {\bf I}'\,\, ({\bf z}_\sigma \,+\,w\,) \ =  \ {\cal P}_\sigma \circ {\bf I}'_o\, ({\bf z}_\sigma\,) \,+\, {\cal P}_\sigma \circ  {\bf I}''_o\, ({\bf z}_\sigma\,)\,[\,w]  \,+\,  {\cal P}_\sigma \circ  \, {\bf G}'\, ({\bf z}_\sigma \,+\,w\,) \,+\, {\cal P}_\sigma \circ  R_{{\bf z}_\sigma }(w)
$$
\hspace*{5.3in}for $\,\Vert w \Vert_\btd\ $ small\,,\,\\[0.1in]
where
 \begin{eqnarray*}
   {\bf I}'_o\,({{\bf z}_\sigma}) \,[\bullet] & = &  \int_{\R^n} \left[\,\langle\, \btd \,({{\bf z}_\sigma})\,,\
   \btd\, [\bullet]\,\rangle  \ - \ n \, (n\,-\,2) ({{\bf z}_\sigma})^{{n\,+\,2}\over {n\,-\,2}}\cdot [\bullet]\,\right]\,,\\[0.1in]
 (\,{\bf I}''_o\,({{\bf z}_\sigma})\,[w]\,\bullet) & = &  \int_{\R^n} \left[\,\langle \btd \,w\,, \  \btd\, [\bullet]\, \rangle  \ - \ n \, (n\,+\,2)\, {{\bf z}_\sigma}^{4\over {n\,-\,2}} \,w \cdot [\bullet]\,\right]\,,\ \ \ \  \ \  \ \ \ \\[0.1in]
 {\bf G}'\,({\bf z}_\sigma \,+\,w\,)\,[\,\bullet\,]  &  = &  - \    \int_{\R^n}
({\tilde c}_n \cdot H)\,({\bf z}_\sigma \,  +\, w\,)_+^{{n\,+\,2}\over {n\,-\,2}} \cdot [\,\bullet\,]\ ,
\end{eqnarray*}

\newpage

and \\[0.05in]
(\,3.17)
 \begin{eqnarray*}
 R_{{\bf z}_\sigma }(w)\,[\bullet] \!\! & = & \!\!\int_{\R^n} \left[\  n \, (n \,-\,2) \, {{\bf z}}^{{n\,+\,2}\over {n \,-\,2}}_\sigma \cdot [\bullet]  \, -\, n \, (n \,-\,2) ({{\bf z}_\sigma}+ w)_+^{{n\,+\,2}\over {n \,-\,2}} \cdot [\bullet] \,\,+\,\, n\,  (n\,+\,2)\, {{\bf z}_\sigma}^{4\over {n \,-\,2}}\, w \cdot [\bullet]\,\right] \\[0.15in]
& = & - \,n \, (n \,-\,2) \int_{\R^n} \left[     ({{\bf z}_\sigma} \,+\, w)_+^{{n\,+\,2}\over {n \,-\,2}}
\, -\, {{\bf z}}^{{n\,+\,2}\over {n \,-\,2}}_\sigma \,  - \, {{n\,+\,2}\over {n \,-\,2}} \
{{\bf z}_\sigma}^{4\over {n \,-\,2}}\, w \,\right]\cdot[\bullet]\ .
\end{eqnarray*}
In order to solve  equation (3.14)\,,\, that is
 \begin{eqnarray*}
{\cal P}_\sigma \circ  {\bf I}''_o\, ({\bf z}_\sigma\,)\,[\,w]  & = & - \,{\cal P}_\sigma \circ {\bf I}'_o\, ({\bf z}_\sigma\,)    \ -\   {\cal P}_\sigma \circ  {\bf G}'\,({\bf z}_\sigma \,+\,w\,) \,-\, {\cal P}_\sigma \circ  R_{{\bf z}_\sigma }(w)\,,\\[0.07in]
{\mbox{``\,interaction\  term\,"}} \!\!\!\!\!\!\!\!\!\!\!\!\!\!\!\!\!\!\!\!& \ & \ \ \ \ \ \ \uparrow \ \ \ \ \ \ \ \ \ \ \ \ \ \ \ \ \ \ \ \ \ \uparrow \ \ {\mbox{depending  \ \ on\ \ }} w \ \ \uparrow
\end{eqnarray*}
we first seek a solution $\,w\,$ to
$$
{\cal P}_\sigma \circ  {\bf I}''_o\, ({\bf z}_\sigma\,)\,[\,w]  \ = \ - {\cal P}_\sigma \circ {\bf I}'_o\, ({\bf z}_\sigma\,)\,.
$$
({\Large{\bf $\star$}}) \  {\it The interaction term\,.\,} \ \  From the Weak Interaction Lemma 2.16\,,\,
$$
\Vert \,{\bf I}'_o\, ({\bf z}_\sigma\,)\,\Vert \ \le \ C_1 \cdot  {1\over { {\bf D}^{{{n\,+\,2}\over 2} \ -\ o\,(1\,)} }}\ .\leqno (3.18)
$$
Using Lemma 3.11\,,\, we can find $\,{\bar w}_1 \ \in \ \perp_{\,\sigma}$\, such that
$$
{\cal P}_\sigma  \circ  {\bf I}''_o\, ({\bf z}_\sigma\,)\,[\,{\bar w}_1 ] \ = -\, {\cal P}_\sigma
\circ  {\bf I}'_o\, ({\bf z}_\sigma\,)\,,\,
$$

\vspace*{-0.25in}

and
$$
 \Vert \,{\bar w}_1  \Vert_\btd \  \le \ C_2 \cdot {1\over { {\bf D}^{ \,{{n\,+\,2}\over 2}  \,-\ o\,(1)}}} \ . \leqno (3.19)
$$
In (3.18) and (3.19)\,,\,  $\,o\,(1) \ \to \ 0^+\,$ as $\,{\bf D} \ \to \ \infty\,.\,$

\vspace*{0.2in}

({\Large{\bf $\star$}}) \   {\it Fixed point.} \ \ Next [\,because of the linearity of $\,{\cal P}_\sigma \circ {\bf I}''_o\, ({\bf z}_\sigma\,)\,$], we intend to find $\,w_2\, \in \ \perp_{\,\sigma}\,$ such that
$$
{\cal P}_\sigma  \circ \,{\bf I}''_o\, ({\bf z}_\sigma\,)\,[\,w_2] \ = \ -\, \left\{ \,[\,\,\,{\cal P}_\sigma  \circ \, {\bf G}']\,({\bf z}_\sigma \,+\,{\bar w}_1 \,+\,w_2\,) \,+\, {\cal P}_\sigma  \circ  R_{{\bf z}_\sigma}({\bar w}_1\,+\,w_2)\right\}\,. \leqno (3.20)
$$
[\,Here $\,{\bar w}_1\,$ appears in (3.19)\,.\,]
That is, we seek a fixed point to
$$
\bullet   = {\bf T}\,(\bullet) \ := \ -\, (\, {\cal P}_\sigma \circ {\bf I}''_o\, ({\bf z}_\sigma\,)\,)^{-1} \, \left\{ \,[\,\,\,{\cal P}_\sigma \circ   {\bf G}'\,]\,({\bf z}_\sigma \,+\,{\bar w}_1  \ +\  \bullet\,) \,+\, {\cal P}_\sigma \circ  R_{{\bf z}_\sigma}({\bar w}_1 \,+\,\bullet)\,\right\}\,. \leqno (3.21)
$$

From (3.16) we have
$$
   {\bf G}'\,(\,{\bf z}_\sigma \,+\,{\bar w}_1 \,+\,w_2\,)\,[\,h] \ = \   -\, \int_{\R^n} ({\tilde c}_n \cdot H)\,({\bf z}_\sigma \,+\,{\bar w}_1 \,+\, w_2\,)_+^{{n\,+\,2}\over {n\,-\,2}} \cdot h  \mfor \  h \ \in \ {\cal D}^{1\,,\,2}. \leqno (3.22)
$$
Also\,,\, Lemma 3.12 implies that
$$
\int_{\R^n} |\,H|^{\,{{ 2n}\over {n\,+\,2  }} } \cdot {\bf z}_\sigma^{{2\,n}\over { n\,-\,2 }} \ = \ O \left(\lambda^{\,*  }\,\right) \ + \ \lambda^2 \cdot \,O\,\left({1\over { {\bf D}
^{\,n\,-\ 2 } }}\right)
\ + \  O\,\left({1\over { {\bf D}
^{\,n\,-\ o\,(1) } }}\right)\,.\leqno (3.23)
$$
$$
{\mbox{where}} \ \ \ \ \  \ \ \ \ \  \ \ \ \ \   \ \  \ \ \ \ \  *  \ = \ \left\{
\begin{array}{l@{\ \ \ \ }l}
\displaystyle{{{ 2n}\over {n\,+\,2  }} \cdot \ell  }& \displaystyle{\mbox{if}\ \ } \ \  {{ 2n}\over {n\,+\,2  }} \cdot \ell \ <  \ n \ ,\\[0.15in]
\displaystyle{n \ - \ o\,(1)  } & \displaystyle{\mbox{if}\ \ } \ \ {{ 2n}\over {n\,+\,2  }} \cdot \ell \  \ge  \ n\ . \ \ \ \ \  \ \ \ \ \   \ \ \ \ \  \ \ \ \ \    \ \ \ \ \  \ \ \ \ \   \ \ \ \ \  \ \ \ \ \
\end{array}
\right.
$$
Here $\,o\,(1) \,\to\,0^+$\  as \ $\lambda\,\to\,0\,$ [\,condition (2.12) applies\,]\,.\,  Using the Minkowski inequality and Sobolev inequality we obtain
\begin{eqnarray*}
 (3.24)  \ \ \ \ \ &\  &\!\! \!\! \! \bigg\vert \, {\bf G}'({\bf z}_\sigma \,+\, {\bar w}_1 \,+\,w_2\,)\,[\,h\,]\,  \bigg\vert\\[0.15in]
   &\le &    |\,{\tilde c}_n\,| \cdot \bigg\vert \ \int_{\R^n} H\,({\bf z}_\sigma \,+\,{\bar w}_1 \,+\, w_2\,)_+^{{n\,+\,2}\over {n\,-\,2}} \cdot h \ \bigg\vert\\[0.15in]
& \le & C\cdot \left\{ \ \int_{\R^n} |\,H|\,{\bf z}_\sigma^{{n\,+\,2}\over {n\,-\,2}} \cdot |\,h| \ + \ \int_{\R^n} |\,H|\, |\,{\bar w}_1\,|^{{n\,+\,2}\over {n\,-\,2}} \cdot |\,h| \ + \ \int_{\R^n} |\,H|\,|\,w_2\,|^{{n\,+\,2}\over {n\,-\,2}} \cdot |\,h\,|\, \right\} \ \ \ \ \ \ \ \ \ \ \ \ \ \\[0.15in]
& \le & C\,(n) \cdot \Bigg\{ \, \left[ \ O \left(\lambda^{\,*}\,\right) \ + \ \lambda^2 \cdot \,O\,\left({1\over { {\bf D}
^{\,n\,-\ 2 } }}\right)
\ + \ O\,\left({1\over { {\bf D}^{n\,-\,2} }}\right) \right]^{{n\,+\,2}\over {2n}} \\[0.15in]
& \ & \ \ \ \ \ \  \ \ \ \ \ \ \ \ \ \ \ \ \ \ \ \ \ \  \ \ \ \ \ \ \ \ \ \ \ \  \ \ \ \ \ +  \
\left(\, \Vert\,{\bar w}_1\Vert_\btd \right)^{\,{{n \,+\,2}\over {n\,-\,2}}}
\ + \ \left(\Vert\,w_2\,\Vert_\btd \right)^{\,{{n \,+\,2}\over {n\,-\,2}}} \ \Bigg\} \cdot \Vert\,h\,\Vert_\btd\ .
\end{eqnarray*}

\vspace*{-0.15in}

Likewise \{\,from (3.17)\,,\, see also \cite{III}\,\}\,,\\[0.1in]
(3.25)
$$
  |\, R_{{\bf z}_\sigma} \,({\bar w}_1\,+\,w_2\,)\,(h)\,  | \ \ \le \  C\,(n) \cdot \left[\ {1\over {\, {\bf D}^{ {{n\,+\,2}\over 2} \ - \ o\,(1)\,  }  }} \ + \  \left\{\, \Vert\,{\bar w}_1\Vert_\btd \right\}^{{n \,+\,2}\over {n\,-\,2}} \ + \ \left\{\, \Vert\,w_2\,\Vert_\btd \right\}^{{n \,+\,2}\over {n\,-\,2}} \ \right] \cdot \Vert\,h\,\Vert_\btd\ .
$$
Hence by choosing $\,\lambda\,$ ( $ = \, \sqrt{\lambda_1 \cdot \lambda_2\,}$\ )\, and $\,c\,$ to be small, we obtain
$$
  \Vert\, w_2 \,\Vert_\btd \ \le \  c \ \ \ \Longrightarrow \ \ \ \Vert \,{\bf T}\,(h ) \,\Vert \ \le \ c\,.
$$
Note that, via (3.19),  (3.24) and  (3.25)\,,\,  $\ c\ \to \ 0^+$\, as $\,\lambda \ \to \ 0^+\,$  and $\, {\bf D} \ \to  \ \infty$\,.

\vspace*{0.2in}

({\Large{\bf $\star$}}) \  {\it Contraction map.} \ \ Consider [\,${\bar w}_1\,,\,$ first appeared in (3.19)\,,\, and $w_1$ below are to be distinguished\,]
$$
\Vert \, {\bf T}\,(w_1)\ - \ {\bf T}\,(w_2) \,\Vert\,, \ \ \ \ \ \ {\mbox{where}} \ \ \ \Vert\, w_1  \Vert_\btd \ \le \  c \ \ \ {\mbox{and}} \ \ \ \Vert\, w_2  \Vert_\btd \ \le \  c\,.
$$
Using the inequality (\,cf. for example \cite{Progress-Book}\,,\, $\,n \ \ge \ 6\,$)
$$
\bigg\vert \ a^{{n \,+\,2}\over {n\,-\,2}} \ - \ b^{{n \,+\,2}\over {n\,-\,2}}\,\bigg\vert \ \le \ C\,(n)
\cdot |\,a\,-\,b\,| \cdot (a \ + \ b)^{4\over {\,n\,-\ 2\,}}\mfor a \ > \ 0\ \ \& \ \ b \ > \ 0\,, \leqno (3.26)
$$
we have\\[0.1in]
(3.27)
\begin{eqnarray*}
& \ &   \bigg\vert \ [\,  {\bf G}'({\bf z}_\sigma \,+\,{\bar w}_1 \,+\,w_1\,)\ - \  {\bf G}'({\bf z}_\sigma \,+\,{\bar w}_1 \,+\,w_2\,)]  \,[h] \,\bigg\vert \\[0.15in]
& = &    {\tilde c}_n \cdot \bigg\vert \,     \int_{\R^n} |\,H\,| \ \left[\ ({\bf z}_\sigma \,+\,{\bar w}_1 \,+\,w_1\,)_+^{{n\,+\,2}\over {n \,-\, 2}} \ - \ ({\bf z}_\sigma \,+\,{\bar w}_1 \,+\,w_2\,)_+^{{n\,+\,2}\over {n \,-\, 2}} \right] \cdot h\, \bigg\vert\\[0.15in]
& \le & C\cdot\int_{\R^n}|\,H\,| \cdot |\, w_1 - w_2\,| \cdot \left({\bf z}_\sigma \,+\, |\,{\bar w}_1 |\,+\,|\,w_1 |\,+\,|\,w_2 | \ \right)^{{4\over {n\, -\, 2}} } \cdot |\,h|\\[0.15in]
& \le & C\, (n) \cdot \left[ \  \left[ \ O \left(\lambda^{\,n \,-\,o\,(1)}\,\right)
\ + \ O\,\left({1\over { {\bf D}^{n\,-\,2} }}\right) \right]^{{2}\over { n}}
\ + \ \Vert \,{\tilde w}_1\,\Vert_\btd^{ {4\over {n\,-\,2}} }
\ + \ \Vert \,{w}_1\,\Vert_\btd^{ {4\over {n\,-\,2}} }  \ + \ \Vert \,{w}_2\,\Vert_\btd^{ {4\over {n\,-\,2}} }
\,\right]\,*\\[0.2in]
& \ &  \ \ \ \  \ \ \  \ \      *\, \Vert \, h\,\Vert_\btd \cdot \Vert \, w_1 -
w_2\,\Vert_\btd \ \ \ \  \     (\,{\mbox{taking \   the  \   condition \   }}\ell \ \ge \ 2 \ \ {\mbox{into \  consideration}}\,)\ .\\
\end{eqnarray*}

\vspace*{-0.2in}

Here we use Minkowski inequality and Sobolev inequality again. \bk
Finally, for the remainder\,,\, we make use of the inequality
\begin{eqnarray*}
& \ & \bigg\vert \ \left[\ ({\bf z}_\sigma \,+\,{\bar w}_1 \,+\,w_1\,)_+^{{n\,+\,2}\over {n \,-\, 2}} \ - \ ({\bf z}_\sigma \,+\,{\bar w}_1 \,+\,w_2\,)_+^{{n\,+\,2}\over {n \,-\, 2}}\ \right] \ - \ {{n\,+\,2}\over {n\,-\,2}} \cdot  {\bf z}_\sigma^{4\over {n\,-\,2}} \cdot [\,w_1 \ - \ w_2\,]\ \bigg\vert \\[0.15in]
& \le &  C\,(n) \cdot \left\{ \ |\, w_1 \ - \ w_2\,| \left[ \ |\,{\bar w}_1 \,|^{4\over {n\,-\,2}} \ + \ |\, w_1\,|^{4\over {n\,-\,2}} \ + \ |\, w_2 \,|^{4\over {n\,-\,2}} \ \right] \ + \ |\, w_1 \ - \ w_2\,|^{{n\,+\,2}\over {n\,-\,2}} \ \right\}\ .
\end{eqnarray*}
This comes from inequalities of the following forms.
$$
\bigg\vert \  1 \ - \ (1 \,+\,T)^{{n \,+\,2}\over {n\,-\,2}} \ + \ {{n\,+\,2}\over {n\,-\,2}}
\cdot T \ \bigg\vert \ \le \ C_1 \cdot  |\,T|^{{n\,+\,2}\over {n\,-\,2}}\ \ \ \  \ \ \  \ \  \  \
(n \ \ge \ 6\,,\, \ \ |\,T\,| \ {\mbox{is \ \ small}}\,)\ , \leqno (3.28)\
$$

\vspace*{-0.2in}

$$
\bigg\vert \  a^{4\over {n \,-\,2}} \ - \ (a \,+\,s)^{4\over{n\,-\,2}} \ \bigg\vert
\ \le \ C_2 \cdot  s^{{4}\over {n\,-\,2}}\ \ \ \mfor \ \ a \ > \ 0\,, \
\ a\ +\ s \ >  \ 0\,, \ \ \ \ \   (n \ \ge \ 6)\,.\leqno (3.29)
$$
Here $\,C_1\,$ and $C_2\,$ depend on $\,n\,$ only. It follows that
\begin{eqnarray*}
(3.30)  \   & \ & |\ R_{{\bf z}_\sigma}({\bar w}_1\,+\,w_1) \ - \ R_{{\bf z}_\sigma}({\bar w}_1\,+\,w_2)  \ | \\[0.15in]
  & \ & \!\!\!\!\!\!\!\!\!\!\!\! \!\!\!\!\!\!\!\!\!\!\! \le \,  C\,(n) \!\cdot\! \Vert \, w_1 \, - \, w_2\,\Vert_\btd   \cdot \!\left[ \, \Vert\,{\bar w}_1\,\Vert_\btd^{4\over {n \,-\,2}}   + \ \Vert\,w_1\,\Vert_\btd^{4\over {n \,-\,2}}  + \ \Vert\,w_2\,\Vert_\btd^{4\over {n \,-\,2}}   + \ \Vert \, w_1 \, - \, w_2\,\Vert_\btd^{4\over {n\,-\,2}}\ \right] \cdot \Vert\,h\,\Vert_\btd\,. \ \ \ \ \ \
\end{eqnarray*}
Hence we can find a positive number $\,\gamma \ <  \ 1\,$ so that
$$
\Vert \, {\bf T}\,(w_1)\, -\, {\bf T}\,(w_2) \,\Vert \ \le \ \gamma \cdot \Vert \,  w_1 \, - \,w_2  \,\Vert_\btd\ , \ \ \  \ {\mbox{where}} \ \ \ \ \Vert\, w_1  \Vert_\btd \ \le \  c\ \ \ {\mbox{and}} \ \ \  \Vert\, w_2  \Vert_\btd \ \le \  c\,.
$$
Via the contraction mapping theorem, (3.21) has a unique fixed point $\ {\bar w}_2\ $ with
$$
\Vert \,{\bar w}_2 \,\Vert_\btd \ \le \ c\ .
$$
Thus we find a solution
$$
w_{{\bf z}_\sigma} \ = \ {\bar w}_1 \ +\  \ {\bar w}_2
$$
to auxiliary equation (3.14). Moreover,
$$
\Vert \, w_{{\bf z}_\sigma}  \Vert_\btd \ \le \ \Vert \, {\bar w}_1   \Vert_\btd \,+\, \Vert \, {\bar w}_2   \Vert_\btd \ \
\Longrightarrow \ \ \Vert \, w_{{\bf z}_\sigma}  \Vert_\btd \ = \ 2 \cdot  c\ .
$$
Thus we can take $\,{\tilde \epsilon}  \ = \ 2 \cdot c  \,$ in (3.15)\,.\,
With this, together with (3.27) and (3.30), we can also show that the small $\perp -\,$solution \,$w_{{\bf z}_\sigma}$\, depends on \,${\bf z}_\sigma$\, in a $\,C^1\,$ manner.
Cf. \cite{Progress-Book} and the proof of Proposition 4.2 in \cite{Pino}\,.\, This completes the proof of the theorem.  \qedwh

\vspace*{-0.2in}

{\bf{\S\ 3\,a}.}  \ \ {\it  Finite dimension reduction.}\ \
Let $\,w_{{\bf z}_\sigma}\,$ be the unique  small  $\perp\!-\,$solution of
 the auxiliary equation (as described in Theorem 3.13)\,.\, Consider the reduced functional, which depends on $\,(\lambda_1\,, \,\lambda_2\,; \  \xi_1\,, \ \xi_2)\ \in \ (\R^+ \times \R^+) \times (\R^n \times \R^n)\,:$ \\[0.1in]
(3.31)

\vspace*{-0.4in}

\begin{eqnarray*}
& \ &  {\bf I}_{{\cal R}_H} (\lambda_1\,, \,\lambda_2\,; \  \xi_1\,, \ \xi_2)\\[0.05in]
& = &   {\bf I}_{{\cal R}_H}    (\,{\bf z}_\sigma \,+\,w_{{\bf z}_\sigma})\
  =  \ {1\over 2}\int_{\R^n} \langle \,\btd\, (\,{\bf z}_\sigma\,+\,w_{{\bf z}_\sigma})\,, \ \btd\, (\,{\bf z}_\sigma\,
  +\,w_{{\bf z}_\sigma})\,\rangle \\[0.125in]
  & \ &  \ \ \ \ \ \ \ \ \ \     - \ {1\over 2}\!\cdot  \! (n\,-\,2)^2 \int_{\R^n}\,(\,{\bf z}_\sigma+w_{{\bf z}_\sigma}\,)^{{2n}\over {n\,-\,2}}_+ \,+\,\
   \left[\,-\,{{n\,-\,2}\over {2\,n}} \,\right] \,\cdot \int_{\R^n} ({\tilde c}_n \cdot H)\,(\,{\bf z}_\sigma\,+\,w_{{\bf z}_\sigma})_+^{{2n}\over {n\,-\,2}}\ .
\end{eqnarray*}
This finite dimensional reduced functional forms the main object in our study. We first show its link to the full functional (2.3)\,.

\vspace*{0.2in}

{\bf Lemma 3.32.} \ \ {\it Under the conditions in Theorem}\, 3.13\,,\,
{\it if}\, $\,{\bf z}_\sigma \ = \ V_{\lambda_1\,, \ \xi_1} \ + \ V_{\lambda_2\,, \ \xi_2}\,$ {\it is a critical point of the reduced functional in\,} (3.31)\,,\, {\it that is}
$$
\ \ {{\partial \,{\bf I}_{\cal R}}\over {\partial \lambda_k}}\,\bigg\vert_{(\lambda_1\,, \,\lambda_2\,; \  \xi_1\,, \ \xi_2)}   \
= \ {{\partial\, {\bf I}_{\cal R} }\over {\partial \xi_{k_\ell}}}\,\bigg\vert_{(\lambda_1\,, \,\lambda_2\,; \  \xi_1\,, \ \xi_2)}  \
= \ 0\leqno (3.33)
$$
 {\it{for}} $\,k \ = \ 1\,, \, 2\,,$ {\it and}  $\ \ell\  = \ 1\,,\, \cdot \cdot \cdot\,, \   n\,, $
{\it then\,}  $\,{\bf z}_\sigma \ +\ w_{{\bf z}_\sigma}\,$ {\it is a critical point of}\, (\,{\it the full functional}\,)\, $\,{\bf I}\,$,\, {\it that is\,,\,}

\vspace*{-0.345in}

$$
{\bf I}' \,(\,{\bf z}_\sigma \ +\  w_{{\bf z}_\sigma}) \ = \ 0\,. \leqno (3.34)
$$

\hspace*{0.5in}Using the smallness of $\Vert \,w_{{\bf z}_\sigma}\,\Vert_\btd\,$ provided by  (3.15) when $\,{\bf D}\,$ is large enough, the proof of Lemma 3.44 is similar to the proof of Theorem 2.8 in   \cite{III}\,. There one can also find information on the regularity of the solution $\,{\bf z}_\sigma \ +\  w_{{\bf z}_\sigma}\,$ in (3.34)\,, as well as the property that it can be transferred back to $\,S^n\,$ via (1.3)\,.

\newpage

{\bf \S\,3\,b.} \   {\it Degree and gradient.}  \ \ It is convenient and natural to work  with the  coupled quasi\,-\,hyperbolic gradient\,,\, denoted by  $\,(\lambda \cdot \btd)\,$ (introduce in \cite{III}) and defined by  \\[0.1in]
(3.35)
$$
({\lambda \cdot \btd})\,{\cal T} \,:=\,\left(\,\lambda_1 \cdot D_{1_o}\,{\cal T}\,,\,\lambda_1 \cdot D_{1_1}\,{\cal T}\,,\,\cdots\,,\,\lambda_1 \cdot D_{1_n}\,{\cal T}\,;\ \lambda_2 \cdot D_{2_o}\,{\cal T}\,,\,\lambda_2 \cdot D_{2_1}\,{\cal T}\,,\,\cdots\,,\,\lambda_2 \cdot D_{2_n}\,{\cal T}\,\right)
$$
for $\,{\cal T}\ \in\ C^1\,(\,\overline{\Omega}\,)\,$, where \\[0.1in]
(3.36) \ \ \ \ $\,\,\Omega\, \subset\, (\R^+ \times \R^+) \times (\R^n \times \R^n)\,$ is a bounded domain with smooth boundary $\,\partial\,\Omega\,$, and\\[0.075in]
\hspace*{0.67in} $\overline \Omega \ \subset\ (\R^+ \times \R^+) \times (\R^n \times \R^n)\,.\,$ \\[0.1in]
In (3.35),
\vspace*{-0.02in}
$$
\displaystyle{  D_{k_o}\ =\ {{\partial}\over {\partial \lambda_k}}\ \ \ {\mbox{and}}\ \ \ D_{k_\ell}\ =\ {{\partial}\over{\partial\,{\xi_k}_{|_\ell}}}\ \ \ \ \ {\mbox{for}}\ \ \ \ k\ =\ 1\,,\ 2\ \ \ {\mbox{and}}\ \ \ \ell\ =\ 1\,,\ \cdots\ ,\ n\,. }   \leqno (3.37)
$$
As usual
$$
\Vert\, {(\lambda \cdot \btd)\,{\cal T}}\,\Vert  \, = \, \sqrt{\,\left( \lambda_1 \cdot   D_{1_o}\,{\cal T} \right)^{2}\ + \ \sum_{\ell\,=\,1}^n \left(\lambda_1 \cdot D_{1_\ell}\,{\cal T} \right)^{2}\ + \ \left( \lambda_2 \cdot D_{2_o}\,{\cal T} \right)^{2}\ + \ \sum_{\ell\,=\,1}^n \left( \lambda_2 \cdot D_{2_\ell}\,{\cal T} \right)^{2}\,}\ \,.
$$

The following theorem can be shown by using the homotopy invariance of the degree. See \cite{III}\

\vspace*{0.2in}

{\bf Theorem 3.38.}\ \ \ {\it Let $\,\Omega\,$ be as described in the above, and } $\,{\cal F}\,,\ {\cal G}\ : \ (\R^+ \times \R^+) \times (\R^n \times \R^n) \ \to \ \R^{2\,(n\,+\,1)}\,$ {\it be of class $\,C^o (\,\overline{\Omega}\,)\,$, which satisfy }
$$
\!\!\!\!\!\!\!\!\!\!\!\!\!   \ \  \min_{\partial\,\Omega}\,\left\{\ \Vert\, {\cal F}\,\Vert \,\right\} \ >  \ 0\ , \leqno (3.39)
$$

\vspace*{-0.3in}

$$
\min_{\partial\,\Omega}\,\left\{\ \Vert\,{\cal F}\,\Vert\ \right\} \ \,>\, \ \max_{{\partial\,\Omega}}\left\{\ \Vert\ {\cal F} \ - \ {\cal G}\,\Vert\ \right\} \ . \leqno (3.40) \ \ \ \ \ \ \ {\it{and}}$$
{\it Then we have}
$$
{\bf{Deg}}\,\left[\,  {\cal G}\,,\ \Omega\,,\ {\vec{\,0}}\ \right] \ \,=\, \ {\bf{Deg}}\,\left[\, {\cal F}\,,\ \Omega\,,\ {\vec{\,0}}\ \right]\,.\leqno (3.41)
$$
{\it In particular} [\,{\it under  conditions}\, (3.39) \,{\it and}\, (3.40)\,]\,,\, {\it if\,} $\,{\bf{Deg}}\,\left[\, {\cal F}\,,\ \Omega\,,\ {\vec{\,0}}\ \right]\ \not=\ 0\,$,\, {\it then there is a point \,$\,p\ \in\ \Omega\,$ such that}
\vspace*{-0.025in}
$$
 {\cal G} \, (p) \ = \ {\vec{\,0}}\ . \leqno (3.42)
$$

\newpage

{\bf \S\,3\,b.} \ \ {\it Estimates on}  $\,w_{{\bf z}_\sigma}$\,.\, \ \  In order to extract effectively from the reduced functional from the key information (Proposition 4.1), we need the following estimates\,,\, which are shown in \cite{III}\,.\,

\vspace*{0.2in}

{\bf Lemma 3.43.}\ \ \  {\it Under the conditions in Theorem\,} 3.13\,,\, {\it let \,$\,w_{{\bf z}_\sigma}$ be the
unique small $\,\perp\!-\,$solution of the auxiliary equation\,} (3.14) [\,{\it which satisfies\,} (3.15)\,]\,. {\it We have}
\begin{eqnarray*}
\Vert\,w_{{\bf z}_\sigma}\Vert_\btd & \le & {\bar C}_2\,\cdot\,\left\{\,
 \left[\,\int_{\R^n} |\,H\,|^{{2n}\over {n\,+\,2}} \cdot {\bf z}_\sigma^{{2n}\over{n\,-\,2}}\,
 \right]^{{{n\,+\,2}\over 2n}}\ + \ \, O\,\left(\,{1\over { {\bf D}^{ \,{{n\,+\,2}\over 2}  \,-\ o\,(1)}}}\right)\,\right\}\ ,
\end{eqnarray*}
{\it and}
$$
\Vert\,\lambda_k \cdot D_{k_\ell}\,w_{{\bf z}_\sigma}\Vert_\btd \ \le \  {\bar C}_3\,\cdot\,\left\{\,  \left[\,\int_{\R^n} |\,H\,|^{{2n}\over {n\,+\,2}} \cdot {\bf z}_\sigma^{{2n}\over {n\,-\,2}}\,\right]^{{{n\,+\,2}\over 2n} \cdot {4\over {n\,-\,2}}}\ + \ O\,\left(\,{1\over { {\bf D}^{ \,2  \,-\ o\,(1)}}}\right) \,\right\} \ .
$$
{\it Here\,} $\,k\ =\ 1\,,\ 2\,;\ \ell\ = \ 0\,,\ 1\,,\ \cdots\,,\ n\,$. {\it In addition, $\,{\bar C}_2\,$ and \,$\,{\bar C}_3\,$ are positive constants independent on \,$\,(\,\lambda_1\,,\, \lambda_2\,;\ \xi_1\,,\ \xi_2\,)\,$ as long as the conditions in Theorem\,} 3.13 {\it are satisfied\,.}  [\,{\it Here \,$\,o\,(1)\ \to\ 0^{\,+}\,$ as \,$\,{\bf D}\ \to\ \infty\,$.}\, {\it Cf.} (2.18)\,.\,]\,

\vspace*{0.2in}

\hspace*{0.5in}Combining with Lemma 3.12\,,\, we obtain the following result.

\vspace*{0.2in}

{\bf Lemma 3.44.}\ \ \  {\it Under the conditions in Theorem\,} 3.13\,,\,
{\it assume also that }
$$
\lambda^\ell \ \le\ {{C_o }\over { {\bf D}^{\,n\,-\ 2} }}\ \ \ \ \ {\it{and}} \ \ \ \ \ 2 \ \le \ \ell \ < \ n \ - \ 2\,,\leqno (3.45)
$$

\vspace*{-0.15in}

{\it then we have}
 \begin{eqnarray*}
 (3.46) \ \ \ \ \ \ \ \ \ \ \ \ \ \ \ \ \Vert\,w_{{\bf z}_\sigma}\Vert_\btd & \le & {\bar C}_4 \cdot {1\over { {\bf D}^{ \,{{n\,+\,2}\over 2}  \,-\ o\,(1)}}}  \ \ \ \ \ {\it{and}} \ \ \ \\[0.125in]
 (3.47) \ \ \ \ \  \,\Vert\,\lambda_k \cdot D_{k_\ell}\,w_{{\bf z}_\sigma}\Vert_\btd & \  & \!\!\!\!\!\!\!\!\!\!\!\!\!\left(  \, \le \  {{C }
 \over { {\bf D}^{ \,2 \cdot {{n\,+\,2}\over {n\,-\,2}} }}} \ + \  {{C }\over { {\bf D}^4}}
\ + \  {{C }\over { {\bf D}^{ \,2  \,-\, o\,(1) }}}\right) \ \le \ {\bar C}_5 \cdot  {1\over { {\bf D}^{ \,2  \,-\  o\,(1) }}}\ \,. \ \ \ \ \ \ \ \ \ \ \ \ \ \ \ \ \ \
 \end{eqnarray*}
 {\it Here the positive constants $\,C\,$,\, $\,{\bar C}_4\,$ and $\,\,{\bar C}_5\,$ depend on $\,C_o\,$,\, but are independent on    \,$\,(\,\lambda_1\,,\, \lambda_2\,;\ \xi_1\,,\ \xi_2\,)\,$ as long as the conditions in Theorem\,} 3.13 {\it are satisfied\,.}

\vspace*{0.3in}

{\bf \Large {\bf  \S\,4.    }} {\Large{\bf     Extracting main information from the reduced }}\\[0.1in]
 \hspace*{0.5in} {\Large\bf{functional.}}\\[0.15in]
In the following we show that the main contributions  to  reduced functional (3.31) come from the critical points (the  twin  pseudo\,-\,peaks) and the interaction of the two bubbles [\,cf. (2.23) \ \& \ (2.25)\,]\,. Cf. \cite{Yan+}\,.\, In separating the key terms, it is here that the dimension restriction $\,n \ < \ 10\,$ comes in\,.\, In the following we denote the reduced functional by $\,{\bf I}_{{\cal R}_H}$,\, highlighting the dependence on $\,H\,.\,$


\newpage

{\bf Proposition 4.1.}\ \  {\it For \,$\,10 \ > \ n \ \ge \ 6\,$,\, assume the conditions in Theorem\,}  3.13\,. {\, \it There is a positive constant $\,{\bar D}_5\ (\,\ge \ {\bar D}_4\,)$ such that if }
$$
{\bf D} \ \ge \  {\bar D}_5\,,
$$

\vspace*{-0.25in}

{\it then we have}\\[0.075in]
(4.2)
\begin{eqnarray*}
 & \ &  \!\!\!\!\!\!\!\!\!\!\!\!\Bigg\vert \! \left(\lambda_1\!\cdot\!{{\partial}\over {\partial \lambda_1}} \right) \!{\bf I}_{{\cal R}_H}\,({\bf z}_\sigma)     \ - \    \left\{ \,
 {\bar C}^+_a\! \cdot\! [\,-\, {\varpi}_1\,] \!\cdot\! \lambda_1^\ell
  \ - \    {\bar C}^+_b \!\cdot\! {{\lambda_1^{{ n\,-\,2 }\over 2} \!\!\cdot\! \lambda_2^{{ n\,-\,2 }\over 2} }
 \over  { \gamma^{\,n\,-\ 2} }} \,\right\}\!\cdot\!  [\,1  \, + \, o\,(1)\,] \,\Bigg\vert \, \le \,   o\,(1) \cdot {1\over  { {\bf D}^{\,n\,-\ 2} }}\ ,\\
\end{eqnarray*}
(4.3)
\begin{eqnarray*}
& \ &\!\!\!\!\!\!\!\!\!\!\!\! \Bigg\vert\! \left(\lambda_1\!\cdot\!{{\partial}\over {\partial {\,\xi_{1_j}}}}\!         \right)  \!{\bf I}_{{\cal R}_H} ({\bf z}_\sigma)   \ -\
\left\{  {\bar C}^+_c \! \cdot  {\varpi}_1 \!\cdot\! \lambda_1^\ell \cdot {{ \xi_{1_j} }\over {\lambda_1}}  \,+\,      {{  {\bar C}^+_d }\over {   \lambda_2}} \!\cdot\! [\,\xi_{1_j} \!-  \xi_{2_j}\,]
  \!\cdot\!  {1\over {{\bf D}^{n} }}    \right\}  \!\cdot\! [\, 1   \,+\,  o\,(1)\,]  \, \Bigg\vert  \, \le \,     o\,(1) \cdot \,{1\over  { {\bf D}^{\,n\,-\ 2} }}\\
 \end{eqnarray*}
{\it for $\,j\ = \ 1\,, \ 2\,, \ \cdot \cdot \cdot\,, \ n\,.\,$ Here\ }({\it as usual}\,) {\it\,$\,{\bf z}_\sigma\ =\ V_{\lambda_1\,,\ \xi_1}\ +\ V_{\lambda_2\,,\ \xi_2}\,$.\, }{\it In the above\,,\,}  {\it\,$\,o\,(1)\ \to\ 0^{\,+}\,$ when \,$\,{\bf D}\ \to\ \infty\,$. Similar estimates hold  for derivatives in $\,\lambda_2\,$ and $\ \xi_2\,,\,$ with the same positive  constants $\,  {\bar C}^+_a\,, \  {\bar C}^+_b\,, \  {\bar C}^+_c\,$ and  $\  {\bar C}^+_d \,$} ({\it which depend only on $\,n\,$ and $\,\ell\,$}) {\it  as in\,} (4.2) {\it and\,} (4.3)\,.\,

\vspace*{0.2in}

{\bf Proof.} \ \ Arguing as in the proof of Proposition 2.11 in \cite{III}\,,\, and using (2.23)\,--\,(2.25)\,,\, we have (\,see {\bf \S\,A.2}\, in the \,{\bf e}\,-\,Appendix)\\[0.1in]
(4.4)
\begin{eqnarray*}
& \ &\!\!\!\!\!\!\!\!\!\!\!\!\!\!\! \Bigg\Vert  \left(\lambda_1\cdot{{\partial}\over {\partial \lambda_1}} \right) \,{\bf I}_{{\cal R}_H} (\,{\bf z}_\sigma) \ -
 \left\{  \left(\lambda_1\cdot{{\partial}\over {\partial \lambda_1}} \right)\, {\bf G}\,({\bf z}_\sigma)
 \ - \ \omega_n \cdot {{(n\,-\,2)^2}\over {2 }} \cdot \,{{\lambda_1^{{ n\,-\,2 }\over 2} \cdot \lambda_2^{{ n\,-\,2 }\over 2} }
 \over  { \gamma^{\,n\,-\ 2} }} \cdot [\,1 \ + \ o\,(1)\,]\,\right\} \Bigg\Vert
\\[0.2in]
&  \  &   \ \ \ \ \ \le \ o\,(1) \cdot \,{1\over  { {\bf D}^{\,n\,-\ 2} }}\ ,
\end{eqnarray*}
(4.5)
  \begin{eqnarray*}
  & \ & \!\!\!\!\!\!\!\!\!\!\!\!\!\!\!\!\!\Bigg\Vert\, \left(\lambda_1\cdot{{\partial}\over {\partial {\,\xi_{1_1}}}}         \right) \,{\bf I}_{{\cal R}_H} (\,{\bf z}_\sigma) \ - \
 \left\{\, \left(\lambda_1\cdot{{\partial}\over {\partial {\,\xi_{1_1}}}} \right) {\bf G}  \,({\bf z}_\sigma)
 \ + \   \lambda_1 \cdot  {{ C^+(n)}\over { \lambda_1 \cdot \lambda_2}} \cdot [\,\xi_{1_1} - \ \xi_{2_1}\,]
  \cdot  {1\over {{\bf D}^{n} }}   \,\right\}\ \Bigg\Vert\\[0.2in]
&  \  &    \le \ o\,(1) \cdot \,{1\over  { {\bf D}^{\,n\,-\ 2} }}\ \,.
 \end{eqnarray*}
Here $\,o\,(1) \ \to \ 0\,$ as $\,{\bf D} \ \to \ \infty\,.\,$  We continue from (4.4)\,:\\[0.1in]
(4.6)
\begin{eqnarray*}
& \ & \!\!\!\!\!\!\!\!\!{\partial\over {\partial \lambda_1}} \, {\bf G}  \,({\bf z}_\sigma) \ = \  -\, {{n\,-\,2 }\over { 2\,n }} \cdot {\partial\over {\partial \lambda_1}} \,\int_{\R^n} ({\tilde c}_n \cdot H) \cdot  \left( \,V_1 \ + \ V_2\,
 \right)^{{{2n}\over {n\,-\,2}}   }   \  \ \  \ (\,V_1 \ = \ V_{\lambda_1\,, \ \xi_1}\,, \   V_2 \ = \   V_{\lambda_2\,, \ \xi_2}\,)\\[0.15in]
 & =  &   -  \int_{\R^n} ({\tilde c}_n \cdot H) \!\cdot\! \left( \,V_1 \ + \ V_2\, \right)^{{{n\,+\,2}\over {n\,-\,2}}   }
\cdot {{\partial V_1}\over {\partial \lambda_1}}\\[0.15in]
& = & - \int_{\R^n} ({\tilde c}_n \cdot H) \!\cdot\!  V_1^{{{n\,+\,2}\over {n\,-\,2}}   }
\cdot {{\partial V_1}\over {\partial \lambda_1}} \cdot [\,1\ + \ o\,(1)\,] \ \ \   ( {\mbox{similar \ \ to     \ \ Weak  \ \ Interaction \ \  Lemma \ \ 2.16}}\,)\\[0.15in]
& = &   -  \int_{\R^n} [\, ({\tilde c}_n \cdot H)\,] \!\cdot\!  V_1^{{{n\,+\,2}\over {n\,-\,2}}   }
\cdot  \left\{\ -\, {{n - 2}\over 2}   \cdot \lambda_1^{{n -\, 4}\over 2} \cdot {{ \lambda_1^2 \ - \ |\,y\,-\,\xi_1\,|^2  }\over {(\,\lambda^2
\,+\, |\, y \,-\, \xi_1|^{\,2})^{{n }\over 2} }} \,\right\}\ dy\cdot [\,1\ + \ o\,(1)\,]\\[0.15in]
& \ & \hspace*{.7in} [\  {\mbox{direct \ \ calculation \ \ from  \ \ (2.9)}} \ \ \uparrow \ ; \ \ \ \ \,o\,(1)\ \to\ 0^{\,+}\,\ \  {\mbox{when}} \ \ \,{\bf D}\ \to\ \infty\,]\\[0.15in]
& = &  {{n - 2}\over 2}   \cdot  \int_{\R^n} [\, ({\tilde c}_n \cdot H)\,] \!\cdot\!  V_1^{{{n\,+\,2}\over {n\,-\,2}}   }
\cdot  \left\{\  \lambda_1^{{n \,-\, 4}\over 2} \cdot {{ 2\,\lambda_1^2 \ - \ [\ \lambda^2_1 \ + \ |\,y\,-\,\xi_1\,|^2 \,] }\over {(\,\lambda^2
\,+\, |\, y \,-\, \xi_1|^{\,2})^{{n }\over 2} }} \,\right\}\ dy\cdot [\,1\ + \ o\,(1)\,]\\[0.15in]
& = &  {{n - 2}\over 2} \cdot   \int_{\R^n} [\, ({\tilde c}_n \cdot H)\,] \left( {{ \lambda_1}\over {  \lambda^2_1 \, + \, |\,y\,-\,\xi_1\,|^2 }} \right)^{\!\!{{n\,+\,2}\over {\,2}}} \!\!
* \\[0.15in]
& \ &\ \ \ \ \ \ \ \ \ \ \ \  \ \  *\, \left\{  2 \cdot \left( {{ \lambda_1}\over {  \lambda^2_1 \, + \, |\,y\,-\,\xi_1\,|^2 }} \right)^{\!\!{n\over {\,2}}}\!
  \ - \  {1\over {\lambda_1}} \cdot  \left( {{ \lambda_1}\over {  \lambda^2_1 \ + \ |\,y\,-\,\xi_1\,|^2 }}
\right)^{\!\!{{n\,-\ 2}\over 2}}\,\right\} \!\cdot [\,1\, + \, o\,(1)\,] \ dy  \\[0.15in]
&  = &    {{n\,-\,2}\over 2} \cdot   \!\! \int_{\R^n}    [\,  ({\tilde c}_n \cdot H)\,]\cdot \left[ \ 2 \cdot \left( {{\lambda_1}\over { \lambda_1^2 \ + \ |\,y\,|^2}}
\right)^{\!\!{{n\,+\,2}\over 2} \,+\, {n\over 2}}     -    {1\over {\lambda_1}} \cdot  \left( {{\lambda_1}\over { \lambda_1^2 \ + \ |\,y\,|^2 }}
\right)^{\!\!n} \ \right] \, dy \cdot [\,1\ + \ o\,(1)\,] \\[0.2in]
& \ &  \!\!\!\!\!\!\!\!\!\!\!\!\left[\,   {\mbox{cf. \  (1.8)\,, \ \ (1.10) \  and \   (2.4) \ }}\uparrow\ ; \ \    \ \ \ \ \ \ \  \ \ \uparrow \ \ {\mbox{similar \ \ to \ \ (A.10)}} \ \ \ \ \ \ \ \ \  \downarrow \ \  Y \ = \ {y\over \lambda}  \ \right]\\[0.2in]
&  = &   {{n\,-\,2}\over 2} \cdot   {{\lambda_1^\ell }\over { \lambda_1 }} \cdot \int_{\R^n}  \ [\,{\bf P}^\ell_1\,(Y)\,] \cdot \left[ \,2 \cdot \left( {1\over { 1 \ + \ |\,Y\,|^2 }}
\right)^{\!\!n\,+\,1}   \ - \     \left( {1\over { 1 \ + \ |\,Y\,|^2 }}
\right)^{\!\!n} \ \right] \cdot\, dY \cdot [\,1\ + \ o\,(1)\,]  \\[0.15in]
& \ & \hspace*{4in} \ \ \ \ \ \ \ \ \ \ [\,{\mbox{via \ \ (1.8) \ \& \ (1.10)}}\,]\\[0.15in]
& = &  -\,{{n\,-\,2}\over 2} \cdot   {{\lambda_1^\ell }\over { \lambda_1 }} \cdot \int_{\R^n}  \ [\,{\bf P}^\ell_1\,(Y)\,]  \cdot \left[ \ \left( {1\over { 1 \ + \ |\,Y\,|^2 }}
\right)^{\!\!n}    - \,2 \cdot     \left( {1\over { 1 \ + \ |\,Y\,|^2 }}
\right)^{\!\!n\,+\,1} \ \right] \cdot\, dY \cdot [\,1\ + \ o\,(1)\,]\,.
\end{eqnarray*}
Here $\,0\,(1) \ \to \ 0\,$ as $\,{\bf D} \ \to \ \infty\,.\,${{\hfill {$\rlap{$\sqcap$}\sqcup$}}

\newpage

\hspace*{0.5in}Now we apply the following result, whose proof is similar to the proof of Lemma A.6.8 in \cite{2nd}\,;\, see also \S\,{\bf A.2} in the {\bf e}\,-\,Appendix\,. \\[0.2in]
  {\bf Lemma 4.7 (Reduction Lemma)\,.} \ \  {\it In} $\,\R^n\,,\,$ $n \ \ge \ 3$\,,\,  {\it consider a homogeneous polynomial $\,{\cal Q}_{\ell}\,$}  {\it of even degree\,} $\,\ell \ \le \ n \,-\, 1$\,.\,  {\it  We have}
 $$
\int_{\R^n} \   {\cal Q}_{\ell}\, (y) \cdot \left( {1\over {1+ |\,y|^{\,2}}} \right)^{\!n} \,d y \ = \ {{J_n}\over {   \ell \cdot (\ell \,-\,2) \cdot \cdot \cdot \,2\cdot 1 }}\cdot  [\,\Delta_{\,o}^{\!(h_\ell)}\ Q_\ell\,]\ \ \ \ \ \ (h_\ell \ = \ \ell/\,2)\,, \leqno (4.8)
$$
 {\it where}
 $$
 \ \ \ \ \ \ \ \ \ \ J_n \ = \ \int_{R^n} y_1^2 \cdot \cdot \cdot y_{h_\ell}^2 \cdot \left( {1\over {1+ |\,y|^{\,2}}} \right)^{\!n}  d y \ \ \ \ \ \ \ \  [\,y \ = \ (y_1\,, \cdot \cdot \cdot\,, \ y_{h_\ell}\,, \ \cdot \cdot \cdot\,, \ y_n)\,]\ . \leqno (4.9)
 $$
{\it Likewise, }
$$
\int_{\R^n} \   {\cal Q}_{\ell}\, (y) \cdot \left( {1\over {1+ |\,y|^{\,2}}} \right)^{\!n\,+\,1}  d y \ = \ {{J_{n\,+\,1}}\over {   \ell \cdot (\ell \,-\,2) \cdot \cdot \cdot \,2\cdot 1 }}\cdot [\, \Delta_{\,o}^{\!(h_\ell)}\ Q_\ell\,]\ , \leqno (4.10)
$$
 {\it where}
 $$
 J_{n\,+\,1} \ = \ \int_{R^n} y_1^2 \cdot \cdot \cdot y_{h_\ell}^2 \cdot \left( {1\over {1+ |\,y|^{\,2}}} \right)^{\!n\,+\,1}  d y \ \ \ \ (\,\Longrightarrow \ \ J_n \ - \  J_{n\,+\,1} \ > \ 0\,)\,. \leqno (4.11)
 $$

\vspace*{0.15in}

{\it Proof of Proposition} 4.1 {\it continues...} \ \ It follows from (4.4) and Lemma 4.5 that
$$
 \left(\lambda_1\cdot{{\partial}\over {\partial \lambda_1}} \right)  {\bf G}_H\,({\bf z}_\sigma)  \ = \ {\bar C}^+_1\,(n) \cdot [\,-\, {\varpi}_1\,] \cdot \lambda_1^\ell \cdot [\,1  \ + \ o\,(1)\,]  \ \ \ \ \ \ .  \leqno (4.12)
$$
where
 $$
{\bar C}^+_1\,(n\,,\,\ell)  \ = \ {{n\,-\,2}\over 2}\cdot {{ J_n \ -\ 2 \, J_{n\,+\,1} }\over { \, \ell \cdot (\ell \,-\,2) \cdot \cdot \cdot \,2\cdot 1  \,}}\,. \leqno (4.13)
 $$
To show that
$$
 J_n \ -\ 2 \, J_{n\,+\,1} \ > \ 0\,,\, \leqno (4.14)
$$
consider the stereographic projection
\begin{eqnarray*}
(4.15) \ \ \ \ \   \ \ \   \dot{\cal P} : S^n \setminus \{ {\bf N} \}\, &\to &  {\R}^n\,\\
x & \mapsto & y \ = \ \dot{\cal P}(x)\,,\ \ \ \
{\mbox{where}} \ \    y_i \ = \ {{x_i}\over{1 \,-\, x_{n + 1}}}\,,
\ \ \ 1 \ \le\  i \ \le \ n\,.\ \ \ \ \   \ \ \ \ \ \ \ \ \ \   \ \ \ \ \
\end{eqnarray*}
Here $\,x \,=\,(x_1\,,\, \cdot \cdot \cdot, \ x_{n + 1})
\in \,S^n\, \setminus \{ {\bf N} \}\,$ and $\,{\bf N} \,=\, (\,0, \ \cdot \cdot \cdot, \ 0, \ 1)\,$ is the north pole. Conversely\,,
\begin{eqnarray*}    x_i \ = \  {{2\,y_i}\over{1 \,+\, r^2}}\,, \ \ \ \ \  1 \ \le\  i \ \le\ n\,,\ \ \ \ \
{\mbox{and}} \ \ \ \ \ \ \ \ \
x_{n + 1} \ = \  {{r^2 \,-\,1 }\over{r^2 \,+\, 1}}\,,  \ \ \ \ \ \ \
{\mbox{where}} \ \  r \ = \ |\,y|\,.\,
\end{eqnarray*}
It is known that $\,\dot{\cal P}\,$ is a conformal map between $(S^n \setminus \{ {\bf N} \}\,, \ g_1)$ and $\,(\R^n,\ g_o)\,.\,$ The conformal factor is given by
$$
g_1 \,(x) = \left[ {{4}\over {(1 + r^2)^2}} \right]\cdot g_o \,(y) \ \ \mfor \ \ y = \dot{\cal P}(x)\,.
$$
Back to checking (4.14)\,:
 \begin{eqnarray*}
 J_n \ - \ 2\,J_{n\,+\,1} \ & = &  \int_{R^n} y_1^2 \cdot \cdot \cdot y_{h_\ell}^2 \cdot \left( {1\over {1\,+\, r^{\,2}}} \right)^{\!n} \cdot  \left[\ 1\  - \ {2\over {1\ + \ r^2 }}\,\ \right] \ d y\\[0.15in]
 & = &  \int_{R^n} y_1^2 \cdot \cdot \cdot y_{h_\ell}^2 \cdot \left( {1\over {1\,+\, r^{\,2}}} \right)^{\!n} \cdot  \left[\  {{r^2 \ - \ 1}\over { r^2\ + \ 1 }}\,\ \right] \ d y\\[0.12in]
 & = & {1\over {2^n}} \cdot \int_{S^n} \left[\, {{x_1^2}\over {(1\,-\,x_{n\,+\,1})^2 }} \right]\cdot \cdot \cdot \left[\, {{x_{h_\ell^2}}\over {(1\,-\,x_{n\,+\,1})^2 }} \,\right] \cdot x_{n\,+\,1} \ dS\\[0.12in]
 & = & {1\over {2^n}} \cdot \int_{S^n_+} \left[\, {{x_1^2}\over {(1\,-\,x_{n\,+\,1})^2 }} \right]\cdot \cdot \cdot \left[\, {{x_{h_\ell^2}}\over {(1\,-\,x_{n\,+\,1})^2 }} \,\right] \cdot x_{n\,+\,1} \ dS\\[0.12in]
 &  \  & \ \ \ \ \ + \ {1\over {2^n}} \cdot \int_{S^n_-} \left[\, {{x_1^2}\over {(1\,-\,x_{n\,+\,1})^2 }} \right]\cdot \cdot \cdot \left[\, {{x_{h_\ell^2}}\over {(1\,-\,x_{n\,+\,1})^2 }} \,\right] \cdot x_{n\,+\,1} \ dS\\[0.12in]
 & = & {1\over {2^n}} \cdot \int_{S^n_+} \left[\, {{x_1^2}\over {(1\,-\,x_{n\,+\,1})^2 }} \right]\cdot \cdot \cdot \left[\, {{x_{h_\ell^2}}\over {(1\,-\,x_{n\,+\,1})^2 }} \,\right] \cdot x_{n\,+\,1} \ dS\\[0.12in]
 &  \  &   + \ {1\over {2^n}} \cdot \int_{S^n_+} \left[\, {{x_1^2}\over {(1\,-[\,-\,\,x_{n\,+\,1}])^2 }} \right]\cdot \cdot \cdot \left[\, {{x_{h_\ell^2}}\over {(1\,-[\,-\,\,x_{n\,+\,1}])^2 }} \,\right] \cdot [\,-\,x_{n\,+\,1}] \ dS\,.\\
 \end{eqnarray*}
Here $\,S^n_+\,$ is the upper hemisphere\,,\, and $\,S^n_-\,$   the lower\,.\, Consider a fixed point
\begin{eqnarray*}
& \ & \ \  (\,x_1\,,\,\cdot \cdot \cdot\,,\ x_{h_\ell}\,,\,\cdot \cdot \cdot\,,\,x_{n\,+\,1}\,) \ \in \ S^n_+ \\[0.12in]
{\mbox{and \ \ its \ \ ``reflection\," \ :}} \  \ \ & \ & \ \  (\,x_1\,,\,\cdot \cdot \cdot\,,\ x_{h_\ell}\,,\,\cdot \cdot \cdot\,,\,[\,-\,x_{n\,+\,1}]\,) \ \in \ S^n_-\ . \ \ \ \ \ \  \ \ \  \  \ \ \ \ \ \ \ \ \ \ \ \ \ \ \ \ \ \ \ \ \ \ \
\end{eqnarray*}
For $\, 0 \ <  \ x_{n\,+\,1} \ <  \ 1\,,\,$ we have
$$
\left[\, {{x_1^2}\over {(\,1\,-\,x_{n\,+\,1})^2 }} \right] \,\ge \ \left[\, {{x_1^2}\over {(\,1\,-[\,-\,x_{n\,+\,1}])^2 }} \right]\,, \ \cdot \cdot \cdot\,,\ \ \left[\, {{x_{h_\ell}^2}\over {(\,1\,-\,x_{n\,+\,1})^2 }} \right] \,\ge \ \left[\, {{x_{h_\ell}^2}\over {(\,1\,-[\,-\,x_{n\,+\,1}])^2 }} \right]\,.
$$
Thus (4.14) holds.

\newpage

\hspace*{0.5in}For  derivative  in $\,\lambda_2\,,\,$ we have similar expression  with the same constant as in (4.13)\,.\,
Likewise (see \S\,A.4 in the {\bf e}\,-\,Appendix),
$$
\left(\lambda_1\cdot{{\partial}\over {\partial \xi_{1_j} }} \right)\, {\bf G}_H\,({\bf z}_\sigma) \ = \ {\bar C}^+_2\,(n\,,\,\ell) \cdot  {\varpi}_1 \cdot  \lambda_1 \cdot {{ \xi_{1_j} }\over {\lambda_1}} \cdot {{\lambda_1^\ell }\over { \lambda_1 }} \cdot [\, 1 \ + \ o\,(1)\,] \   , \leqno (4.16)
$$
with similar expression for derivatives in $\,\xi_2\,.\,$  \qed

{\bf \Large {\bf  \S\,5. \ \   The target - solving the equations - balancing the }}\\[0.1in]
\hspace*{0.68in}{\bf \Large {\bf  local and global contributions.}}

\vspace*{0.2in}

In view of     (4.2) and (4.3), and Theorem 3.38\,,\, our attention is drawn to the terms in the brackets in (4.2) and (4.3)\,.\, Thus consider the map
\begin{eqnarray*}
(5.1) \ \ \ \ \ \  \     {\cal T}\,( \, {\bf P}\, ) & = & (\,T_{1_o}\,, \ T_{2_o}\,;\, \ T_{1_1}\,, \ \cdot \cdot \cdot\,, \  T_{1_n}\,;\, \  T_{2_1}\,, \ \cdot \cdot \cdot\,, \  T_{2_n}\,) \ \in\ \R^{2\,(n\,+\,1)}\ , \\[0.15in]
{\bf P} & = &  (  \lambda_1\,,\    \lambda_2\,; \ \xi_{1_1}\,, \ \cdot \cdot \cdot\,, \ \xi_{1_n}\,; \,   \ \xi_{2_1}\,, \ \cdot \cdot \cdot\,, \ \xi_{2_n}\, ) \ \in \ (\R^+ \times \R^+) \times (\,\R^n \times \R^n)\ ,  \ \ \ \ \ \  \ \ \  \ \   \ \ \ \ \ \  \ \ \  \ \
\end{eqnarray*}
where
\begin{eqnarray*}
(5.2) \ \ \ \ \ \ \ \ \ \  T_{1_o} ( \, {\bf P}\, ) & = &  C^+_1  \cdot   |\, {\varpi}_1\,| \cdot  \lambda_1^{\ell }    \ - \,     {{ (\sqrt{\lambda_1 \cdot \lambda_2\,}\ )^{n\,-\,2}
}\over { \gamma^{\,n\,-\,2}  }} \ \ \ \ \ \    \left(   C^+_1   \ = \ {{ {\bar C}_a^+ }\over { {\bar C}_b^+ }}  \right)\ ,\\[0.15in]
(5.3) \ \ \ \ \ \ \ \ \ \  T_{2_o}( \, {\bf P}\, ) & = &  C^+_1  \cdot    |\, {\varpi}_2\,| \cdot   \lambda_2^{\ell }     \ - \,   {{ (\sqrt{\lambda_1 \cdot \lambda_2\,}\ )^{n\,-\,2}
}\over { \gamma^{\,n\,-\,2}  }} \ \ \     \left(\downarrow\ \,   C^+_2   \ = \ {{ {\bar C}_c^+ }\over { {\bar C}_b^+ }}  \ \ \,  {\mbox{and}} \ \    C^+_3   \, = \, {{ {\bar C}_d^+ }\over { {\bar C}_b^+ }} \right),\ \ \ \ \ \ \ \ \ \ \\[0.2in]
(5.4) \ \ \ \ \ \ \ \ \ \  T_{1_{\,j}}( \, {\bf P}\, ) & = &   C^+_2  \cdot    {\varpi}_1  \cdot  \lambda_1^{\ell } \cdot {{ \xi_{1_j} }\over {\lambda_1}} \ + \ \, C^+_3  \cdot
{1\over { \lambda_2}} \cdot {{ (\sqrt{\lambda_1 \cdot \lambda_2\ }\ )^{n}
}\over {\gamma^{\,n}  }} \cdot  (\,{\xi}_{\,1_j} \, - \,{\xi}_{\,2_j}\,)\ , \ \ \ \  {\mbox{and}} \ \ \ \ \  \ \ \ \ \
\\[0.2in]
(5.5) \ \ \ \ \ \ \ \ \ \  T_{2_{\,j}}( \, {\bf P}\, ) & = &  C^+_2   \cdot    {\varpi}_2  \cdot   \lambda_2^{\ell }  \cdot {{ (\xi_{2_j} \ - \ {\bf q}_{\,2_j}\,) }\over {\lambda_2}} \ + \ \, C^+_3  \cdot {1\over {\lambda_1  }} \cdot {{ (\sqrt{\lambda_1 \cdot \lambda_2\,}\ )^{n}
}\over { \gamma^{\,n}  }} \cdot  (\,\xi_{2_j} \, - \,{\xi}_{\,1_j}\,)\,, \ \ \ \ \  \ \ \ \ \   \ \ \ \ \  \ \ \ \ \
\end{eqnarray*}
\hspace*{4in}\ \ \ \ \ \ \ \ \ \ \ \ \ where $\,j \ = \ 1\,,\ 2\,, \ \cdot \cdot \cdot\,, \ n\,.\,$ \\[0.01in]To find
$$\,\lambda_1 \,=\,\lambda_{1_\tau} \ \ \ \ {\mbox{and}} \ \ \ \ \,\lambda_2 \,=\,\lambda_{2_\tau} \ \ \ \ {\mbox{so \ \ that}} \ \ \,T_{1_o} \ = \ T_{2_o} \ = \ 0\,,\, \leqno (5.6)
$$

\vspace*{-0.1in}

we let

\vspace*{-0.3in}

$$
\lambda_2 \ = \ \alpha \cdot \lambda_1\,. \leqno (5.7)
$$
Here the positive constant $\,\alpha \,$ is chosen so that
\begin{eqnarray*}
(5.8) \ \ \   |\, {\varpi}_1\,| \ = \ |\, {\varpi}_2\,|  \cdot \alpha ^\ell  \ \ \ \Longleftrightarrow \ \
\alpha  \ = \ \left({{    {\varpi}_1 }\over {  {\varpi}_2 }}  \right)^{\!\!{1\over {\,\ell}}}
&  \Longrightarrow &  \lambda_{1_\tau}  \ = \ \left[\,\gamma^{\,n\,-\,2} \cdot {{ C_1^+  \cdot   |\, {\varpi}_1\,|}\over {\alpha ^{{n\,-\,2}\over 2} }}\,
 \right]^{1\over {\,(n\,-\,2) \ - \ \ell\,}}  \ \ \    \ \ \    \ \ \    \ \ \   \\[0.1in]
 & \ &  \ \ \ \ \ [ \ \uparrow \ \ {\mbox{via \ \ (5.2)}}\,] \\[0.05in]
 &  \Longrightarrow & \lambda_{2_\tau} \ = \ \left({{   {\varpi}_1  }\over {   {\varpi}_2 }}  \right)^{\!\!{1\over {\,\ell}}}  \cdot \lambda_{1_\tau}\ .
\end{eqnarray*}
Note that $\,\lambda_{1_\tau} \ \ \& \ \,\lambda_{2_\tau} \ \to \ 0^+$\, as $\,\gamma\ \to \ 0^+\,.\,$  (5.7) and (5.8) imply that
$$
C_p^{-\,{1\over \ell}} \cdot \lambda_2 \ <  \ \lambda_1 \ <  \ C_p^{1\over \ell} \cdot \lambda_2  \leqno (5.9)
$$
$$
\lambda^\ell \ \le\ C_o \cdot {1\over { {\bf D}^{\,n\,-\ 2} }}\ \ \ \ \   \ \ \ \ \ [\ 2 \ \le \ \ell \ < \ n \ - \ 2\,]\leqno (5.10) \ \ \ \ \ \ \ \ \ \ \ {\mbox{and}}
$$
for $\,(\lambda_1\,, \ \lambda_2\,; \ \xi_1\,,\ \xi_2)\,$ close to $\,(\lambda_{1_\tau}\,, \ \lambda_{2_\tau}\,; \ 0\,,\ {\bf p}_2\,)\,$.\, Cf. (3.45)\,.\, Here the positive constant $\displaystyle{\,C_p\,\approx\,{{   {\varpi}_2 }\over {   {\varpi}_1 }} \,}$.\,   \,\bk
Next, to seek
$$
\,\xi_{1_j} \,=\,\xi_{1_{j_\tau}} \ \ \ \ {\mbox{and}}  \ \ \ \  \xi_{2_j} \,=\,\xi_{2_{j_\tau}} \ \ \ \  {\mbox{so \ \  that}} \ \ \ \ \,T_{1_j} \ = \ T_{2_j} \ = \ 0\,
$$
for  $\,j \ = \ 1\,,\ 2\,, \ \cdot \cdot \cdot\,, \ n\,,\,$ we let
$$
\xi_{1_{j_\tau}} \ = \ -\, \beta  \cdot (\,\xi_{2_{j_\tau}} \ - \ {\bf q}_{\,2_j}\,) \ \ \ \ \ \ {\mbox{for}} \ \ \ \ j \ = \ 1\,, \ 2\,, \ \cdot \cdot \cdot\,, \ n\ , \leqno (5.11)
$$
where the positive constant $\,\beta\,$ (\,independent on $\,j\,$\,) is chosen so that
$$
  {\varpi}_1  \cdot  \lambda_{1_\tau}^{\ell } \cdot {{ \lambda_{2_\tau} }\over {\lambda_{1_\tau}}} \cdot \beta  \ = \ {\varpi}_2  \cdot   \lambda_{2_\tau}^{\ell }  \cdot {{ \lambda_{1_\tau}  }\over {\lambda_{2_\tau}}}\ \ \ \left[ \ \Longrightarrow \ \ (5.4) \ ``=" \ (5.5)\,\right]\ . \leqno (5.12)
$$
Note that $\,\beta \ = \ O(1)\,$\, via (5.7) and  (5.8)\,.\,
Thus we find $\,\xi_{1_{j_{\,\tau}}}\,$ via (5.4) by writing
\begin{eqnarray*}
(5.13) \ \ \   C_2^+  \cdot  {\varpi}_1  \cdot  \lambda_{1_\tau} ^{\ell } \cdot {{ \xi_{1_{j_\tau}} }\over {\lambda_{1_\tau}}} \ = \  -\  C^+_3  \cdot
{1\over { \lambda_{2_\tau}}} \!\!\!\!\!&\cdot&\!\!\!\!\!{{ \left(\sqrt{\lambda_{1_\tau} \cdot \lambda_{2_\tau}\ }\ \right)^{n}
}\over {\gamma^{\,n}  }} \cdot  (\,\xi_{1_{j_\tau}} \  - \ {\bf q}_{\,2_j} \, - \,[\, \xi_{2_{j_\tau}}\ - \ {\bf q}_{\,2_j}\,]\,)  \\[0.15in]
\Longrightarrow \ \ \   {{ C_2^+}  \over {  C_1^+ }}  \cdot    {{ \left(\sqrt{\lambda_{1_\tau}\cdot \lambda_{2_\tau}\,}\ \right)^{n\,-\,2}
}\over { \gamma^{\,n\,-\,2}  }}  \cdot {{\xi_{1_{j_\tau}} }\over {\lambda_{1_\tau}}}   & = &  \!\!C^+_3  \cdot {1\over \alpha} \cdot
  {{ \left(\sqrt{\lambda_{1_\tau} \cdot \lambda_{2_\tau}\ }\ \right)^{n}
}\over {\gamma^{\,n}  }} \cdot  \left[\, {{ \xi_{1_{j_\tau}} } \over { \lambda_{1_\tau}}} \cdot \left(1  \ +  {1\over {\beta}} \right) \,  - \,  {{ {\bf q}_{\,2_j}  } \over { \lambda_{1_\tau}}} \,\right]\ \ \ \ \ \ \ \ \ \ \ \ \ \ \ \ \ \ \ \ \ \ \ \   \ \ \ \ \ \ \ \  \ \ \ \     \ \ \\[0.15in]
\Longrightarrow \ \   \ \ \ \ \  {{\xi_{1_{j_\tau}} }\over {\lambda_{1_\tau}}} \cdot \Bigg\{ \,{{ C_2^+}  \over {  C_1^+ }}  \ - \  C_3^+ \cdot {1\over { {\bf D}^2 }} \cdot {1\over \alpha} \!\!\!\!\!&\cdot&\!\!\!\!\! \left(1  \ +  {1\over {\beta}} \right) \Bigg\} \  = \ - \  C^+_3  \cdot {1\over { {\bf D}^2 }} \cdot {1\over \alpha} \cdot {{ {\bf q}_{\,2_j}  } \over { \lambda_{1_\tau}}}    \ \ \ \ \\[0.15in]
& \ & \hspace*{-1.5in} \ \ \ \ \ \ \ \ \ \ \ \   \uparrow \ \left(\,\to \ 0 \ \ {\mbox{as}} \ \, {\bf D} \ \to \ \infty\,, \ \,  {\bf D} \ \, {\mbox{takes \  \,  the \ \ form}} \ \   {{ \gamma }\over { \sqrt{\lambda_{1_\tau} \cdot \lambda_{2_\tau}\,}   }}\right).
\end{eqnarray*}
From here we can find $\,\xi_{1_{j_\tau}}\,,\,$ and hence \,$\xi_{2_{j_\tau}}\,.\,$ We observe that
$$
\xi_{1_{j_\tau}} \ = \ O\left( {1\over { {\bf D}  }} \right) \cdot \lambda_{1_\tau}\ = \ o\,(1) \cdot \lambda_{1_\tau}\  \ \ \ {\mbox{and}} \ \ \ \ (\,\xi_{2_{j_\tau}}\ - \ {\bf q}_{\,2_j}\,) \ = \ -\,{{ \xi_{1_{j_\tau}} }\over { \beta }} \ = \ o\,(1) \cdot \lambda_{1_\tau}\,,\leqno (5.14)
$$
[\,$o\,(1) \ \to  \ 0^+$\, as $\,{\bf D}\ \to \ \infty$\,]\,,\, and the solution
$$
 {\bf P}_\tau \ := \ \, (  \lambda_{1_\tau}\,,\    \lambda_{2_\tau} \,;\, \  \xi_{1_{1_\tau}}\,, \ \cdot \cdot \cdot\,, \ \xi_{1_{n_{\,\tau}}}\,; \   \ \xi_{2_{1_\tau}}\,, \ \cdot \cdot \cdot\,, \ \xi_{2_{n_{\,\tau}}}\, ) \ \ \ {\mbox{to}} \ \ {\cal T} \ = \ \vec{\,0} \ \ \ {\mbox{is \ \ unique\,.}} \leqno (5.15)
$$

\vspace*{0.2in}

{\bf \S\,5\,a.}\ \  \it Jacobian matrix.} \ \ In order to compute the degree of $\,{\cal T}\,$ at the image $\,\vec{\,0}\,$,\, (\,in a small\\[0.1in] neighborhood of $\,{\bf P}_\tau \,$)\,,\, let us consider
the Jacobian matrix of the map $\,{\cal T}\,$:\\[0.1in]
{\scriptsize{\[  {\large{ }} \ \ \  \ \ \  \ \ \  \ \ \ \   \ \ \  \ \ \  \ \ \  \ \ \  \ \ \  \ \ \  \ \ \  \ \ \  \ \ \  \ \ \  \left( \begin{array}{cccc}
{{\partial T_{1_o}}\over {\lambda_1}}  & {{\partial T_{1_o}}\over {\lambda_2}} & { 0\  \cdot \cdot \cdot } & {0 \  \cdot \cdot \cdot }\\[0.15in]
{{\partial T_{2_o}}\over {\lambda_1}}  & {{\partial T_{2_o}}\over {\lambda_2}} & { 0\  \cdot \cdot \cdot } & {0 \  \cdot \cdot \cdot }\\[0.15in]
{{\partial T_{1_1}}\over {\lambda_1}}  & {{\partial T_{1_1}}\over {\lambda_2}} & { {{\partial T_{1_1}}\over {\xi_{1_1}}} \cdot \cdot \cdot } & { {{\partial T_{1_1}}\over {\xi_{2_1}}} \cdot \cdot \cdot }\\[0.15in]
\cdot & \\[0.15in]
\cdot & \\[0.15in]
{{\partial T_{2_1}}\over {\lambda_1}}  & {{\partial T_{2_1}}\over {\lambda_2}} & { {{\partial T_{2_1}}\over {\xi_{1_1}}} \cdot \cdot \cdot } & { {{\partial T_{2_1}}\over {\xi_{2_1}}} \cdot \cdot \cdot }\\[0.15in]\cdot & \\[0.15in]
\cdot & \\[0.15in]
\end{array} \right)\ \ \  \ \ \  \ \ \  \ \ \  \ \ \  \ \ \  \ \ \  \ \ \  \ \ \  \ \ \  \ \ \ \ \ \  \ \ \  \ \ \  \ \ \  \ \ \  \ \ \  \ \ \  \ \ \  \ \ \  \ \ \  \ \ \   \]}}

\vspace*{0.1in}

At $\,{\bf P}_\tau\,,\,$
we have [\,using (5.2)\,--\,(5.5)\, and \ $\,{\cal T}\,({\bf P}_\tau) \ = \ {\vec{\,0}}\ $]
\begin{eqnarray*}
{{\partial T_{1_o}}\over {\partial \lambda_1}} \,\bigg\vert_{{\bf P}_\tau}\!\!\!& = &  C_1^+ \cdot   |\, {\varpi}_1\,| \cdot \ell \cdot  ( \lambda_{1_\tau})^{\ell\ -\,1 }    \, - \,    {{n\,-\,2}\over 2} \cdot {1\over { \lambda_{1_\tau}}} \cdot {{ (\sqrt{ \lambda_{1_\tau} \cdot  \lambda_{2_\tau}\,}\ )^{n\,-\,2}
}\over { \gamma^{\,n\,-\,2}  }}\\[0.15in]
& = & \left[ \ \ell \ - \ {{n\,-\,2}\over 2} \ \right] \cdot C_1^+ \cdot   |\, {\varpi}_1\,| \cdot  ( \lambda_{1_\tau})^{\,\ell\ -\,1 } \\[0.1in]
  & \ & \ \ \ \  \ \  \ \ \ \ \ \  \ \  \ \ \ \ \ \  \ \  \ \ \ \ \ \  \ \  \ \ \ \ \ \ \ \  (\,\uparrow \ \ {\mbox{`unit' \ \ for \ \ this \ \ calculation}}\,)\\[0.2in]
{{\partial T_{2_o}}\over {\partial \lambda_2}}\,\bigg\vert_{{\bf P}_\tau}\!\!\!   & = &  C_1^+ \cdot   |\, {\varpi}_2\,| \cdot  \ell \cdot ( \lambda_{2_\tau})^{\ell\ -\,1 }    \, - \,    {{n\,-\,2}\over 2} \cdot {1\over { \lambda_{2_\tau}}} \cdot {{ (\sqrt{ \lambda_{1_\tau} \cdot  \lambda_{2_\tau}\,}\ )^{n\,-\,2}
}\over { \gamma^{\,n\,-\,2}  }} \\[0.15in]
 & = & \left[ \ \ell \ - \ {{n\,-\,2}\over 2} \ \right] \cdot  C_1^+ \cdot   |\, {\varpi}_2\,| \cdot \alpha^{\ell\ -\,1} \cdot ( \lambda_{1_\tau})^{\,\ell\ -\,1 } \ \ \ \ \ \  [\ {\mbox{using \ \ }} \lambda_{2_\tau} \ = \ \alpha \cdot \lambda_{1_\tau}\,]\ ,\\[0.15in]
{{\partial T_{1_o}}\over {\partial \lambda_2}}\,\bigg\vert_{{\bf P}_\tau}\!\!\!   & = &   - \    {{n\,-\,2}\over 2} \cdot {1\over { \lambda_{2_\tau}}} \cdot {{ \left(\sqrt{ \lambda_{1_\tau} \cdot  \lambda_{2_\tau}\,}\ \right)^{n\,-\,2}
}\over { \gamma^{\,n\,-\,2}  }}\ = \ - \    {{n\,-\,2}\over 2} \cdot C_1^+  \cdot   |\, {\varpi}_2\,| \cdot \alpha^{\ell\ -\,1} \cdot  ( \lambda_{1_\tau})^{\,\ell\ -\,1 } \ ,\\[0.15in]
{{\partial T_{2_o}}\over {\partial \lambda_1}}\,\bigg\vert_{{\bf P}_\tau}   & = &   - \    {{n\,-\,2}\over 2} \cdot {1\over { \lambda_{1_\tau}}} \cdot {{ \left(\sqrt{ \lambda_{1_\tau} \cdot  \lambda_{2_\tau}\,}\ \right)^{n\,-\,2}
}\over { \gamma^{\,n\,-\,2}  }}\ = \ - \    {{n\,-\,2}\over 2} \cdot C_1^+  \cdot   |\, {\varpi}_1\,| \cdot    ( \lambda_{1_\tau})^{\ell\ -\,1 } \ ,\\[0.2in]
{{\partial T_{1_j}}\over {\partial \xi_{2_j}}} \,\bigg\vert_{{\bf P}_\tau}\!\!\!& = &   \ o\,(\,[\, \lambda_{1_\tau}\,]^{\,\ell\ -\,1 }\,)\ , \ \ \ \ \  \ \ \ \ \ \ \  {{\partial T_{2_j}}\over {\partial \xi_{1_j}}}\,\bigg\vert_{{\bf P}_\tau}\!\!\! \ =\   \ o\,(\,[\, \lambda_{1_\tau}\,]^{\,\ell\ -\,1 }\,)\ ,\\[0.2in]
{{\partial T_{1_j}}\over {\partial \xi_{1_j} }} \,\bigg\vert_{{\bf P}_\tau}\!\!\! & = &   C_2^+ \cdot  {\varpi}_1 \cdot  ( \lambda_{1_\tau})^{\ell\ -\,1 }  \ + \ o\,(\,[\, \lambda_{1_\tau}\,]^{\,\ell\ -\,1 }\,) \  \ \ \ \ \ \  \  \  \ \ \ \ \ \  \ [\,{\mbox{via}} \ \ (5.4) \ \& \ (5.14)\,]\,,\\[0.2in]
{{\partial T_{2_j}}\over {\partial \xi_{2_j} }} \,\bigg\vert_{{\bf P}_\tau}\!\!\! & = \ &  C_2^+ \cdot  {\varpi}_2  \cdot  ( \lambda_{1_\tau})^{\ell\ -\,1 }  \ + \ o\,(\,[\, \lambda_{1_\tau}\,]^{\,\ell\ -\,1 }\,)\ ,\\[0.2in]
{{\partial T_{1_j}}\over {\partial \lambda_1}}  \,\bigg\vert_{{\bf P}_\tau}\!\!\!& = &   o\,(\,[\, \lambda_{1_\tau}\,]^{\,\ell\ -\,1 }\,)\ ,\ \ \ \
{{\partial T_{1_j}}\over {\partial \lambda_2}}  \  = \   o\,(\,[\, \lambda_{1_\tau}\,]^{\,\ell\ -\,1 }\,)\ ,  \\[0.2in]
{{\partial T_{2_j}}\over {\partial \lambda_1}}  \,\bigg\vert_{{\bf P}_\tau}\!\!\!& = &    o\,(\,[\, \lambda_{1_\tau}\,]^{\,\ell\ -\,1 }\,) \ ,\ \ \ \ {{\partial T_{2_j}}\over {\partial \lambda_2}} \,\bigg\vert_{{\bf P}_\tau}\!\!\! \ = \     o\,(\,[\, \lambda_{1_\tau}\,]^{\,\ell\ -\,1 }\,) \ \ \ \ \ \  \ \ \ \ \ \ (\,j \ = \ 1\,, \ 2\,, \ \cdot \cdot \cdot\, \ n\,)\,.
\end{eqnarray*}
It follows that the Jacobian matrix\,,\,   evaluated at $\,{\bf P}_\tau\,,$\,  can be written as

\vspace*{0.1in}

(5.16)
$\!\!\!\!\!\!\!\!\!\!\!\!\!\!\!\!\!\!\!\!\!\!\!$
{\scriptsize{ \[  \left( \begin{array}{cccc}
\!\!\!\left[ \, \ell \, - \, {{n\,-\ 2}\over 2} \, \right] \cdot C_1^+  \!\cdot\!  |\, {\varpi}_1\,|       & - \,{{n\,-\,2}\over 2} \cdot C_1^+   \!\cdot\!  |\, {\varpi}_2\,|   \!\cdot\! \alpha^{\ell\ -\,1}  & { 0\  \cdot \cdot \cdot } & \!\!{0 \  \cdot \cdot \cdot }\\[0.15in]
- \    {{n\,-\,2}\over 2} \cdot C_1^+  \!\cdot\!  |\, {\varpi}_1\,|      &  \left[ \, \ell \ - \ {{n\,-\,2}\over 2} \ \right] \cdot C_1^+  \!\cdot\!  |\, {\varpi}_2\,|   \!\cdot\!  \alpha^{\ell\ -\,1}   & { 0\  \cdot \cdot \cdot } & {0 \  \cdot \cdot \cdot }\\[0.15in]
o\,(1) & o\,(1) & C_2^+ \!\cdot\!  |\, {\varpi}_1\,|   \!\cdot\!  [\,1 \ + \ o\,(1)\,]& {o\,(1) \cdot \cdot \cdot }\\[0.15in]
\cdot &   \\[0.15in]
\cdot & \\[0.15in]
o\,(1) & o\,(1) \cdot \cdot  \cdot &   \ o\,(1)  \cdot \cdot \cdot    & \!\!\!\!\!\!\!{C_2^+  \cdot\!  |\, {\varpi}_2\,|   \!\cdot\!  [\,1 \ + \ o\,(1)\,] \ \ \  \ \  o\,(1) \cdot \cdot \cdot }\\[0.15in]
\cdot & \\[0.15in]
\cdot & \\[0.15in]
\end{array} \right)\ \ \ \ \ \ \ \ \ \ \] }}

\vspace*{-1.7in}

\hspace*{6.35in}$\scriptsize{\ *\ {\bf I}_{\lambda_{1_\tau}}\, \,.}$

\vspace*{1.5in}

In the above\,,\,  ${\bf I}_{\lambda_{1_\tau}}\,$ is the $\,[\,2\,(n\,+\,1)\,] \times [\,2\,(n\,+\,1)\,]\,$ diagonal matrix with each diagonal entry equal to $\,  ( \lambda_{1_\tau})^{\ell\ -\,1 }  \,$. Focusing on the four terms at the top left hand corner of the matrix in (5.16)\,,\,  the Jacobian determinant at $\,{\bf P}_\tau\,$  is given by

\vspace*{0.1in}

(5.17)

\vspace*{-0.25in}

\begin{eqnarray*}
 {\mbox{Jacobian \ \ Det.}}  \!\! & = &\!\!  \left[\,(\,C_1^+)^2  \cdot  |\, {\varpi}_1\,|  \cdot  |\, {\varpi}_2\,|  \cdot \alpha^{\ell \,-\,1} \cdot \left\{\,\ell \!\cdot\! \left[\,\ell \ - \ 2 \cdot  {{n\,-\,2}\over 2}\,\right]\, \right\}  \,\right] \times \! \\[0.15in]
 & \ &      \times  \left\{  \,  \left[ \,(\,C_2^+)^{2n} \cdot  |\, {\varpi}_1\,|^n  \cdot  |\, {\varpi}_2\,|^n    \,\right] \,\right\}  \cdot [ \ 1 \ + \ o\,(1)\,]  \  *\, [\, ( \lambda_{1_\tau})^{\ell\ -\,1 } \,]^{\,2\,(n\,+\,1)}\ <  \ 0\,,
\end{eqnarray*}
as $\,\ell \ < \ n\,-\,2\,.\,$
Here $\,o\,(1) \ \to \ 0\,$ as $\,{\bf D} \ \to \ 0\,.\,$  Together with (5.15), we conclude that
$$
{\bf{{Deg}}}\,(\,{\cal T}\,,\ B_{\,{\bf P}_\tau} (\,c \cdot  \lambda_{1_\tau})\,, \ {\vec{\,0}}\,) \ = \  -\,1\ \ \ \ \ {\mbox{when}} \ \ {\bf D}_\tau \ = \ {{\gamma }\over {  \lambda_\tau}} \ \ge \ {\bar{D}}_6\ ,  \leqno (5.18)
$$
$$\!\!\!\!\!\!\!\!\!\!\!\!\!\!\!\!\!\!\!\!\!\!\!\!\!\!\!\!\!\!\!\!\!\!\!\!\!\!\!\!\!\!\!\!\!\!\!\!\!\!\!\!\!\!\!\!\!\!\!\!\!\!\!\!\!\!
\lambda_{\,\tau} \ = \ \sqrt{ \lambda_{1_\tau} \cdot \lambda_{2_\tau} \,}\ \,. \leqno (5.19) \ \ \ \ \  \ \ \ \ \ \ {\mbox{where}}
$$
With regard to condition (3.39)\,,\, using the regularity of the map $\,{\cal T}\,$\, and (5.11)\,,\, one deduces that\,,\, for $\,\mu\,$ to be small enough and $\,{\bf P} \ \in (\R^+ \times \R^+) \times (\R^n \times \R^n)\ ,$
\begin{eqnarray*}  (5.20)   \ \ \ \ \  \ \ \ \ \ \    & \ &  \Vert \,{\bf P}\, \ - \ {\bf P}_\tau\,\Vert \  = \    \mu \cdot  \lambda_{1_\tau} \ \ \ \ \   [\  B_{\,{\bf P}_\tau} (\,\mu \cdot  \lambda_{1_\tau}) \, \subset\, (\R^+ \times \R^+) \times (\R^n \times \R^n)  \,] \\[0.15in]
\Longrightarrow \ & \ &   \Vert \, {\cal T}\,( \,{\bf P}) \Vert \   = \  \Vert \, {\cal T}\,(\,{\bf P}) \ - {\cal T}\,(\,{\bf P}_\tau) \,\Vert \ \ \ \ \ \ \ \ \ \ \ \ \ \ \ \ \ \ \ \ \ \ \ \ [\ {\mbox{as}} \ \ {\cal T}\,(\,{\bf P}_\tau) \ = \ 0 \,]  \\[0.15in]
 & \  & \ \ \ \  \ \ \ \ \ \ \ \ \ \ \ \    \ge \  {\tilde c}'\cdot (\lambda_{1_\tau})^{\ell\,-\,1}  \cdot ( \mu \cdot  \lambda_{1_\tau}  )\\[0.15in]
 &  \  & \ \ \ \ \ \ \ \ \ \ \ \ \ \ \ \ \ge \ {\tilde c}'' \cdot (\lambda_{\tau})^{\ell}  \ \ge \ {\tilde c}  \cdot (\lambda_{1_\tau})^{\ell} \ \ \  \left[\,
 = \, O\left({1\over { {\bf D}^{\,n\,-\,2} }} \right)\,; \ \ (2.14) \ \ {\mbox{is \ \ used}}\,\right]\,,  \ \ \ \ \  \ \ \ \  \ \ \
\end{eqnarray*}

\vspace*{-0.15in}

$$\!\!\!\!\!\!\!\!\!\!\!\!\!\!\!\!\!\!\!\!\!\!\!\!\!\!\!\!\!\!\!\!\!\!\!\!\!\!\!\!\!\!\!\!\!\!\!\!\!\!\!\!\!\!\!\!\!\!\!\!\!\!\!\!\!\!
\lambda_{\,\tau} \ = \ \sqrt{ \lambda_{1_\tau} \cdot \lambda_{2_\tau} \,}\ \,.   \ \ \ \ \  \ \ \ \ \ \ {\mbox{where}}
$$
Here $\,\mu\,$ and
$\ {\tilde c}'\,$ are independent on $\,{\bf P}\,$ as long as the conditions in (5.18) \ \& \ (5.20) are satisfied.
Note that in $\,B_{\,{\bf P}_\tau} (\,\mu \cdot  \lambda_{1_\tau})\,,\,$ via (5.4) and (5.9)\,,\, condition (3.45) is satisfied\,.\, We summarize the discussion in this subsection in the following. \\[0.2in]
{\bf Lemma 5.21.} \ \ {\it The map\,} $\,{\cal T}\,: (\R^+ \times \R^+) \times (\,\R^n \times \R^n) \ \to \ \R^{2\,(n\,+\,1)}$,\, {\it as defined in\,} (5.1)\,--\,(5.5)\,,\, {\it has an unique point\,} $\,{\bf P}_\tau \, \in \, (\R^+ \times \R^+) \times (\,\R^n \times \R^n)\,$  [\,{\it  given via\,} (5.15)\,] {\it so that}
$$
{\cal T}\,({\bf P}_\tau\,) \ = \ {\vec{\,0}}\,.
$$
{\it In addition, for a small enough positive number $\,\mu\,$,\, we have } [\,{\it refer to\,} (5.8)\,]
$$
{\bf{{Deg}}} \,(\,{\cal T}\,,\ B_{\,{\bf P}_\tau} (\,\mu \cdot  \lambda_{1_\tau})\,, \ {\vec{\,0}}\,) \ = \  -\,1\,,
$$

\vspace*{-0.3in}

$$
\min_{ \partial B_{ {\bf P}_\tau} (\,\mu \cdot  \lambda_{1_\tau})} \{ \ \Vert \,{\cal T} \,\Vert \,\} \ \ge \ {\tilde c} \cdot (\lambda)^\ell \ . \leqno {\it{and}}
$$
{\it Moreover, condition\,} (2.12) {\it and\,} (3.45) {\it are satisfied by all\,}  $\,{\bf P} \,=\, (\lambda_1\,,\,\lambda_2\,; \ \xi_1\,,\,\xi_2\,) \, \in  B_{\,{\bf P}_\tau} (\,\mu \cdot  \lambda_{1_\tau})\,.$

\newpage

\S \,{\bf 5\,b.} \ \ {\it  Proof of the main theorem.} \ \  In view of (4.2) and (4.3)\,,\, we consider
$$\ \ \ \  \ \ \ \ \ \ \ \  \ \ \ \ \ \ \ \  \ \ \ \ \ \ \ \ \
{\bar C}_b^+ \cdot {\cal T} \ \ \ \  \ \ \ \ \ \ \ \  \ \ \ \ \left[\ {\mbox{recall \ \ that \ \ }} {\bar C}_b^+\ = \  \omega_n \cdot {{(n\,-\,2)^{\,2}}\over {2 }}  \ , \ \ {\mbox{cf. \ \ (4.2)}}\,\right]\,.
 $$
Accordingly, denote the terms within $\,\{ \cdot \cdot \cdot \ \}\,$ [\,referring to (4.2) and (4.3)\,]\, by
 $$
{\bar C}_b^+ \cdot {\cal T}_{o\,(1)} \ \{\ = \ {\bar C}_b^+ \cdot {\cal T} \cdot [\ 1 \ + \ o\,(1)\ ]\,\}\,. \leqno (5.22)
 $$

From (5.1)\,--\,(5.5)\,  and $\,{\cal T}\,({\bf P}_\tau) \ = \ {\vec{\,0}}\,,\,$ we obtain\\[0.1in]
 (5.23)
 $$
\bigg\Vert \   {\bar C}_b^+  \cdot \left[\,{\cal T}_{o\,(1)}\,({\bf P})  \ - \ {\cal T}\,({\bf P}) \,\right]  \bigg\Vert \ \le \ o\,(1) \cdot \lambda_\tau^\ell  \ \ \left[\ \approx \ o\,(1) \cdot {1\over { {\bf D}^{n\,-\,2} }}  \  \right]  \ \ \ \ \mfor \ \ {\bf P} \, \in  \ \overline{B_{\,{\bf P}_\tau} (\,\mu \cdot  \lambda_{1_\tau})}\ \ .
$$
Via (4.2) and (4.3)\,, we have
\begin{eqnarray*}
(5.24) \ \ \ \ \ \ & \ & \ \ \ \ \Vert \ (\lambda \cdot \btd)\,{\bf I}_{{\cal R}_H}\,({\bf z}_\sigma) \ -  \ {\bar C}_b^+  \cdot {\cal T}_{o\,(1)}\ \Vert \ \le \  o\,(1) \cdot {1\over { {\bf D}^{n\,-\,2} }} \ \ \left[\ \approx \ o\,(1) \cdot (\lambda_{\tau})^{\ell}  \  \right] \ \ \ \ \ \ \  \ \ \ \  \ \ \ \\[0.1in]
\Longrightarrow & \ & \Vert \ (\lambda \cdot \btd) \,{\bf I}_{{\cal R}_H}\,({\bf z}_\sigma) \ - \  {\bar C}_b^+  \cdot {\cal T} \,({\bf P}) \,\Vert \ \le \ o\,(1) \cdot {1\over { {\bf D}^{n\,-\,2}  }} \ \approx \ o\,(1) \cdot \lambda_\tau^\ell\\[0.1in]
& \ &  \ \ \ \ \ \ \ \ \  \ \ \ \ \ \ \ \  \ \ \ \ \ \ \ \ \ \ \ \ \ \ \    \ \ \ \ \ \ \  \ \ \ \ \ \ \ \ \ \ \ \ \ \ \ \ \ \ \ \ \ \ \ \  \mfor \ \ {\bf P} \,\in \  \overline{B_{\,{\bf P}_\tau} (\,\mu \cdot  \lambda_{1_\tau})}\ \ .
\end{eqnarray*}
In (5.23) and (5.24)\,,\, $\,o\,(1)\,\to \ 0^+\,$ as $\displaystyle{\,{\bf D}_\tau  \ := \ {\gamma \over \lambda_\tau} \ \to \ \infty\,,\,}$ and $\,{\bf z}_\sigma$\, corresponds to $\,{\bf P}\,$ via
$$
{\bf z}_\sigma \ = \ V_{\lambda_1\,,\, \xi_1} \, + \ V_{\lambda_2\,,\, \xi_2} \ \ \longleftrightarrow \ \ {\bf P} \ = \ (\,\lambda_1\,, \ \lambda_2\,; \ \xi_1\,,\ \xi_2\,)\,.
$$
Thus there exists a positive constant $\,{\bar D}_6 \ (\,\ge {\bar D}_5)\,$,\,  which depends on   $\,n\,$,\,   $\,\ell\,,\,$ $\,{\bar C}_b\,$,\,  $\hbar\,,\,$ ${C}_{R_1}$\,,\, ${C}_{R_2}$\,,\, ${C}_{P_1}$\,, ${C}_{P_2}$\,,\,  and $\,C_p\,,$\, such that if
$$
{\bf D}_\tau \   \ge \ {\bar D}_7\,,\, \leqno (5.25)
$$
then the degree properties [\,as shown in Theorem 3.38\,]\,,\,  (5.18)\,,\, (5.19)  and  (5.24) imply that

\vspace*{-0.1in}

$$
{\bf{Deg}} \left(\ (\lambda \cdot \btd) \, {\bf I}_{{\cal R}_H}\,,\ B_{\,{\bf P}_\tau} (\,\mu \cdot  \lambda_{1_\tau})\,, \ {\vec{0}}\,\right) \ = \   {\bf{Deg}}\,\left({\cal T}\,,\ B_{\,{\bf P}_\tau} (\,\mu \cdot  \lambda_{1_\tau})\,, \ {\vec{0}}\,\right) \ = \ -\,1 \ (\,\not= \  0\,)\,.\leqno (5.26)
$$

\vspace*{0.1in}

Hence the reduced functional has a critical point at $\,{\,{\bf P}}_R\,=\,(\,{\bar\lambda}_1\,, \ {\bar\lambda}_2\,; \ {\bar \xi}_1\,,\ {\bar\xi}_2\,) \, \in \, B_{\,{\bf P}_\tau} (\,\mu \cdot  \lambda_{1_\tau})\,.\,$
 Via Lemma 3.32\,,\, equation (1.2) has a solution of the form $\,{\bar{\bf z}}_\sigma \, +\,  w_{{\bar{\bf z}}_\sigma}$\,,\, where
 $$
 {\bar{\bf z}}_\sigma \ = \ V_{{\bar\lambda}_1\,, \  {\bar \xi}_1  } \ + \ V_{{\bar\lambda}_2\,, \  {\bar \xi}_2  } \ .
 $$

\vspace*{0.1in}

\S \,{\bf 5\,c.} \ \ {\it  Estimate on $\,\gamma_o\,.$} \ \   The following consideration provides an idea on the size of $\,\gamma_o\,$ [\,cf. Remark (1) in {\bf \S\,1\,a}\,]\,.\, From condition (5.25)
\begin{eqnarray*}
& \ & \!\!\!\!\!\!\!\!\!\!{\bf D}_\tau \ = \  {{\gamma}\over {\sqrt{ \lambda_{1_\tau} \cdot \lambda_{2_\tau}  }  }} \ \ \ge \ {\bar D}_7 \ \Longleftrightarrow\   \left( {{\gamma}\over {\sqrt{ \lambda_{1_\tau} \cdot \lambda_{2_\tau}  }  }} \right)^{\!\!n\,-\,2}
 \ \ \ge \ {\bar D}_7^{\,n\,-\,2}\\[0.15in]
 &  \Longleftrightarrow  &    C^+_1 \cdot {1\over { |\, {\varpi}_1\,|  }} \cdot {1\over {  \lambda_{1_\tau}^\ell } }
 \ \ \ge \ {\bar D}_7^{\,n\,-\,2} \ \ \ \ \ \ \ \ \ \ \ \ \ \ \ \ \ \ [\,{\mbox{via \ \ (5.2) \ \ with}} \ \ {\cal T} \,({\bf P}_\tau) \ = \ \vec{\,0}\ ]\\[0.15in]
 & \Longleftrightarrow &   C\,\cdot{1\over {|\, {\varpi}_1\,| \cdot  |\, {\varpi}_1\,|^{ {{   \ell }\over { (n\,-\,2) \ - \ \ell }}  }  }} \cdot {1\over {  \gamma^{ {{ \ell \cdot (n\,-\,2) }\over { (n\,-\,2) \ - \ \ell }}  } } }
 \ \ \ge \ {\bar D}_7^{\,n\,-\,2}\ \ \ \ \ \ \ \ \ \ \ \ \ \ \ \ \   [\,{\mbox{using\ \ (5.8) }}\ ]\\[0.15in]
  & \Longleftrightarrow &  C\,\cdot{1\over { |\, {\varpi}_1\,|^{ {{ n\,-\,2}\over { (n\,-\,2) \ - \ \ell }}  }  }} \cdot {1\over { {\bar D}_7^{\,n\,-\,2}    } }
 \ \ \ge \ \gamma^{ {{ \ell \cdot (n\,-\,2) }\over { (n\,-\,2) \ - \ \ell }}  }\\[0.15in]
  & \Longleftrightarrow &  C\,\cdot {1\over { {\bar D}_7^{{ (\,n\,-\,2) \ - \ \ell}\over {   \ell }} } } \cdot {1\over { |\, {\varpi}_1\,|^{1\over {\ell}} }}
 \ \ \ge \ \gamma \ \ \ \ \ \ \ \ \ \ \ \ \ \ \ \
  \ \ \ \ \ \ \ \ \ \ \ \ \ \ \ \ \ \ \ \ \ \ \ \ \ \ [\,{\mbox{cf.}} \ \ (1.15)\,]\ . \\[0.15in]
\end{eqnarray*}

\vspace*{0.2in}

{\bf \Large {\bf  Appendix.}}

\vspace*{0.2in}

{\bf Lemma A.1.} \ \ {\it For positive numbers \,$\,a\,$, \,$\,b\,$, \,$\,\tau\,$, \,$\,\beta\,$, and \,$\,M\,$, we have the following.}
$$
(\,0 \ <\ \,)\,\ \beta\ \le\ 2\ \ \Longrightarrow\ \ \Big|\,(\,a\,+\, b\,)^{\,\beta}\ -\ (\,a^\beta\,+\,b^{\,\beta}\,)\,\Big|\ \le\ {\bar C}_\beta\cdot (\,a\cdot b\,)^{\beta\over 2}\ . \leqno (A.2)
$$
{\it Here $\,{\bar C}_\beta\,$ is a positive constant depending on $\,\beta\,$, but not on $\,a$ and $\,b\,$.}
$$
M\ \ge\ 2\ \ \Longrightarrow\ \ \Big|\,(\,a \,+\,b\,)^{\,M}\ -\ (\,a^M\,+\,b^{\,M}\,)\,\Big|\ \le\ {\bar C}_M \cdot \left(\,a^{M-1} \cdot b\ +\ b^{M-1}\cdot a\,\right)\,. \leqno (A.3)
$$
{\it Here $\,{\bar C}_M\,$ is a positive constant depending on $\,M\,$, but not on $\,a$ and $\,b\,$. Moreover\,,}
$$
(\,0\ <\ \,)\,\ \tau\ \le\ 1\ \ \Longrightarrow\ \ (\,a\,+\, b\,)^{\,\tau}\ \le\ a^\tau\ +\ b^\tau\ . \leqno (A.4)
$$

\vspace*{0.1in}

{\bf Lemma A.5.} \ \ {\it For \,$\,n\ \ge\ 6\,$, let \,$\,P\,$ and \,$\,Q\,$ be positive numbers with}
$$
P\ \,+\,\ Q\ \,=\,\ n\ .
$$
{\it For \,}$\,(\,\lambda_1\,,\ \lambda_2\,;\ \xi_1\,,\ \xi_2\,)\ \in \ (\R^+ \times \R^+) \times (\R^n \times \R^n)\,$,\, {\it let}
$$
\mathcal{S}\,(\lambda_1\,,\ \lambda_2\,;\ \xi_1\,,\ \xi_2\,;\ P\,,\ Q)\ :=\ \int_{\R^n} \left({{\lambda_1}\over {\lambda_1^2\,+\,|\,y\,-\,\xi_1\,|^2}}\right)^{\!\!P} \cdot\, \left({{\lambda_2}\over {\lambda_2^2\,+\,|\,y\,-\,\xi_2\,|^2}}\right)^{\!\!Q}\ dy\ .
$$
{\it Assume that}
$$
 {\bar C}^{\,-1}  \cdot \lambda_1\ \le\ \lambda_2\ \le\ {\bar C} \cdot \lambda_1\ ,\ \ \ \ \ \ |\,\xi_1\,| \ \le \ \lambda \ \ \ \ {\it{and}} \ \ \ \  |\ \xi_2 \ - \ {\bf q}_2\,| \ \le \ \lambda\,. \leqno (A.6)
 $$
{\it There is a positive number \,$\bar{D}_1$ such that for\ \,}$\,(\,\lambda_1\,,\ \lambda_2\,;\ \xi_1\,,\ \xi_2\,)\ \in \ (\R^+ \times \R^+) \times (\R^n \times \R^n)\,$ {\it satisfying} (A.6)
$$
 {\it if}\ \ \ \ {\bf  D}\ \ge\ \bar{D}_1 \ , \leqno (A.7)
$$

\vspace*{-0.25in}
$$
{\it then}\ \ \ \ \mathcal{S}\,(\lambda_1\,,\ \xi_1\,;\ \lambda_2\,,\ \xi_2\,;\ P\,,\ Q)\ \leq\ \left\{
\begin{array}{l@{\ \ \ \ }l}
\displaystyle{{\bar{C}_3}\over{\bf D}^{2\,P}} & \displaystyle{\it if\ \ }P\ <\ {n\over2}\ ,\\[0.15in]
\displaystyle{{\bar{C}_3}\over{\bf D}^{n}}\cdot\ln\,{\bf D} & \displaystyle{\it if\ \ }P\ =\ Q\ =\ {n\over2}\ ,\\[0.15in]
\displaystyle{{\bar{C}_3}\over{\bf D}^{2\,Q}} & \displaystyle{\it if\ \ }P\ >\ {n\over2}\ .
\end{array}
\right. \leqno (A.8)
$$
{\it Here \,$\,\bar{C}_3\,$ depends on \,$\,n\,$, \,$\,{\bar C}\,$, and \,$\,Q\,$. In particular, $\,\bar{D}_1\,$ and $\ \bar{C}_3\,$ do  not depend on \,$\,(\,\lambda_1\,,\ \lambda_2\,;\ \xi_1\,,\ \xi_2\,)\,$ as long as\ }(A.6) {\it is satisfied.}

\vspace*{0.2in}

\hspace*{0.5in}For the proofs of Lemmas A.1 and A.5\,,\, we refer to  \cite{III}\,.

\vspace*{0.2in}

{\bf Proof of Lemma 3.12.} \ \ Let $\,V_1\ = \ V_{\lambda_1\,,\ \xi_1}\,$ and $\,V_2\ = \ V_{\lambda_2\,,\ \xi_2}\ $.\, Via Lemma A.5\,,\, we have
$$
\bigg\vert \  (V_1 \ + \ V_2)^{{2n}\over {n\,-\,2}} \ - \ \left( V_1^{{2n}\over {n\,-\,2}}
\ + \  V_2^{{2n}\over {n\,-\,2}} \right) \, \bigg\vert\  \le \   C\cdot \left( V_1^{{n\,+\,2}\over {n\,-\,2}} \cdot V_2 \ + \ V_2^{{n\,+\,2}\over {n\,-\,2}} \cdot V_1 \right) \ . \leqno (A.9)
$$

\vspace*{-0.12in}

As
$$\,|\,\xi_1\,| \ \le \ {\bar c} \cdot \lambda_1\,,$$

\vspace*{-0.2in}

\begin{eqnarray*}
 (A.10) \ \ \ {{  1}\over {(\,\lambda^2_1
\,+\, |\, y \,-\, \xi_1|^{\,2}\,)^{\,n}  }} & = & {{  1}\over {(\,\lambda^2
\,+\, |\, y \,|^{\,2}\,)^{\,n}  }} \cdot \left\{ 1\ - \ n \cdot {{ |\,\xi_1|^2 \ - \ 2 \,y \cdot \xi_1 }\over {\lambda^2
\,+\, |\, y \,|^{\,2}    }} \ + \ \cdot \cdot \cdot \ \right\}\\[0.2in]
& \Longrightarrow & \ \ \left( {{\lambda_1 }\over {\lambda^2_1 \ + \ |\, y \ - \ \xi_1\,|^2   }}\right)^{\!\!n} \ = \ \left( {{\lambda_1 }\over {\lambda^2_1 \ + \ |\, y \,|^2   }}\right)^{\!\!n}
\cdot [\,1 \ + \ o\,(1)\,]\ . \ \ \ \
\end{eqnarray*}
Here $\,o\,(1) \ \to \ 0^+\,$ as $\ {\bar c} \ \to \ 0^+\,.\,$ Consider the partition
$$
\R^n \ = \ B_{o} ({\rho}) \ \cup \ B_{{\bf q}_2} ({\rho}) \ \cup \ \{\,\R^n \,\setminus\,  [\,B_{o} ({\rho}) \,\cup\,B_{{\bf q}_2} ({\rho})\,]\, \}\ . \leqno (A.11)
$$
Here $\,\rho\,$ is the same number as in $\,{\mbox{(1.9)}}_{(i)}\,.\,$  The leading terms in the estimate in Lemma 3.12 can be obtained by integral of the form\\[0.1in]
(A.12)

\vspace*{-0.3in}

$$
 \int_0^\rho \,r^{\,L} \cdot \left( {{ \lambda }\over {  \lambda^2 \ + \ r^2}} \right)^n \cdot r^{n\,-\,1} \cdot dr  \ \leq\ \left\{
\begin{array}{l@{\ \ \ \ }l}
\displaystyle{O\,(\lambda^{L}\,) } & \displaystyle{\mbox{if}\ \ } \ \ \  L \ <\ n\ ,\\[0.11in]
\displaystyle{O\,(\lambda^{\,n\,-\,o\,(1)}\,)  } & \displaystyle{\mbox{if}\ \ } \ \ \  L \ = \ n \ \ \ \ \ \left( \, {\bf D} \ \ge \ {\bar D}_3 \right)\,,\\[0.11in]
\displaystyle{  O\,(\lambda^n\,) } & \displaystyle{\mbox{if}\ \ } \ \ \  L \  > \ n\ .\\[0.11in]
\end{array}
\right.
$$
Cf. (A.14) for the contribution outside $\,B_{o} ({\rho}) \ \cup \ B_{{\bf q}_2} ({\rho}) \,.\,$
Observe that
$$\ \ \ \ \ \  \
{\rho\over \lambda} \ =\  \hbar \cdot {\gamma\over \lambda} \ = \ \hbar \cdot {\bf D}\ \ \ \ \ \  \ \ \ \ \ \ \  \ \ \ \ \ \ \  \  [\,{\mbox{cf}}. \ \ {\mbox{(1.9)}}_{(i)}\,]\,.
$$
The ``\,cross over\," terms are estimated as in {\bf \S\,2\,b} [\,cf. (2,21)\,]\,,\, giving rise to
 $$
 \int_{B_o(\rho)} V_2^{{2n}\over {n\,-\,2}} \ = \ O\left( {1\over { {\bf D}^{n} }} \right)\,.
 $$
Next, we consider  the mixed terms that come from the right hand side of (A.9)\,,\, again making use of (2.21)\,.

\vspace*{-0.3in}

\begin{eqnarray*}
(A.13)& \ & \int_{\R^n} |\,H\,(y)\,|^m \cdot \left( {{ \lambda_1}\over {\lambda_1^2  \ + \ |\,y\ - \ \xi_1|^2  }}\right)^{{n\,+\,2}\over 2}
\cdot \left( {{ \lambda_2}\over {\lambda_2^2  \ + \ |\,y\ - \ \xi_2|^2  }}\right)^{{n\,-\,2}\over 2} \\[0.1in]
& \le &  C_1\,(n) \cdot {1\over {\lambda_1^{{n\,-\,2}\over 2} }} \cdot {1\over {{\bf D}^{n\,-\,2} }} \,\cdot  \, [\,1 \ + \ o\,(1)\,] \cdot \int_{B_{o} (\rho)} r^{\,m \cdot \ell} \cdot \left( {{ \lambda_1}\over {\lambda_1^2  \ + \ |\,y\ - \ \xi_1|^2  }}\right)^{{n\,+\,2}\over 2} \  dy   \ \,+ \\[0.1in]  & \ &  \ \ \ + \  C_2 (n) \cdot \!\!{1\over {\lambda_2^{{n\,-\,2}\over 2} }} \!\cdot\! {1\over {{\bf D_o}^{n\,+\,2} }} \cdot  \!  [\,1 \, + \, o\,(1)\,] \cdot \int_{B_{{\bf q}_2} (\rho)} r^{\,m \cdot \ell} \cdot \left( {{ \lambda_2}\over {\lambda_2^2  \ + \ |\,y \ - \ \xi_2\,|^2  }}\right)^{{n\,-\,2}\over 2} \ dy \ \\[0.1in]
& \ & \  \ \ \  \ \ \ \ + \ O \left( {1\over { {\bf D}^{\,n\,-\,o\,(1)} }} \right)   \ \   \left( \leftarrow \ \ \, {\mbox{see \ \ below\,}}\,; \ \  \ \ \  {\mbox{shift \ \ center \ \ to   \ \ }} {\bf q}_2 \ \ \uparrow \right) \ \ \ \ \ \ \ \ \ \ \ \ \ \ \ \ \ \ \ \ \ \  \ \ \ \ \ \ \ \ \ \ \  \left(  \uparrow \right) \\[0.11in]
& \le& C_3\,(n) \cdot \lambda^2 \cdot {1\over {{\bf D_o}^{n\,-\,2} }} \,\cdot  \, [\,1 \ + \ o\,(1)\,] \ + \ O \left( {1\over { {\bf D}^{\,n\,-\,o\,(1)} }} \right) \ .
\end{eqnarray*}
In the above we apply (1.10)\,,\, (1.11) and  the following.
\\[0.1in]
(A.14)
\vspace*{-0.3in}
\begin{eqnarray*}
& \ & \bigg\vert\ \int_{\R^n \setminus\,[\,B_{o}\,({\rho }) \ \cup \ B_{{\bf q}_2}\,({\rho })\,]}  \ |\,H\,|^{\,m} \cdot  V_1^{{n\,+\,2}\over {n\,-\,2}} \cdot V_2\\[0.1in]
& \le & C_1^m \cdot \int_{\big\{\,\R^n \,\setminus\ [\,B_{o}\,({\rho }) \ \cup \ B_{{\bf q}_2}\,({\rho })\,]\,\big\} \ \cap \ \{\,V_1 \ > \ V_2\,\}} V_1^{{n\,+\,2}\over {n\,-\,2}} \cdot V_2 \ + \  \\[0.1in]
 & \ & \ \ \ \ \ \  \  \ \ \ \   \ \ \ \ \ \  \  \ \ \ \  \  \ \ \ \  + \  C_1^m \cdot  \int_{\big\{\,\R^n \,\setminus\ [\,B_{o}\,({\rho }) \ \cup \ B_{{\bf q}_2}\,({\rho })\,]\,\big\} \ \cap \ \{\,V_2 \ > \ V_1\,\}}   V_1^{{n\,+\,2}\over {n\,-\,2}} \cdot V_2 \\[0.1in]
& \le &   C_1^m \cdot C\,(n) \cdot   \int_{\rho }^\infty \left( {\lambda\over {\lambda^2 \ + \ r^2}} \right)^{\!\!n}\,r^{n\,-\,1} \cdot dr\ \le \    C_1^m \cdot C'\,(n) \cdot {1\over { {\bf D}^{\,n\,-\ o\,(1)}  }} \ .
\end{eqnarray*}
Similarly we estimate the corresponding integral for $\,V_2^{{n\,+\,2}\over {n\,-\,2}} \cdot V_1\,$.   {{\hfill {$\rlap{$\sqcap$}\sqcup$}}

\newpage

 {\bf e\,-Appendix} \, starts at pp. 36.

\newpage

\centerline{\bf \LARGE {\bf e\,-Appendix to the Article  }}

\vspace*{0.18in}

\begin{center} {\LARGE {\bf ``\,Conformal Scalar Curvature Equation on $S^n$      }} \medskip \medskip   \\
{\LARGE {\bf   -   Functions  With Two Close Critical Points } } \medskip \medskip  \medskip  \medskip   \\
{\LARGE {\bf     (Pseudo\, Twin\,-\,Peaks)}\," }\medskip \medskip \smallskip \\

\vspace{0.2in}

\centerline{\Large  {Man Chun  {\LARGE L}EUNG${\,} $\ \ \ \&\ \ \ Feng \,{\LARGE Z}HOU${\,} $ }}

\vspace*{0.2in}

{\large {National University of Singapore}} \\[0.1in]
\end{center}


\vspace*{0.68in}

{\it In this {\bf e}\,-Appendix we follow the notations, conventions, equation numbers, section numbers, lemma, proposition and theorem numbers as  used in the main article\,}  \cite{Twin}\,,\,  {\it unless otherwise is specifically mentioned} (\,{\it  for instances, those equation numbers starting with}\, `A'\,).  \\[0.4in]
%
%
{\bf \large \S\,A.1.} \ \ {\bf \large Estimates (2.23)   \,\&\, (2.25).}\\[0.2in]
See also Lemma B.2 and Lemma B.4 in \cite{Yan+} and formulas  2.119 and 2.206 in \cite{Bahri}\,.\,
For
$$
{\bf z}_\sigma \ = \ V_{\lambda_1\,, \ \xi_1} \ + \ V_{\lambda_2\,, \ \xi_2}\,,
$$
where [\,cf. (2.12)\,]
$$
{\bar C }^{-1}  \cdot \lambda_2 \ < \ \lambda_1 \ <\  {\bar C} \cdot \lambda_2\,, \ \ \ \ |\,\xi_1\,| \ <  \   {\bar c} \,\cdot \lambda \ \  \ \ \  {\mbox{and}} \  \ \ \ \ \ |\,\xi_2 \ - \ {\bf q}_{\,2}\,| \ <  \  {\bar c} \,\cdot \lambda\,, \leqno (A.1.1)
$$
let
$$
\rho_\mu \ =  \ \mu \cdot |\,\xi_1 \ - \ \xi_2\,| \ . \leqno (A.1.2)
$$
Here $\,\mu\,$ is a chosen small positive number so that
$$
\mu \ \to \ 0^+\ \ ({\mbox{slowly}}) \ \ {\mbox{and}} \ \ \ \ \mu^M \cdot {\bf D} \ \to \ \infty \ \ \ \ {\mbox{when}} \ \ \ \ {\bf D} \ \to \ \infty\,.\leqno (A.1.3)
$$
Here $\,M\,$ is a (fixed) large integer\,.\,
See (A.1.12), (A.1.19) and (A.1.28)\,.\, For most particular purpose one can take
$$
\mu \ = \ {1\over 2} \cdot {1\over { {\bf D}^{\,\epsilon} }} \ \ \ \ \ \ {\mbox{for}} \ \ \ \ {\bf D} \ \gg \ 1\,,\leqno (A.1.4)
$$
where $\,\epsilon\ < \ 1\,$ is any fixed small positive number\,.\,

\newpage

\hspace*{0.5in} Consider the partition
 $$
 \R^n \ = \ B_{\xi_1} (\,\rho_\mu) \ \cup \ B_{\xi_2} (\,\rho_\mu) \ \cup \ \{\,\R^n \setminus [\, B_{\xi_1} (\,\rho_\mu) \ \cup \ B_{\xi_2} (\,\rho_\mu)  \,]\,\}\,.\leqno (A.1.5)
 $$
 We compare $\,V_1\,$ and $\,V_2\,$ on $\, B_{\xi_1} (\,\rho_\mu)\,.\,$ In order to do so,
 we make use of the following inequalities ($A$ and $B$ are positive numbers)
 $$
 \left({1\over {1 \ + \ t}}\right)^{\!\!\!{{n\,-\,2}\over 2}} \ = \ 1 \  - \ {{n\,-\,2}\over 2} \cdot t \cdot [\,1 \ + \ o\,(1)\,] \ \ \ \ \mfor t \ >  \ 0 \ \ {\mbox{small}}\,, \leqno (A.1.6)
 $$
$$
(A \ + \ B)^{{n\,+\,2}\over {n \,-\,2}} \ = \ A^{{n\,+\,2}\over {n \,-\,2}} \ + \ {{n\,+\,2}\over {n \,-\,2}}\cdot
A^{4\over {n \,-\,2}} \cdot B \ + \ O\,(1) \cdot B^{{n\,+\,2}\over {n \,-\,2}} \ \ \ \  \ \ \ \mbox{for}  \ \ \ \ {B\over A} \ \
{\mbox{small}}\,,\leqno (A.1.7)
$$
\hspace*{4.7in}where $\,o\,(1) \ \to \ 0\ $ as $\ t \ \to  \ 0\,.\,$\\[0.05in]
It follows from (A.1.7) that\\[0.1in]
(A.1.7\,$'$)
$$
(A \ + \ B)^{{n\,+\,2}\over {n \,-\,2}} \ - \ \left( A^{{n\,+\,2}\over {n \,-\,2}} \ + \ B^{{n\,+\,2}\over {n \,-\,2}} \right)
\ =  {{n\,+\,2}\over {n \,-\,2}}\cdot
A^{4\over {n \,-\,2}} \cdot B \ + \ O\,(1) \cdot B^{{n\,+\,2}\over {n \,-\,2}} \ \ \ \  \mbox{for}  \ \ \ \ {B\over A} \ \
{\mbox{small}}\,.
$$
In $\,B_{\xi_1} (\,\rho_\mu)\,,$\, we have\\[0.1in]
(A.1.8)
 \begin{eqnarray*}
 V_{\lambda_2\,,\ \xi_2}(y)   & = & \left( {{\lambda_2}\over { \lambda^2_2 \ + \ |\,y\ - \ \xi_2\,|^2  }}\right)^{\!\!{{n\,-\,2}\over 2}  }\\[0.125in]
& = &  \left(\ {{{1\over\lambda_1}}\over { \left( {{\lambda_2}\over {\lambda_1}}  \right)\ + \ {{|\,y\ - \ \xi_2\,|^2}\over
{\lambda_1 \cdot \lambda_2}}  }}\right)^{\!\!{{n\,-\,2}\over 2}  } \ = \  {1\over {\lambda_1^{{n\,-\,2}\over 2} }} \cdot \left(\  {{{1}}
\over { \left( {{\lambda_2}\over {\lambda_1}}  \right)\ + \ {{|\,(\,y\ - \ \xi_1)\ + \ (\,\xi_1\ - \ \xi_2)\,|^2}\over
{\lambda_1 \cdot \,\lambda_2}}  }}\right)^{\!\!{{n\,-\,2}\over 2}  } \\[0.125in]
& = &  {1\over {\lambda_1^{{n\,-\,2}\over 2} }} \cdot \left(\ {{{1}}
\over { \left( {{\lambda_2}\over {\lambda_1}}  \right) \ + \ {{|\, \,\xi_1\ - \ \xi_2\,|^2}\over
{\lambda_1 \cdot \,\lambda_2}}   \ + \ {{|\,\,y\ - \ \xi_1\,|^2}\over
{\lambda_1 \cdot \,\lambda_2}}   \ + \ {{\,2\,(\,y\ - \ \xi_1)\cdot(\,\xi_1\ - \ \xi_2)\, }\over
{\lambda_1 \cdot \,\lambda_2}}}}\ \right)^{\!\!{{n\,-\,2}\over 2}  }\\[0.125in]
& = &  {1\over {\lambda_1^{{n\,-\,2}\over 2} }} \cdot {1\over {{\bf d}^{n\,-\,2} }} \cdot \left(\ {{{1}}
\over { 1 \ + \ {1\over {{\bf d}^{2} }} \cdot \left( {{\lambda_2}\over {\lambda_1}}  \right) \ + \ {1\over {{\bf d}^{2} }} \cdot {{|\,\,y\ - \ \xi_1\,|^2}\over
{\lambda_1 \cdot \,\lambda_2}}   \ + \ {1\over {{\bf d}^{2} }} \cdot {{\,2\,(\,y\ - \ \xi_1)\cdot(\,\xi_1\ - \ \xi_2)\, }\over
{\lambda_1 \cdot \,\lambda_2}}}}\ \right)^{\!\!{{n\,-\,2}\over 2}  }\\[0.125in]
& \ & \ \ \ \ \ \ \ \ \ \ \ \ \ \ \ \ \ \ \ \ \ \  \ \ \ \ \ \ \ \   \leftarrow \ \ \ \ \ \ \ \ \ \ \ \  \ \ \ \ \ \ \ = \ {\bf t}\   \ \ (\,\leftarrow \ \ {\mbox{small}}\,) \ \ \ \ \  \ \ \ \ \ \ \ \ \ \ \,      \rightarrow \\[0.125in]
& = & {1\over {\lambda_1^{{n\,-\,2}\over 2} }} \cdot {1\over {{\bf d}^{n\,-\,2} }} \cdot \left\{\, \  1 \   - \ {{n\,-\,2}\over 2} \cdot {\bf t} \cdot [\,1\,+\,o\,(1)\,]\,\right\} \ \ \ \ \ \ \ \ \ \ \ \ \ \ \ \ \ [\,{\mbox{via}} \ \ (A.1.6)\,]\ .
\end{eqnarray*}
Recall that
$$
{\bf d}  \  :=   \  {{|\,\xi_1\ -\ \xi_2\,|}\over{\,\sqrt{\,\lambda_1 \,\cdot\,\lambda_2\,}\,}}\ .
$$

\newpage
In (A.1.8)\,,\, $\,o\,(1) \ \to \ 0\,$ as
$$
{\bf t} \ = \ {1\over {{\bf d}^{2} }} \cdot \left[\ \left( {{\lambda_2}\over {\lambda_1}}  \right) \ + \   {{|\,\,y\ - \ \xi_1\,|^2}\over
{\lambda_1 \cdot \,\lambda_2}}   \ + \   {{\,2\,(\,y\ - \ \xi_1)\cdot(\,\xi_1\ - \ \xi_2)}\over
{\lambda_1 \cdot \,\lambda_2}}\,\right] \ \ \to \ 0\ . \leqno (A.1.9)
$$
Note that
$$
 {1\over {{\bf d}^{2} }} \cdot {{|\,\,y\ - \ \xi_1\,|^2}\over
{\lambda_1 \cdot \,\lambda_2}}   \  \le \    {1\over {{\bf d}^{2} }} \cdot {{\rho^2_\mu}\over
{\lambda_1 \cdot \,\lambda_2}} \ = \ {1\over {{\bf d}^{2} }} \cdot {{|\,\xi_1\ - \ \xi_2\,|^2}\over
{\lambda_1 \cdot \,\lambda_2}} \cdot \mu^2 \ = \ \mu^2\,  \leqno (A.1.10)
$$
and
 \begin{eqnarray*}
 (A.1.11) \ \ \ \ \ \  \ \ \ \    \ \ \ \ \   {1\over {{\bf d}^{2} }} \cdot {{\, (\,y\ - \ \xi_1)\cdot(\,\xi_1\ - \ \xi_2)\, }\over
{\lambda_1 \cdot \,\lambda_2}}
& \le & {1\over {{\bf d}^{2} }} \cdot {{\,|\,y\ - \ \xi_1\,|\, }\over
{\sqrt{\lambda_1 \cdot \,\lambda_2\,}}}  \cdot {{\,|\,\xi_1\ - \ \xi_2\,|\, }\over
{\sqrt{\lambda_1 \cdot \,\lambda_2\,}}} \\[0.2in]
& \le & {1\over {{\bf d}^{2} }} \cdot {{\,\rho_\mu\, }\over
{\sqrt{\lambda_1 \cdot \,\lambda_2\,}}}  \cdot {\bf d}  \ =  \ \mu\,. \ \ \ \ \ \  \ \ \ \  \ \ \ \ \ \  \ \ \ \   \ \ \ \ \ \  \ \ \ \
\end{eqnarray*}
Recall (A.1.3)\,.\, Hence we obtain
$$
V_2 \ = \ {1\over {\lambda_1^{{n\,-\,2}\over 2} }} \cdot {1\over {{\bf D}^{n\,-\,2} }} \cdot [\,1 \ + \ o\,(1)\,]
\ \ \ \ \ \ \ \  {\mbox{in}} \ \ \ \ B_{\,\xi_1} (\rho_\mu) \ \ \   \left( \,{\bf D} \ = \ {{ \gamma}\over { \sqrt{ \,\lambda_1 \cdot \lambda_2\,}  }} \right)\,.\leqno (A.1.12)
$$
Here $\,o\,(1) \ \to \ 0\,$ as $\,{\bf D}  \ \to \ \infty\,$.\,  In the above we apply (A.1.3)\,,\, (2.13) \,and\, (2.14)\,.\,\bk
As in (2.22)\,,\,
$$
{\cal I} \ = \ \left\{\, \left( V_{\lambda_1\,,\ \xi_1}^{{n\,+\,2}\over{n\,-\,2}}\ +\ V_{\lambda_2\,,\ \xi_2}^{{n\,+\,2}\over{n\,-\,2}} \right)  \ \,- \  (\,V_{\lambda_1\,,\ \xi_1} \ + \ V_{\lambda_2\,,\ \xi_2})^{{n\,+\,2}\over{n\,-\,2}} \,\right\}\ \ ( \, < \ 0\,)\ . \leqno (A.1.13)
$$
According to the partition in (A.1.5)\,,\, we break down the integral as
\begin{eqnarray*}
(A.1.14) \ \ \ \ \ \ \ \ \   \int_{\R^n}  {\cal I} \cdot [\,\partial_{\lambda_1} \,V_1\,]  & =  & \int_{B_{\xi_1} ({\rho_\mu})}  {\cal I} \,\cdot\, [\,\partial_{\lambda_1} \,V_1\,]  \ + \  \int_{B_{\xi_2} ({\rho_\mu})}  {\cal I} \,\cdot\,[\,\partial_{\lambda_1} \,V_1\,]  \\[0.15in]
  & \ & \ \ \ \ \ \  \ \ \ \ \ \ \ \ \ \ \ \ \ \ \ \ + \  \int_{\R^n \,\setminus\,  [\,B_{\xi_1} ({\rho_\mu}) \,\cup\,B_{\xi_2} ({\rho_\mu})\,] } {\cal I} \,\cdot\,[\,\partial_{\lambda_1} \,V_1\,]\ .\ \ \ \ \ \ \ \ \ \ \ \
\end{eqnarray*}
A direct calculation shows that
\begin{eqnarray*}
(A.1.15) & \ &
  \!\!\!{{\partial V_{\lambda_1\,, \, \xi_1}}\over {\partial \lambda_1}} \ = \  {{\partial}\over {\partial \lambda_1}} \left[\left({\lambda_1\over {\lambda_1^2 \,+\, |\, y \,-\, \xi_1|^{\,2}}} \right)^{\!\!\!{{n \,-\, 2}\over 2}} \,\right] \\[0.15in]
   & = & \!\!-\,{{n \,-\, 2}\over 2} \cdot \lambda_1^{{n \,-\, 4}\over 2} \!\cdot \! {{ (\,\lambda^2_1 \,-\, |\, y \,-\, \xi_1|^{\,2})}\over {(\,\lambda^2_1 \ + \ |\, y \,-\, \xi_1|^{\,2})^{{n }\over 2} }}\ \ \Longrightarrow \   \bigg\vert \,{{\partial V_{\lambda_1\,, \, \xi_1}}\over {\partial \lambda_1}} \,\bigg\vert \, \le \, {{n\,-\,2}\over 2}\cdot {1\over {\lambda_1}} \cdot V_{\lambda_1\,, \, \xi_1} \ .\ \ \ \ \  \ \ \ \ \ \ \ \ \ \  \ \ \ \ \ \ \ \ \ \  \ \ \ \ \ \ \ \ \ \  \ \ \ \ \ \ \ \ \ \  \ \ \ \ \   \end{eqnarray*}

   \newpage

 In the following, we play careful attention on the ($\,\pm\,$) sign. Applying (A.1.7'), we have\\[0.1in]
 (A.1.16)
 \begin{eqnarray*}
& \ & \!\!\!\!\!\!\!\!\!\!\int_{B_{\xi_1} ({\rho_\mu})}{\cal I}  \cdot [\,\partial_{\lambda_1} \,V_1\,]
\  \left( =  \int_{B_{\xi_1} (\rho_\mu)} \!\left\{\, \left( V_1^{{n\,+\,2}\over{n\,-\,2}}  +\ V_2^{{n\,+\,2}\over{n\,-\,2}} \right) \ \,-(\,V_1 \ + \ V_2)^{{n\,+\,2}\over{n\,-\,2}}  \,\right\} \cdot [\,\partial_{\lambda_1} \,V_1\,] \right)\\[0.15in]
& = & \int_{B_{\xi_1} (\rho_\mu)} \!\left\{\,-\,  {{n\,+\,2}\over {n\,-\,2}} \cdot V_1^{4\over {n\,-\,2}} \cdot V_2 \, + \, O \,(1) \cdot V_2^{{n \,+\,2}\over {n\,-\,2}}
\,\right\} \cdot [\,\partial_{\lambda_1} \,V_1\,] \ \ \ \   (\,V_1\,   = \ V_{\lambda_1\,,\ \xi_1}\,, \ \   V_2 \, = \ V_{\lambda_2\,,\ \xi_2})\\[0.15in]
& = &   +\,{{n + 2}\over 2} \cdot  \int_{B_{\xi_1}\,({\rho_\mu})} V_1^{4\over {n\,-\,2}}
\cdot \left( {{ \lambda_2^{{n\,-\,2}\over 2 } }\over {\gamma^{n\,-\,2}   }} \,\cdot  \, [\,1 \ + \ o\,(1)\,]  \right) \cdot
\left\{   \lambda_1^{{n \,-\, 4}\over 2}
\cdot {{ (\,\lambda_1^2\, -\, |\, y \,-\, \xi_1|^{\,2})}\over {(\,\lambda_1^2 \,+\, |\, y \,-\, \xi_1|^{\,2})^{{n }\over 2} }}\right\}
\\[0.15in]
& \ & \ \ \ \  \ \  \ \ \ \ \  \ \  \  \ \ \ \  \ \  \  \ \ \ \  \ \  \ \ \ \ \  \ \  \  \ \ \ \  \ \  \  \ \  \  \ \ \ \  \ \  \  \ \  \  \ \ \ \  \ \  \  +   \  \, O\,(1) \cdot \int_{B_{\xi_1}\,({\rho_\mu})} V_2^{{n \,+\,2}\over {n\,-\,2}} \cdot [\,\partial_{\lambda_1} \,V_1\,]  \\[0.15in]
& = &   \omega_n\cdot {{n + 2}\over 2} \!\cdot  {{\lambda_1^{{n\,-\,2 }\over 2}\!\cdot\lambda_2^{{n\,-\ 2}\over 2}}\over {\gamma^{n\,-\,2} }}  \cdot \!{1\over {\lambda_1}}\! \times\!\!\!
\int_0^{\rho_\mu}\!
\left[    \left( {{\lambda_1}\over {\lambda^2_1 \,+\, r^2}} \right)^{\!\!2}\!\!
\cdot {{ (\,\lambda_1^2 \ -\  r^{\,2})}\over {(\,\lambda_1^2 \ + \  r^{\,2})^{{n }\over 2} }}\right] \cdot r^{n\,-\,1} \cdot dr    \cdot    [\,1 \, + \, o\,(1)\,]   \\[0.15in]
& \ & \ \ \ \  \ \ \ \ + \ \ {\bf II} \ \ \ \  \ \ \ \ \ \ \ \ \  \ \ \ \  \ \ \ \  \ \ \ \  [\,\uparrow \ \ {\mbox{to \ \ determine \ \ the \ \ sign}}\,]\\[0.15in]
    & \ & \!\!\!\!\!\!\!\!\!\!  \bigg\{ \ {\bf II} \ = \  O\,(1) \cdot \int_{B_{\xi_1}\,({\rho_\mu})} V_2^{{n \,+\,2}\over {n\,-\,2}} \cdot [\,\partial_{\lambda_1} \,V_1\,] \ \ \ \ {\mbox{is \ \ estimated \ \ in}} \ \ (A.1.21)\,;   \ \   \omega_n \ = \ {\mbox{Vol}\,} (\,S^{n \,-\,1}\,) \ \bigg\} \ .
    \end{eqnarray*}
    Via the change of variables $\,r \ = \ \lambda_1 \cdot \tan\,\theta \,$,\, we continue with
    \begin{eqnarray*}
\!\!\!\!\!\!\!(A.1.17) \ \ \     \   & \ & \int_0^{\rho_\mu}\!
\left[    \left( {{\lambda_1}\over {\lambda^2_1 \,+\, r^2}} \right)^{\!\!2}\!\!
\cdot {{ (\,\lambda_1^2 \ -\  r^{\,2})}\over {(\,\lambda_1^2 \ + \  r^{\,2})^{{n }\over 2} }}\right] \cdot r^{n\,-\,1} \cdot dr\\[0.15in]
& = &
\int_0^{\arctan\left( {{\rho_\mu}\over \lambda_1}\right)} \left\{\, [\,\sin\,\theta\,]^{\,n\,-\,1}  \cdot [\,\cos\,\theta\,]^{\,3} \ - \
[\,\sin\,\theta\,]^{\,n\,+\,1}  \cdot [\,\cos\,\theta\,] \right\} \,d\theta  \\[0.15in]
& = &
\int_0^{\arctan\left( {{\rho_\mu}\over \lambda_1}\right)} \left\{\, [\,\sin\,\theta\,]^{\,n\,-\,1}  \cdot [\,\cos\,\theta\,]  \ - \
2\,[\,\sin\,\theta\,]^{\,n\,+\,1}  \cdot [\,\cos\,\theta\,] \right\} \,d\theta \\[0.15in]
& = &
 [\,1 \ + \ o\,(1)\,] \cdot  \int_0^{ \,{\pi\over 2}} \left\{\, [\,\sin\,\theta\,]^{\,n\,-\,1}     \ - \
2\,[\,\sin\,\theta\,]^{\,n\,+\,1}  \,\right\} \,d\,[\,\sin \,\theta\,]  \ \ \ \   \ \ \ \left( {{ \rho_\mu}\over {  \lambda_1 }} \ \gg \ 1\right) \\[0.1in]
& = &  [\,1 \ + \ o\,(1)\,]   \cdot
\left[ \ {1\over n} \ - \ {2\over {n\,+\,2}}\,\right] \ = \ [\,1 \ + \ o\,(1)\,]   \cdot  (-1) \cdot {{n\,-\,2}\over {n\,(n\,+\,2)  }}   \ \ \ \ \ \ \ \\[0.15in]
& \ & \ \ \  \ \ \ \  \ \ \ \ \  \ \ \ \  \ \ \ \ \  \ \ \ \  \ \ \ \ \  \ \ \ \  \ \ \ \ \  \ \ \ \  \ \ \ \ \  \ \ \ \  \ \ \ \ \  \ \ \ \  \ \ \ \ \    [\,\uparrow \ \ {\mbox{-\,ve}}\,]\ .
\end{eqnarray*}

\newpage

It follows that
$$
\int_{B_{\xi_1} ({\rho_\mu})} {\cal I} \cdot [\,\partial_{\lambda_1} \,V_1\,]  \ = \   \,-\ \left( \omega_n \,\cdot  \,  {{n - 2}\over {2\,n}} \,\right) \cdot {{\lambda_1^{{n\,-\ 2}\over 2} \cdot \lambda_2^{{n\,-\ 2}\over 2}}\over {\gamma^{n\,-\,2} }}
 \cdot {1\over {\lambda_1}} \cdot [\,1\ + \ o\,(1)\,]  \ + \   {\bf II} \ . \leqno (A.1.18)
$$
Here $\,o\,(1)\ \to \ 0\,$ as $\,{\bf D} \ \to \ \infty\,$   and
$$
{{\rho_\mu}\over {\lambda_1}} \  \to \ \infty\,. \leqno (A.1.19)
$$
In (A.1.18), we recognize the leading term in (2.23) after multiplying by the number $\,n(n\,-\,2)\,.\,$ \bk
As for the remainder, we utilize (A.1.8)\,,\, (A.1.15),  and continue with\\[0.1in]
 (A.1.20)
 \begin{eqnarray*}
|\, {\bf II}\,| & \le & C\cdot {1\over {\lambda_1}} \cdot \int_{B_{\xi_1}({\rho_\mu})} V_2^{{n \,+\,2}\over {n\,-\,2}} \cdot V_1 \ \ \ \ \ \ \ \  \ \ \ \ \ \ \ \ \ \ \ \    \left( \ \bigg\vert \ {{\partial V_1}\over {\partial \lambda_1}} \,\bigg\vert  \ \le \ {{n\,-\,2}\over 2} \cdot {1\over {\lambda_1}} \cdot V_1\,\right) \\[0.1in]
& \le & C\,(n) \cdot {1\over {\lambda_1}} \cdot {1\over {\lambda_1^{{n\,+\,2}\over 2} }} \cdot {1\over { {\bf d}^{n\,+\,2} }}  \cdot [\,1 \ + \ o\,(1)\,]
\cdot \int_0^{\rho_\mu} \left( {\lambda_1\over { \lambda^2_1 \  +  \ r^2 }}\right)^{\!\! {{n\,-\,2}\over 2} } \cdot r^{n\,-\,1}\cdot dr\ .\\[0.15in]
\end{eqnarray*}

\vspace*{-0.3in}

Let $\,r \ = \ \lambda_1 \cdot \tan\,\theta\,.\, $ We come to
\begin{eqnarray*}
   (A.1.21) \ \ \ \ \ \ \ \ \ \ \ \ \ \ \  & \ & \int_0^{\rho_\mu}\!\left( {\lambda_1\over { \lambda^2_1 \  +  \ r^2 }}\right)^{\!\! {{n\,-\,2}\over 2} } \!\!\! \cdot r^{n\,-\,1}\cdot dr \\[0.15in]
 & = &
  \int^{\arctan \left( {{\rho_\mu}\over \lambda_1}\right)}_0 {{ \lambda_1^{{n\,-\,2}\over 2} \cdot \lambda_1 \cdot \lambda_1^{n\,-\,1}  }
\over { \lambda^{n\,-\,2}_1 }} \cdot {{[\,\sin\,\theta\,]^{n\,-\,1} }\over {[\,\cos \,\theta\,]^{n\,-\,1\,+\,2}  }}
\cdot [\, \cos \,\theta\,]^{n\,-\,2} \ d\,\theta\\[0.15in]
& \le  & C \cdot  \lambda_1^{{n\,+\,2}\over 2}   \cdot \bigg\vert\,\int^{\arctan \left( {{\rho_\mu}\over \lambda_1}\right)}_0 {{d\,[\,\cos \theta\,]}\over {\ [\,\cos \,\theta\,]^3\    }}\,\bigg\vert\\[0.15in]
& \le & \lambda_1^{{n\,+\,2}\over 2}   \cdot \bigg\vert \ {{1}\over {[\,\cos \,\theta\,]^2   }}\ \bigg\vert^{\,\arctan \left( {{\rho_\mu}\over \lambda_1}\right)}_{\,0} \\[0.15in]
& \le & C \cdot \lambda_1^{{n\,+\,2}\over 2}   \cdot \bigg\vert \   {1\over { \cos^2 \left( \arctan \left( {{\rho_\mu}\over \lambda_1} \right)  \right) }} \ - \ 1 \ \bigg\vert  \ \ \ \ \ \left[ \ \arctan \left( {{\rho_\mu}\over \lambda_1} \right)  \ \approx \ {\pi\over 2} \,\right] \\[0.15in]
& \le & C \cdot  \lambda_1^{{n\,+\,2}\over 2}   \cdot \left[\, \tan \left( \arctan \left( {{\rho_\mu}\over \lambda_1} \right)  \right)\,\right]^2  \ \ \ \ \ \left[ \ \sin\left( \arctan \left( {{\rho_\mu}\over \lambda_1} \right)  \right) \ \approx \ 1 \,\right] \\[0.15in]
& \le & C \cdot  \lambda_1^{{n\,-\,2}\over 2}   \cdot \rho_\mu^2\\[0.15in]
\Longrightarrow \ \ |\, {\bf II}\,| & = & O\,(\mu^2) \cdot {1\over {\lambda_1}} \cdot {1\over { {\bf D}^n}}    \ \ \ \ \  \ \ \ \ \  \ \left[\,{\mbox{recall \ \ that \ \ }} \rho_\mu \ = \ \mu \cdot |\,\xi_1\ - \ \xi_2\,|\,\right]\,.\ \ \ \ \ \ \ \ \ \ \ \\ \ \ \ \ \ \ \ \ \ \ \
\end{eqnarray*}
 Here we apply (2.14)\,.

 \newpage

\hspace*{0.5in}In a similar way, we consider the  integral on the other ball.\\[0.1in]
 (A.1.22)
 \begin{eqnarray*}
& \ & \bigg\vert\ \int_{B_{\xi_2}({\rho_\mu})} {\cal I} \cdot [\,\partial_{\lambda_1} \,V_1\,]\ \bigg\vert \ \le \   C\,(n) \cdot {1\over {\lambda_1}} \cdot  \int_{B_{\xi_2}({\rho_\mu})} V_2^{4\over {n\,-\,2}} \cdot V_1 \cdot V_1\ \
 + \  {\bf{II}}{\,'}\\[0.15in]
& \le & C\,(n) \cdot {1\over {\lambda_1}} \cdot
\left( {1\over {\lambda_2^{{n\,-\,2}\over 2} }} \cdot {1\over {{\bf d}^{n\,-\,2} }} \right)^{\!\!2}  \cdot
\int^{\rho_\mu}_0 \left( {{\lambda_2}\over {\lambda^2_2 \ + \ r^2}} \right)^{\!\!2}\,r^{n\,-\,1} \cdot dr\ \
 + \  {\bf{II}}{\,'}\\[0.15in]
 & \ & \ \ \ \ \ \ \ \ \ \ \ \ \ \ \  \ \ \ \ \ [ \ \uparrow \ \ {\mbox{similar \ \ to  \ \ (A.1.12) \ \ via \ \ symmetry}}\,] \\[0.15in]
& \le &   C (n) \cdot \!{1\over {\lambda_1}} \cdot \! \!
\left( {1\over {\lambda_2^{{n\,-\,2}\over 2} }}   \cdot  \!{1\over {{\bf d}^{n\,-\,2} }} \right)^{\!\!2}   \!\cdot \! \! \int^{\arctan \left( {{\rho_\mu}\over \lambda_2}\right)}_0 \ {{ \lambda^2_2 \cdot \! \lambda^{n\,-\,1}_2 \cdot \! \lambda_2   }
\over { \lambda^4_2 }} \cdot {1\over {[\,\cos\,\theta\,]^{n\,-\,1}}} \cdot \!
{{[\,\cos \,\theta\,]^{\,4}}
\over { [\,\cos \,\theta\,]^{\,2} }}  \ d\,\theta\\[0.15in]
& \ & \ \ \ \ \ \ \
 + \  \ {\bf{II}}{\,'}\\[0.15in]
& \le & -\,  C\,(n) \cdot {1\over {\lambda_1}} \cdot
\left(   {1\over {{\bf d}^{n\,-\,2} }} \right)^{\!\!2}  \cdot \int^{\arctan \left( {{\rho_\mu}\over \lambda_2}\right)}_0 \
\left[\,{{1}\over {[\,\cos\,\theta\,]^{n\,-\,3}}} \,\right]   \ d\,[\,\cos\,\theta\,]\ \
 + \  {\bf{II}}{\,'}\\[0.15in]
& \le &   C\,(n) \cdot {1\over {\lambda_1}} \cdot
\left(   {1\over {{\bf d}^{n\,-\,2} }} \right)^{\!\!2}  \cdot \left[\,{{1}\over {[\,\cos\,\theta\,]^{n\,-\,4}}} \,
\right]^{\arctan \left( {{\rho_\mu}\over \lambda_2}\right)}_0      \ \
 + \  {\bf{II}}{\,'}\\[0.15in]
 & \le &   C\,(n) \cdot {1\over {\lambda_1}} \cdot
\left(   {1\over {{\bf d}^{n\,-\,2} }} \right)^{\!\!2}  \cdot \left[\,{{1}\over {[\,\cos\,\arctan \left( {{\rho_\mu}\over \lambda_2}\right)\,]^{n\,-\,4}}} \ \ - \ 1
\right]       \ \
 + \  {\bf{II}}{\,'}\\[0.15in]
 & \le &   C \,(n) \cdot {1\over {\lambda_1}} \cdot
\left(   {1\over {{\bf d}^{n\,-\,2} }} \right)^{\!\!2}  \cdot \left[\,{{[\,\sin\,\arctan \left( {{\rho_\mu}\over \lambda_2}\right)\,]^{n\,-\,4}}\over {[\,\cos\,\arctan \left( {{\rho_\mu}\over \lambda_2}\right)\,]^{n\,-\,4}}}
\right]       \ \
 + \  {\bf{II}}{\,'} \\[0.15in]
 & \ & \ \ \ \ \  \ \ \ \  \ \ \  \ \ \ \ \  \ \ \ \  \ \ \  \ \ \ \ \  \ \ \ \  \ \ \  \left\{ \ \sin\,\arctan \left( {{\rho_\mu}\over \lambda_2}\right) \ \approx \ 1\,,\ \ \cos\,\arctan \left( {{\rho_\mu}\over \lambda_2}\right) \ \approx \ 0 \  \right\} \\[0.15in]
& \le &  C\,(n) \cdot {1\over {\lambda_1}} \cdot
\left(   {1\over {{\bf d}^{n\,-\,2} }} \right)^{\!\!2}  \cdot [\,\tan\,\arctan \left( {{\rho_\mu}\over \lambda_2}\right)\,]^{n\,-\,4}  + \  \ {\bf{II}}{\,'} \\[0.15in]
&\le &  C\,(n) \cdot {1\over {\lambda_1}} \cdot
\left(   {1\over {{\bf d}^{n\,-\,2} }} \right)^{\!\!2}  \cdot \left( {{\rho_\mu}\over \lambda_2}\right)^{\!\!n\,-\,4}\ \ + \  \ {\bf{II}}{\,'}\\[0.15in]
& \le & O\,(\mu^{n\,-\,4}) \cdot  {1\over {\lambda_1}} \cdot {1\over {{\bf d}^{n } }} \ \
 + \  {\bf{II}}{\,'} \\[0.15in]
& \le & o\,(1) \cdot  {1\over {\lambda_1}} \cdot {1\over {{\bf D}^{n } }} \ \
 + \  {\bf{II}}{\,'} \ \ \  \ \ \   \
  \ \ \  \ \ \    \ \ \  \ \ \    \ \ \  \
  \ \ \  \ \ \    \ \ \  \ \      [\,o\,(1) \ \to \ 0^+ \ \ {\mbox{as}} \  \
 {\bf D} \ \to \ \infty\,] \ . \\
\end{eqnarray*}
In (A.1.22), \\[0.1in]
 (A.1.23)
 \begin{eqnarray*}
 \!\!\!\!\!\!{\bf II}^{\,'} & = & O\,(1) \cdot \int_{B_{\xi_2}\,({\rho_\mu})} V_1^{{n \,+\,2}\over {n\,-\,2}} \cdot [\,\partial_{\lambda_1} \,V_1\,] \ \le \ C \cdot \int_{B_{\xi_2}\,({\rho_\mu})} V_1^{{n \,+\,2}\over {n\,-\,2}} \cdot {{  V_1}\over {\lambda_1}} \ \ \ \ \ \  [\,{\mbox{via}} \ \ (A.1.15)\,] \ \ \ \ \ \ \\[0.15in]
\Longrightarrow \ \  |\,{\bf II}^{\,'}\,| & \le  & C\,(n) \cdot {1\over {\lambda_1}} \cdot \left( {1\over {\lambda_2^{{n\,-\,2}\over 2} }} \cdot {1\over {{\bf d}^{n\,-\,2} }} \right)^{\!\!{{2n}\over {n\,-\,2}}} \cdot {\rho_\mu}^n \ = \ C\,(n) \cdot {1\over {\lambda_1}} \cdot  {1\over {{\bf d}^{2\,n} }} \cdot \left( {{{\rho_\mu}}\over {\lambda_2}}\right)^{\!\!n} \\[0.15in]
& = & O\,(\mu^{n}) \cdot {1\over {\lambda_1}} \cdot  {1\over {{\bf d}^{n} }} \ = \ o\,(1) \cdot  {1\over {\lambda_1}} \cdot  {1\over {{\bf D}^{n} }} \ .
\end{eqnarray*}


\vspace*{0.3in}

{\bf Outside} $\,B_{\xi_1}\,({\rho_\mu}) \ \cup \ B_{\xi_2}\,({\rho_\mu})\,.$\,  \ \ Recall the following inequalities.
$$
\bigg\vert \,(a \ + \ b)^{{n\,+\,2}\over {n\,-\,2}} \ - \ \left( a^{{n\,+\,2}\over {n\,-\,2}} \
+ \ b^{{n\,+\,2}\over {n\,-\,2}}\right) \,\bigg\vert \ \le \ C(n) \cdot  \left[\ a^{{4}\over {n\,-\,2}} \
+ \ b^{{4}\over {n\,-\,2}}\,\right]\cdot \min \ \{ a, \ b\}\,. \leqno (A.1.24)
$$
and
$$
\bigg\vert \,(a \ + \ b)^\alpha \ - \ \left( a^\alpha \
+ \ {{\alpha\,+\,2}\over {\alpha\,-\,2}} \cdot a^{\alpha \,-\,1} \cdot b \,\right) \,\bigg\vert \ \le \ C(n\,,\,\alpha) \cdot
\left[\ b^{\,\alpha} \ + \
 \ a^{\alpha\,-\,2}\cdot \min \ \{ a^2, \ b^2\}\,\right]\,.\leqno (A.1.25)
$$
$$
M\ \ge\ 2\ \ \Longrightarrow\ \ \Big|\,(\,a \,+\,b\,)^{\,M}\ -\ (\,a^M\,+\,b^{\,M}\,)\,\Big|\ \le\ {\bar C}_M \cdot \left(\,a^{M-1} \cdot b\ +\ b^{M-1}\cdot a\,\right)\,. \leqno (A.1.26)
$$
Here $\,a\,$ and $\,b\,$ are positive numbers\,,\, and $\,\alpha\,$ is a fixed positive number.
See (1.1) and (1.2) in \cite{Twin}\,.\, Cf. also \cite{Progress-Book} and \cite{Yan}\,. We continue with  \\[0.1in]
(A.1.27)
\begin{eqnarray*}
& \ & \bigg\vert\ \int_{\R^n \setminus\,[\,B_{\xi_1}\,({\rho_\mu}) \ \cup \ B_{\xi_2}\,({\rho_\mu})\,]}  \left\{\, (V_{\lambda_1\,,\ \xi_1} \ + \ V_{\lambda_2\,,\ \xi_2})^{{n\,+\,2}\over{n\,-\,2}} \ \,- \
\left( V_{\lambda_1\,,\ \xi_1}^{{n\,+\,2}\over{n\,-\,2}}\ +\ V_{\lambda_2\,,\ \xi_2}^{{n\,+\,2}\over{n\,-\,2}} \right)
\,\right\} \cdot [\,\partial_{\lambda_1} \,V_1\,]\ \bigg\vert\\[0.15in]
& \le & \int_{\big\{\,\R^n \setminus\,[\,B_{\xi_1}\,({\rho_\mu}) \ \cup \ B_{\xi_2}\,({\rho_\mu})\,]\,\big\} \ \cap \ \{\,V_1 \ > \ V_2\,\}} \left[\,V_1^{{4}\over {n\,-\,2}} \ + \ V_2^{{4}\over {n\,-\,2}}\,\right] \cdot V_2 \cdot {{V_1}\over {\lambda_1}} \ \ + \ \ \ \ [\,{\mbox{via}} \ \ (A.1.15)\,] \\[0.15in]
&   & \ \ \ \ \ \ \ \ \ \ + \ \int_{\big\{\,\R^n \setminus\,[\,B_{\xi_1}\,({\rho_\mu}) \ \cup \ B_{\xi_2}\,({\rho_\mu})\,]\,\big\} \ \cap \ \{\,V_2 \ > \ V_1\,\}} \left[\,V_1^{{4}\over {n\,-\,2}} \ + \ V_2^{{4}\over {n\,-\,2}}\,\right] \cdot V_1 \cdot {{V_1}\over {\lambda_1}}\\[0.25in]
& \le &  C\,(n) \cdot  {1\over \lambda} \cdot \int_{\rho_\mu}^\infty \left( {\lambda\over {\lambda^2 \ + \ r^2}} \right)^{\!\!n}\,r^{n\,-\,1} \cdot dr\\[0.15in]
& = &    C\,(n) \cdot  {1\over \lambda} \cdot\int_{\arctan \left( {{\rho_\mu}\over \lambda}\right)}^{\pi\over 2} {{ \lambda^n \cdot \lambda \cdot \lambda^{n\,-\,1}  }
\over { \lambda^{2n} }} \cdot [\,\tan \,\theta\,]^{n\,-\,1} \cdot [\, \cos \,\theta\,]^{2\,(n\,-\,1)} \ d\,\theta\\[0.15in]
& \le &    C\,(n) \cdot  {1\over \lambda} \cdot \bigg[\  -\  \int_{\arctan \left( {{\rho_\mu}\over \lambda}\right)}^{\pi\over 2} [\, \cos \,\theta\,]^{n\,-\,1} \ d\,[\,\cos\,\theta\,]\  \bigg]\ \le \
   C\,(n) \cdot  {1\over \lambda} \cdot   \bigg[\ - \,[\, \cos \,\theta\,]^{\,n}\  \bigg]_{\,\arctan \left( {{\rho_\mu}\over \lambda}\right)}^{\,{\pi\over 2} }  \\[0.15in]
& \approx &  C\,(n) \cdot  {1\over \lambda} \cdot    {{\,-1}\over { [\, \tan\,\theta\,]^{\,n} \  \bigg\vert_{\,\arctan \left( {{\rho_\mu}\over \lambda}
\right)}^{\,{\pi\over 2} }  }}   \ \le \
   C\,(n) \cdot  {1\over \lambda} \cdot {1\over {\left( {{\rho_\mu}\over \lambda}
\right)^{\!n}   }} \  = \  C\,(n) \cdot  {1\over \lambda} \cdot {1\over { {\bf d}^{\,n }}} \cdot {1\over {\mu^{n}}}\\[0.15in]
& \le & o\,(1) \cdot {1\over \lambda} \cdot {1\over { {\bf D}^{n\,-\,o\,(1) } }}\  .\\
\end{eqnarray*}
One is led to
$$
{\bf D}^{o\,(1) } \cdot \mu^n \ \to \ \infty \ \  \ \Longrightarrow \ \ \  {1\over { {\bf D}^{o\,(1) } }} \cdot {1\over {\mu^{n}}} \ \to \ 0\ .   \leqno (A.1.28)
$$

\vspace*{0.3in}

Hence we conclude that
\begin{eqnarray*}
(A.1.29)\ \  \ \  & \ &  \int_{\R^n} \left\{\, \left( V_{\lambda_1\,,\ \xi_1}^{{n\,+\,2}\over{n\,-\,2}}\ +\ V_{\lambda_2\,,\ \xi_2}^{{n\,+\,2}\over{n\,-\,2}} \right) \ \,- \ (V_{\lambda_1\,,\ \xi_1} \ + \ V_{\lambda_2\,,\ \xi_2})^{{n\,+\,2}\over{n\,-\,2}}  \,\right\} \cdot [\,\partial_{\lambda_1} V_1\,]\\[0.1in]
& = &  \,-\ \left( \omega_n \,\cdot  \,  {{n - 2}\over {2\,n}} \,\right) \cdot {{\lambda_1^{{n\,-\ 2}\over 2} \cdot \lambda_2^{{n\,-\ 2}\over 2}}\over {\gamma^{n\,-\,2} }}
 \cdot {1\over {\lambda_1}} \cdot [\,1\ + \ o\,(1)\,]  \ + \ {1\over {\lambda_1}} \cdot  O\left(   {1\over {{\bf D}^{\,n\,-\ o\,(1)} }} \right)\ . \ \  \ \  \ \  \ \  \ \  \ \
\end{eqnarray*}
Here $\,o\,(1)\ \to \ 0\,$ as $\,{\bf D} \ \to \ \infty\,$,\, using the property of $\,\mu\,$ described in (A.1.3)\,.\,
Note that
$$\ \ \ \ \ \ \ \
{{\lambda_1^{{n\,-\ 2}\over 2} \!\cdot \lambda_2^{{n\,-\ 2}\over 2}}\over {\gamma^{n\,-\,2} }} \ = \ {1\over { {\bf D}^{n\,-\,2} }} \ \ \ \ \ \ \ \ \ \ \ \ \ \ \ \ \ \ \ \ \ \ \ \  [\,{\mbox{cf}}. \ (2.13)\,]\,.
$$

\newpage

{\bf Extracting the leading term in $\xi$\,-\,derivatives\,.}
\begin{eqnarray*}
 (A.1.30)\ \  \ \ & \ &  {{\partial V_{\lambda_1\,, \, \xi_1}}\over {\partial\,\xi_{1_1}}}\ = \ {{\partial}\over {\partial\,\xi_{1_1}}} \left[\left({\lambda\over {\lambda_1^2 + |\, y \,-\, \xi_{1_1}\,|^{\,2}}} \right)^{\!\!{{n - 2}\over 2}} \,\right]  \\[0.15in]
  & \ & \ \ \ \ \ \ \ \ \ \ \ \ =\,-\,{{n - 2}\over 2} \cdot \lambda_1^{{n - 2}\over 2} \cdot {{ 2\,(\xi_{1_1} \,-\, y_1)}\over {(\,\lambda^2_1 + |\, y\, -\, \xi_{1_1}\,|^{\,2})^{{n }\over 2} }} \\[0.15in]
 \Longrightarrow & \ & \bigg\vert   \  {{\partial V_{\lambda_1\,, \, \xi_1}}\over {\partial\,\xi_{1_1}}}  \,\bigg\vert \ = \ {{n - 2}\over 2} \cdot \lambda_1^{{n - 2}\over 2} \cdot {1\over {\lambda_1}} \cdot {{ 2\,\lambda_1 \cdot |\,\xi_{1_1} - y_1\,|}\over {(\,\lambda^2_1 + |\, y - \xi_1\,|^{\,2})^{{n }\over 2} }} \ \le \  {{n - 2}\over 2} \cdot {1\over {\lambda_1}} \cdot  V_{\lambda_1\,, \  \xi_1} \,. \ \ \ \ \ \ \ \ \ \
\end{eqnarray*}
We apply the same partition as in (A.1.5)\,,\, and also (A.1.13) for the definition of $\,{\cal I}$\,.\, Following (A.1.16) and using (A.1.30)\,,\, we have \\[0.1in]
(A.1.31)
$$
\int_{B_{\xi_1} (\rho_\mu)}   {\bf I} \cdot [\,\partial_{\xi_{1_1}} \,V_1\,]\ = \ \,+\,{{n + 2}\over 2}   \cdot \int_{B_{\xi_1} (\rho_\mu)} V_1^{4\over {n\,-\,2}}
\cdot V_2 \cdot
\left\{   \  \lambda_1^{{n - 2}\over 2} \cdot {{ 2\,(\xi_{1_1} - y_1)}\over {(\,\lambda^2_1 + |\, y \ - \ \xi_1|^{\,2})^{{n }\over 2} }}\right\} \ \ + \ \ {\bf III} \ .
$$
{\bf III}\, is estimated as in (A.1.21)\,.\, Recall (A.1.14) and (A.1.15)\,.\, In particular,
$$
 {\bf I} \ = \ -\,{\cal I} \ .
$$
Consider the main term in (A.1.31)\,:
\begin{eqnarray*}
(A.1.32) \!\!\!\!\!\!\!& \  &  \,+\,{{n + 2}\over 2}   \cdot \int_{B_{\xi_1} (\rho_\mu)} V_1^{4\over {n\,-\,2}}
\cdot V_2 \cdot
\left\{   \  \lambda_1^{{n - 2}\over 2} \cdot {{ 2\,(\xi_{1_1} - y_1)}\over {(\,\lambda^2_1 + |\, y \ - \ \xi_1|^{\,2})^{{n }\over 2} }}\right\}\\[0.15in]
& = &  \,+\,{{n + 2}\over 2} \cdot    \int_{B_{\xi_1} (\rho_\mu)} V_1^{4\over {n\,-\,2}}
\cdot \left( {1\over {\lambda_1^{{n\,-\,2}\over 2} }} \cdot {1\over {{\bf d}^{n\,-\,2} }} \right) \cdot \left\{ \, 1\ - \ {{n\,-\,2}\over 2} \cdot {\bf t} \cdot [\,1 \ + \ o\,(1)\,] \,\right\} \times \\[0.15in]
& \ & \ \ \ \ \  \ \ \ \ \ \  \ \ \ \ \ \  \ \ \ \ \ \  \ \ \ \ \ \  \ \ \ \ \ \  \ \ \ \ \ \  \ \ \ \ \ \  \  \ \ \ \ \  \ \ \ \ \ \  \  \times \,
\left\{   \  \lambda_1^{{n - 2}\over 2} \cdot {{ 2\,(\xi_{1_1} - y_1)}\over {(\,\lambda^2_1 + |\, y \ - \ \xi_1|^{\,2})^{{n }\over 2} }}\right\} \\[0.15in]
& \ & \hspace*{4.7in} [\,{\mbox{via \ \ (A.1.8)\,]}} \\[0.15in]
& \ & \hspace{3.9in}(\,\downarrow \ \ {\mbox{take \   note  \   of \  the\ sign}}\,)\\[0.15in]
&=  & -\,{{n + 2}\over 2} \cdot {{n -  2}\over 2} \cdot \left(   {1\over {{\bf d}^{n\,-\,2} }} \right) \cdot   \int_{B_{\xi_1} (\rho_\mu)} V_1^{4\over {n\,-\,2}} \cdot {\bf t}
\cdot
\left\{   \,  {{ 2\,(\xi_{1_1} - y_1)}\over {(\,\lambda^2_1 + |\, y \ - \ \xi_1|^{\,2})^{{n }\over 2} }}\right\} \cdot [\,1 \ + \ o\,(1)\,] \\[0.15in]
& \ & \hspace*{4.2in} [\,{\mbox{see \ \ (A.1.33) \ \ below\,]}} \ .
\end{eqnarray*}
Here we make use of the ``\,symmetry\,"\\[0.1in]
(A.1.33)
\begin{eqnarray*}
\int_{B_{\xi_1} (\rho_\mu)} V_1^{4\over {n\,-\,2}} \cdot {{(y_1 \,-\,\xi_{1_1} )}\over {(\,\lambda^2_1 + |\, y \ - \ \xi_1|^{\,2})^{{n }\over 2} }} \ dy \ = \ \int_{B_{o} (\rho_\mu)} \left( {{\lambda_1}\over {\lambda_1^2 \ + \ |\,{\bar y}\,|^2 }}\right)^{\!\!2}  \cdot  {{   1 }\over {(\,\lambda^2_1 \,+\, |\,{\bar y}\,|^2 )^{{n }\over 2} }} \cdot {\bar y}_1  \ d{\bar y}  \ = \ 0\ ,
\end{eqnarray*}
where $\, {\bar y} \ = \ y \ - \ \xi_1\,.$\hspace*{2.5in}\ \ \ \ \ \ [\,$+\,/\,-\,$ cancelation \ \ $\uparrow$\,]\, \\[0.1in]
Recall that
$$
{\bf t} \ = \ {1\over {{\bf d}^{2} }} \cdot \left( {{\lambda_2}\over {\lambda_1}}  \right) \ + \ {1\over {{\bf d}^{2} }} \cdot {{|\,\,y\ - \ \xi_1\,|^2}\over
{\lambda_1 \cdot \,\lambda_2}}   \ + \ {1\over {{\bf d}^{2} }} \cdot {{ \,2\,(\,y\ - \ \xi_1)\cdot(\,\xi_1\ - \ \xi_2)}\over
{\lambda_1 \cdot \,\lambda_2}} \ . \leqno (A.1.34)
$$
\hspace*{1.6in}$\leftarrow$\hspace*{1.8in}$\rightarrow$\hspace*{1.5in}$\uparrow$\\[0.1in]
Likewise, the integrations involving the first two terms in (A.1.34) also yield  nothing.  Finally, we come down to  the integral\\[0.1in]
\hspace*{2.2in}\ \ \ \ \ \ [\,change \ \ of \ \ sign \ \ $\downarrow$\,]

\vspace*{-0.35in}

$$ + \cdot{{n + 2}\over 2} \cdot {{n -  2}\over 2} \cdot \left(   {1\over {{\bf d}^{n\,-\,2} }} \right) \cdot \,4\, \cdot {1\over {{\bf d}^{2} }} \cdot {{ 1}\over
{\lambda_1 \cdot \,\lambda_2}}   \int_{B_{\xi_1} (\rho_\mu)} \!\!\!V_1^{4\over {n\,-\,2}} \cdot  {{  ( y_1 \ - \ \xi_{1_1} ) \cdot \left[\,{\displaystyle{\sum_{j\,=\,1}^n}} \,(\,y\ - \ \xi_1)_{|_j} \cdot(\,\xi_1\, - \, \xi_2)_{|_j}\,\right] }\over {(\,\lambda^2_1 + |\, y \ - \ \xi_1|^{\,2})^{{n }\over 2} }}\ .
$$
Via symmetry again,
$$
 \int_{B_{\xi_1} (\rho_\mu)} \!\!\!V_1^{4\over {n\,-\,2}} \cdot  {{   1 }\over {(\,\lambda^2_1 \,+\, |\, y \ - \ \xi_1|^2 )^{{n }\over 2} }} \cdot  ( y_1 \ - \ \xi_{1_1} ) \cdot (\,y\ - \ \xi_1)_{\,|_j} \ = \ 0 \ \ \ \ \ \ \mfor \ \ j \ \not= 1\,.
$$
Thus we are left with\\[0.1in]
(A.1.35)
 \begin{eqnarray*}
 & \ &     {{n\,+\,2}\over {{\bf d}^{n} }} \cdot {{ n\,-\,2}\over
{\lambda_1 \cdot \,\lambda_2}} \cdot  \left[ \ \int_{B_{\xi_1} (\rho_\mu)} \!\!V_1^{4\over {n\,-\,2}} \cdot  {{   1 }\over {(\,\lambda^2_1 \,+\, |\, y \ - \ \xi_1|^2 )^{{n }\over 2} }} \cdot  ( y_1 \ - \ \xi_{1_1} )^2\ dy \ \right] \cdot (\,\xi_{1_1}\ - \ \xi_{2_1}\,)\\[0.15in]
& = &    C (n) \cdot  [\,1 \ + \ o\,(1)\,] \cdot {1\over {{\bf d}^{n} }} \cdot {1\over
{\lambda_1 \cdot \,\lambda_2}} \,\cdot (\,\xi_{1_1}\ - \ \xi_{2_1}\,) \ = \ O \left( {1\over {\lambda_1}} \cdot {1\over { {\bf D}^{\,n\,-\,1} }} \right)\,,\\[0.15in]
& \ & \!\!\!\!\!{\mbox{where}} \ \ \ \  \ \ \ \ \ \  C (n) \ = \ (n \,+\,2) \,(n\,-\,2)\cdot  \int_{\R^n} {{Y_{|_1}^2 }\over {  (1 \ + \ |\,Y\,|^2)^{\,{n\over 2} \,+\,2}   }} \ dY  \ >  \ 0\,,\\[0.15in]
& \ &  {\mbox{via \ \ the \ \ changes \ \ of \ \ variables \ \ }} \ \ {\bar y} \ = \ y \ - \ \xi_1 \ \ \ {\mbox{and}} \ \ \  Y \ = \ {{\bar y}\over {\lambda_1}} \ .
 \end{eqnarray*}
 Here $\,o\,(1) \ \to \ 0\,$ as $\,{\bf D} \ \to \ \infty\,$ and $\,\mu \cdot {\bf D} \ \to \ \infty\,.\,$
\newpage

We estimate the integral on the other ball via\\[0.1in]
(A.1.36)
 \begin{eqnarray*}
& \ & \!\!\!\!\! \bigg\vert \,\int_{B_{\xi_2} (\rho_\mu)} V_2^{4\over {n\,-\,2}} \cdot V_1 \cdot [\,\partial_{\,\xi_{1_1}} \,V_1\,]\,\bigg\vert\\[0.15in]
& \ & \!\!\!\!\!\!\!\!\!\!\!\!\!  = \   {{n - 2}\over 2}   \! \cdot \!\bigg\vert \, \int_{B_{\xi_2}\,(\rho_\mu)} V_2^{4\over {n\,-\,2}}
\cdot \left(  {{ \lambda_1}\over { \lambda_1^2 \, + \, |\,y\, - \, \xi_1\,|^2  }}\right)^{\!\!{{n\,-\,2}\over 2}  } \cdot
\left\{   \  {1\over {\lambda_1}} \cdot {{ 2\,\lambda_1^{n\over 2} }\over {(\,\lambda^2_1 + |\, y \ - \ \xi_1|^{\,2})^{{n }\over 2} }}\right\}
\cdot (\xi_{1_1} - y_1) \,\bigg\vert   \\[0.15in]
& \ & \!\!\!\!\!\!\!\!\!\!\!\!\!  = \  (n\,-\,2) \cdot\!  {1\over {\lambda_1}}  \!\cdot \!\bigg\vert \! \int_{B_{\xi_2} (\rho_\mu)} \!\!V_2^{4\over {n\,-\,2}}
 \!\cdot \!\!\left(  {{ \lambda_1}\over { \lambda_1^2 \, + \, |\,y\, - \, \xi_1\,|^2  }}\right)^{\!\!{{n\,-\,2}\over 2}  }\!\!\!\! \cdot
 \left[\ \left(  {{ \lambda_1}\over { \lambda_1^2 \ + \ |\,y\, - \, \xi_1\,|^2  }}\right)^{\!\!{{n\,-\,2}\over 2}  }
 \, \right]^{\,{n\over {n\,-\,2}} } \!\!\!\!\!\!\cdot(\,\xi_{1_1} - y_1) \,\bigg\vert \\[0.15in]
\\[0.15in]
& \ & \!\!\!\!\!\!\!\!\!\!\!\!\!  = \  C^+(n) \cdot  {1\over {\lambda_1}}  \cdot  \left( {1\over {\lambda_2^{{n\,-\,2}\over 2} }} \cdot {1\over {{\bf d}^{n\,-\,2} }} \right)
\cdot \left( {1\over {\lambda_2^{{n}\over 2} }}
\cdot {1\over {{\bf d}^{n} }} \right) \cdot [\,1\ + \ o\,(1)\,] \,\times\, \\[0.15in]
& \ & \ \ \  \   \times \,\bigg\vert\, \int_{B_{\xi_2} (\rho_\mu)}
\! \left(  {{ \lambda_2}\over { \lambda_2^2 \ + \ |\,y\ - \ \xi_2\,|^2  }}\right)^{\!\!2  }  \cdot
\,[\ (\xi_{1_1} - \,\xi_{2_1} ) \ +\ (\xi_{2_1} - \,y_1 )\ ]\,\bigg\vert  \ \ + \ {\bf{III}}'\\[0.15in]
& \ & \hspace*{4.5in} \ \ \ \ \ \ [\,{\mbox{as \ \ in \ \ (A.1.22)\,}}]\\[0.15in]
& \ & \!\!\!\!\!\!\!\!\!\!\!\!\! \le \    C \!\cdot  \!{1\over {\lambda_1}} \! \cdot\!  \left( {1\over {\lambda_2^{{n\,-\,2}\over 2} }} \cdot {1\over {{\bf d}^{n\,-\,2} }} \right)
\cdot \left( {1\over {\lambda_2^{{n}\over 2} }}
\cdot {1\over {{\bf d}^{n} }} \right) \cdot [\,1\, + \, o\,(1)\,] \cdot\!\! \left[   \int_0^{\rho_\mu} \left(  {{ \lambda_2}\over { \lambda_2^2 \, + \, r^2  }}\right)^{\!\!2  } \!\cdot r^{n\,-\,1} \, dr\,\right] \!\cdot\! |\,\xi_{1_1} - \,\xi_{2_1} \,| \\[0.15in]
& \ & \ \ \ \ \ \ \  \ \ \ \  \ + \ \,{\bf{III}}'\\[0.15in]
& \ & \!\!\!\!\!\!\!\!\!\!\!\!\! \le \   C \cdot  {1\over {\lambda_1}}  \cdot  \left( {1\over {\lambda_2^{{n\,-\,2}\over 2} }} \cdot {1\over {{\bf d}^{n\,-\,2} }} \right)
\cdot \left( {1\over {\lambda_2^{{n}\over 2} }}
\cdot {1\over {{\bf d}^{n} }} \right)  \cdot [\,1\ + \ o\,(1)\,] \,\times\, \lambda_2^{n\,-\,2}
\cdot {\bf d}^{n\,-\,4} \cdot |\,\xi_{1_1} - \,\xi_{2_1} \,|   \\[0.15in]
& \ & \ \ \ \ \ \ \  \ \ \ \  \ + \ \,{\bf{III}}' \hspace*{4in} [\,{\mbox{cf.}} \ \ (A.1.22)\,]\\[0.15in]
& \ & \!\!\!\!\!\!\!\!\!\!\!\!\! = \    C \cdot  {1\over {\lambda_1 \cdot \lambda_2}}  \cdot  {1\over {{\bf d}^{n\,+\,2} }}
 \cdot [\,1\ + \ o\,(1)\,] \,\times\, |\,\xi_{1_1} - \,\xi_{2_1} \,|   \ \ + \ {\bf{III}}' \\[0.15in]
& \ & \!\!\!\!\!\!\!\!\!\!\!\!\! \le \  C \cdot {1\over \lambda}  \cdot {1\over { {\bf D}^{n\,+\,1}}}  \ \ + \ {\bf{III}}'\ .\\
\end{eqnarray*}
The estimates of \,{\bf{III}}'\, and the remaining parts are similar to the case with derivative on $\,\lambda_1$\,,\, using (A.1.30)\,.

\newpage

{\bf \large  \S\,A.2. \ \  Proof of (4.4) \ \& \ (4.5)\,.}\\[0.2in] We  modify the argument presented in the proof of Proposition A.5.27 in  the {\bf e}\,-\,Appendix of  \cite{III}\,.\,
We begin with  (3.31)\,.\, Refer to the notations in (3.37)\,.\,
\begin{eqnarray*}
  I_{\cal R} ({{\bf z}_\sigma})  & = &    I \,({{\bf z}_\sigma} + w_{{\bf z}_\sigma})\\[0.15in]
&\ = \ & \int_{\R^n} \left\{\, {1\over 2}\, \big\langle \btd \,({{\bf z}_\sigma} + w_{{\bf z}_\sigma})\,,\, \btd\,({{\bf z}_\sigma} + w_{{\bf z}_\sigma}) \,\big\rangle   \ - \  {1\over 2} \cdot (n\,-\,2)^2 \!\cdot\![\,{\bf z} \,+\, w_{{\bf z}_\sigma}\,]_+^{{2n}\over {n - 2}} \ \right\}\ \ \\[0.15in]
 & \ & \ \ \ \ \ \ \ \ \ \ \ \  \ \ \ \ \ \ \ \ \ \ \ \ \  \ \  \ \ \ \ \ \ \ \ \ \ \ \ \ \ \  \   + \      \left[\,-\,{{n\,-\,2}\over {2\,n}} \,\right] \cdot \int_{\R^n} (\,{\tilde c}_n \cdot H) \,[\,{{\bf z}_\sigma} + w_{{\bf z}_\sigma}\,]_+^{{2n}\over {n - 2}}\end{eqnarray*}

 \begin{eqnarray*}
\Longrightarrow  & \ & D_{k_\ell} \,I_{\cal R} ({{\bf z}_\sigma}) \ = \   \int_{\R^n} \bigg\{\,   \Big\langle\! \btd \,[\,{{\bf z}_\sigma} \, + \, w_{{\bf z}_\sigma}\,]\,,\ \, \btd\,(D_{k_\ell} \,{{\bf z}_\sigma} \ + \  D_{k_\ell} \,w_{{\bf z}_\sigma}) \,\Big\rangle
 \\
 & \ & \ \  \ \ \  \ \ \ \ \ \ \ \ \  \ \ \  \ \ \ \ \ \ \ \ \ \ \ \ \ \ \ \ \ \ \left. -   \  n\,(n - 2) \,(D_{k_\ell} \,{{\bf z}_\sigma} \, + \, D_{k_\ell} \,w_{{\bf z}_\sigma})\cdot [\,{{\bf z}_\sigma}\, + \,  w_{{\bf z}_\sigma}\,]_+^{{n+2}\over {n - 2}} \,\right\}\\[0.15in]
 & \ & \ \ \ \ \ \ \ \ \  \ \ \ \ \  \    + \      \left[\,-\,{{n\,-\,2}\over {2\,n}} \,\right] \cdot \int_{\R^n} (\,{\tilde c}_n \cdot H)  \cdot (D_{k_\ell} \,{{\bf z}_\sigma} \, + \, D_{k_\ell} \,w_{{\bf z}_\sigma}) \cdot [\,{{\bf z}_\sigma} \, + \, w_{{\bf z}_\sigma}\,]_+^{{n + 2}\over {n - 2}}
\end{eqnarray*}

 \begin{eqnarray*}
\Longrightarrow  & \ & D_{k_\ell} \,I_{\cal R} ({{\bf z}_\sigma})=
 \int_{\R^n} \left\{\,   \langle \,\btd \,{{\bf z}_\sigma}\,,\, \btd\, D_{k_\ell} \,{{\bf z}_\sigma}   \rangle \,- \,n \,(n - 2) (D_{k_\ell} \,{{\bf z}_\sigma})\,{{\bf z}_\sigma}^{{n + 2}\over {n - 2}}\, \right\} \ \ \ \ \ \ \ \ \ \ \ \ \ \\[0.1in]
 & \ &    \hspace*{1in} \leftarrow \ \ \ \ \ \ \ \ \ \ \ \ \ \ \ \ \ \left(\   = \ I'_o\,( {\bf z}_\sigma)\, [\,D_{k_\ell} \,{{\bf z}_\sigma}\,] \,\right)\ \ \ \ \ \ \ \ \ \ \ \ \ \  \rightarrow \\[0.05in]
 & \ & \ \  \ \ \ \ \ +\  \int_{\R^n}   \big\langle \btd \,{{\bf z}_\sigma}\,,\, \btd\, D_{k_\ell} \,w_{{\bf z}_\sigma}   \big\rangle \ +\  \int_{\R^n}   \big\langle \btd \,w_{{\bf z}_\sigma}\,,\, \btd\, D_{k_\ell} \, {{\bf z}_\sigma}    \big \rangle \\[0.15in]
 & \ & \ \ \ \ \ \ \ \ \ \ \ \ \ \ \ \ +\ \int_{\R^n}   \big\langle \btd \,w_{{\bf z}_\sigma}\,,\, \btd\, D_{k_\ell} \,w_{{\bf z}_\sigma}   \big\rangle\\
 & \ & \\
 & \ & \!\!\!\!\!\!\! \!-\, n(n - 2) \,\int_{\R^n} [\,D_{k_\ell} \,{{\bf z}_\sigma} \,] \cdot \left[ ({{\bf z}_\sigma} \,+\, w_{{\bf z}_\sigma})_+^{{n+2}\over {n - 2}} - {{\bf z}_\sigma}^{{n + 2}\over {n - 2}} \right]\\[0.15in]
 & \ &  \ \ \ \ \ - \ n\,(n - 2) \,\int_{\R^n} [\,D_{k_\ell} \,w_{{\bf z}_\sigma}  \, ] \cdot  [\,{{\bf z}_\sigma} + w_{{\bf z}_\sigma}\,]_+^{{n+2}\over {n - 2}}  \\[0.15in]
 & \ & \ \ \ \  \  \ \ \ \ \ \  - \,    \int_{\R^n} ( {\tilde c}_n  \cdot H) \cdot [\, D_{k_\ell} \,{{\bf z}_\sigma} ]\cdot {{\bf z}_\sigma}^{{n \,+\, 2}\over {n \,-\, 2}} \\[0.15in]
 & \ &  \ \ \ \  \ \ \  \        \leftarrow \     (\, {\mbox{key \ \ term}} \ = \ D_{k_\ell} \,[\, G\,({{\bf z}_\sigma})\,] \,)  \rightarrow\\[0.15in]
 & \ &  \ \ \ \ \ \ \ \  \ \ \ \ \ \ - \,      \int_{\R^n} (\,{\tilde c}_n \cdot H)  \cdot[ D_{k_\ell} \,{{\bf z}_\sigma} ]\cdot \left[ \,({{\bf z}_\sigma}+ w_{{\bf z}_\sigma})_+^{{n + 2}\over {n - 2}} - {{\bf z}_\sigma}^{{n + 2}\over {n - 2}}\right]\\[0.15in]
 & \ & \ \ \ \  \ \ \ \ \ \ \ \ \ \ \  \ \ \ \  - \,    \int_{\R^n} (\,{\tilde c}_n \cdot H)  \cdot [\,D_{k_\ell} \,w_{{\bf z}_\sigma}\,]\cdot  [\,{{\bf z}_\sigma} + w_{{\bf z}_\sigma}\,]_+^{{n + 2}\over {n - 2}}
\end{eqnarray*}

\begin{eqnarray*}
  &= & D_{k_\ell} \,[\,  {\bf G}\,({{\bf z}_\sigma})\,]  \ + \  I'_o\,( {\bf z}_\sigma)\, [\,D_{k_\ell} \,{{\bf z}_\sigma}\,]\\[0.2in]
   & \ & \ \ \ \ \ \  +\, \int_{\R^n}   \bigg\{\  \langle \btd \,w_{{\bf z}_\sigma}\,,\, \btd\, D_{k_\ell} \, {{\bf z}_\sigma} \,   \rangle\  -\ n\,(n + 2) \cdot {{\bf z}_\sigma}^{4\over {n - 2}} \,[\,D_{k_\ell}\, {{\bf z}_\sigma}] \cdot w_{{\bf z}_\sigma}\bigg\} \\[0.2in]
 & \ &   \ \ \ \ \  \ \ \ \ \  \left( \leftarrow \ \  \ \ \ \ \  \ \ \ \ \ \ \ \ \ \  = \ (I''_o \,({{\bf z}_\sigma}) \,[\,D_{k_\ell}\, {{\bf z}_\sigma}\,] \,w_{{\bf z}_\sigma} ) \ \ \  \ \ \ \ \  \ \ \ \ \  \ \ \ \  \ \rightarrow \,\right)\\[0.2in]
  & \ &    -\  n
  \,(n - 2) \,\int_{\R^n} [\,D_{k_\ell} \,{{\bf z}_\sigma}\,] \cdot    \left[ \,({{\bf z}_\sigma} + w_{{\bf z}_\sigma})_+^{{n+2}\over {n - 2}}\  - \ {{\bf z}_\sigma}^{{n + 2}\over {n - 2}} \ - \ \left( {{n + 2}\over {n - 2}}\right) {{\bf z}_\sigma}^{4\over {n - 2}}   \,w_{{\bf z}_\sigma}  \,\right]\\[0.2in]
  & \ &    \!\!\!\!\!\!\!\!\!\!\!\!\!  -\, n(n - 2) \,\int_{\R^n} [\,D_{k_\ell} \,w_{{\bf z}_\sigma} \,] \, ({{\bf z}_\sigma} + w_{{\bf z}_\sigma})_+^{{n+2}\over {n - 2}}  \  - \   \int_{\R^n}   \langle \,[\,\Delta \,{{\bf z}_\sigma}\,] \cdot   D_{k_\ell} \,w_{{\bf z}_\sigma}  \\[0.2in]
 & \ &  \ + \ n\,(n + 2) \cdot   \int_{\R^n}     {{\bf z}_\sigma}^{4\over {n - 2}} \,[\,D_{k_\ell}\, w_{{\bf z}_\sigma}] \cdot w_{{\bf z}_\sigma} \ \ \ \ \
   \left\{ \ \uparrow \ \ \Delta {\bf z}_\sigma \ = \ -\,n\,(n - 2) \left[ \,V_1^{{n + 2}\over {n - 2}} \ + \ V_2^{{n + 2}\over {n - 2}} \right] \ \right\}\,\\[0.2in]
   & \ & \ \ \ \ \ \ \ \ +\,  \int_{\R^n}  \bigg\{  \langle \btd \,w_{{\bf z}_\sigma}\,,\, \btd\, D_{k_\ell} \,w_{{\bf z}_\sigma}   \rangle \ -  \ n\,(n + 2) \cdot {{\bf z}_\sigma}^{4\over {n - 2}} \,[\,D_{k_\ell}\, w_{{\bf z}_\sigma}] \cdot w_{{\bf z}_\sigma}\bigg\}\\[0.2in]
 & \ &   \ \ \ \ \   \ \ \left( \, \leftarrow \ \ \ \ \ \ \ \ \ \ \ \ \ \ \ \  \ \ \ \ \ \  = \ ({\bf I}''_o \,({{\bf z}_\sigma}) \,[\,D_{k_\ell}\, w_{{\bf z}_\sigma}\,] \,w_{{\bf z}_\sigma} )  \ \ \ \ \ \ \ \ \ \ \ \ \ \ \ \ \ \ \ \   \rightarrow \,\right)\,\\[0.2in]
   & \ & \ \ \ \ \ \ \ \ \ \ \ \ \ \ \ \  - \     \int_{\R^n}    (\,{\tilde c}_n\cdot  H ) \cdot [ \,D_{k_\ell} \,{{\bf z}_\sigma} ]\cdot \left[ ({{\bf z}_\sigma}+ w_{{\bf z}_\sigma})_+^{{n + 2}\over {n - 2}} \ - \ {{\bf z}_\sigma}^{{n + 2}\over {n - 2}}\right]\\[0.15in]& \ & \ \ \ \ \ \ \ \  \ \ \ \ \ \ \ \ \ \ \ \ \ \ \ - \   \int_{\R^n} ({\tilde c}_n\cdot  H ) \cdot [\, D_{k_\ell} \,w_{{\bf z}_\sigma}\,] \cdot   [\,{{\bf z}_\sigma} + w_{{\bf z}_\sigma}\,]_+^{{n + 2}\over {n - 2}}\\[0.2in]
   &\ = \ &     [\,D_{k_\ell} \, {\bf G}  ({{\bf z}_\sigma})\,] \
   + \    {\bf I}'_o\,( {\bf z}_\sigma)\, [\,D_{k_\ell} \,{{\bf z}_\sigma}\,] \ \ \ \ [\,\leftarrow \ \ {\mbox{refer \ \ to \ \ (2.23 \ \ and \ \ (2.25)}}\,]\\[0.15in]
    & \ &  \ \ \ \ \  \ \ + \  ({\bf I}''_o \,({{\bf z}_\sigma}) \,[\,D_{k_\ell}\, {{\bf z}_\sigma}\,] \,w_{{\bf z}_\sigma}) \ + \ ({\bf I}''_o \,({{\bf z}_\sigma}) \,[\,D_{k_\ell}\, w_{{\bf z}_\sigma}\,] \,w_{{\bf z}_\sigma} ) \ + \  \\[0.2in]
    & \ & \ \ \ \ \   \ \ \ \ \ \  + \ {\bf I} \ + \ {\bf II} \ + \ {\bf III}  \ + \ {\bf IV}\ + \ {\bf V}\,,  \end{eqnarray*}
 where
   \begin{eqnarray*}
    {\bf I} \!& \!\!=\!\!&\! -\,n\,(n - 2) \,\int_{\R^n} [\,D_{k_\ell} \,\,{{\bf z}_\sigma}\,]\!\cdot\!  \left[ ({{\bf z}_\sigma} \ + \ w_{{\bf z}_\sigma})_+^{{n+2}\over {n - 2}} - \ {{\bf z}_\sigma}^{{n + 2}\over {n - 2}} \,-\, \left({{n + 2}\over {n - 2}} \right) {{\bf z}_\sigma}^{4\over {n - 2}}\cdot w_{{\bf z}_\sigma}  \right]\,,\\[0.15in]
     {\bf II}   \!& \!\!=\!\!&\!   -\, n\,(n - 2) \,\int_{\R^n} [\,D_{k_\ell} \,w_{{\bf z}_\sigma}]\!\cdot \!  \left[ \,({{\bf z}_\sigma}\ +\ w_{{\bf z}_\sigma})_+^{{n+2}\over {n - 2}} \,-\, {{\bf z}_\sigma}^{{n + 2}\over {n - 2}} -\left( {{n + 2}\over {n - 2}}\right) {{\bf z}_\sigma}^{4\over {n - 2}} \cdot w_{{\bf z}_\sigma}  \right]\,, \\[0.15in]
     {\bf III}   \!& \!\!=\!\!&\!   - \  \int_{\R^n} ( {\tilde c}_n   \cdot H\,) \cdot [ \,D_{k_\ell} \,{{\bf z}_\sigma} ]\cdot \left[ ({{\bf z}_\sigma}+ w_{{\bf z}_\sigma})_+^{{n + 2}\over {n - 2}} \,-\, {{\bf z}_\sigma}^{{n + 2}\over {n - 2}}\right]\,,\\[0.15in]
     {\bf IV}   \!& \!\!=\!\!&\!   - \    \int_{\R^n} ( {\tilde c}_n   \cdot H\,) \cdot [\,  D_{k_\ell} \,w_{{\bf z}_\sigma} ]\cdot  [\,{{\bf z}_\sigma} \,+\, w_{{\bf z}_\sigma}\,]_+^{\,{{n + 2}\over {n - 2}}}\,,\\[0.15in]
     {\bf V}   \!& \!\!=\!\!&\!   - \ n\,(n + 2)   \int_{\R^n}  [\,  D_{k_\ell} \,w_{{\bf z}_\sigma} ]\cdot  \left[ \,{{\bf z}_\sigma}^{\,{{n + 2}\over {n - 2}}} \ -  \ \left( V_{\lambda_1\,,\,\,\xi_1}^{\,{{n + 2}\over {n - 2}}} \ + \   V_{\lambda_2\,,\,\,\xi_2}^{\,{{n + 2}\over {n - 2}}}\right)\,\right]\ .
 \end{eqnarray*}
For the term
 $$
[\,D_{k_\ell} \, G  ({{\bf z}_\sigma})\,]  \ = \   -\,{{ n\,-\,2 }\over {2\,n}}\cdot D_{k_\ell} \,  \int_{\R^n} (\,{\tilde c}_n\!\cdot H\,) \,{{\bf z}_\sigma}^{{2n}\over {n\,-\,2}}\ ,
 $$
see (4.6), (4.12) and (4.14)\,.\,  As for the term,
 $$
  {\bf I}'_o\,( {\bf z}_\sigma)\, [\,D_{k_\ell} \,{{\bf z}_\sigma}\,] \ = \ n\,(n\,-\,2) \cdot \int_{\R^n}
 \left[ \ \left( V_{\lambda_1\,,\ \xi_1}^{{n\,+\,2}\over{n\,-\,2}}\ +\ V_{\lambda_2\,,\ \xi_2}^{{n\,+\,2}\over{n\,-\,2}} \right) \ - \  (V_{\lambda_1\,,\ \xi_1} \ + \ V_{\lambda_2\,,\ \xi_2})^{{n\,+\,2}\over{n\,-\,2}} \ \right] \cdot [\,D_{k_\ell} \,{{\bf z}_\sigma}\,] \ ,
 $$
refer to (2.23) and (2.25)\,.\, \bk
In the following we take it that [\,see (3.45)\,,\, and compare with (2.14)\,]
 $$
 \lambda^\ell  \ \le \ {{C_o}\over { {\bf D}^{\,n\,-\ \!2}} } \ \ \ \ {\mbox{and}} \ \ \ \ 2 \ \le \ \ell \ < \ n\,-\,2\,.
 $$
Refer to (3.45) in Lemma 3.44.
From Lemma A.5.5 in the  {\bf e}\,-\,Appendix  \cite{III}\,,\, and Lemma 3.44, we have
 \begin{eqnarray*}
(A.2.1) \ \ \ \ \ \ \ \ \ \ \ \  |\,({\bf I}_o'' (\,{\bf z}_\sigma) \,[\, D_{k_\ell} \, {{\bf z}_\sigma}\,] \ w_{ {\bf z}_\sigma} )\,| & = &  O\left( {1\over {{\bf d}^{ 2\,-\, o\,(1)}}}\right) \cdot \max \, \left\{ {1\over {\lambda_1}}\,, \  {1\over {\lambda_2}} \right\} \cdot \,\Vert\, w_{ {\bf z}_\sigma}\Vert_\btd\ \ \ \ \ \ \ \ \ \ \ \ \ \ \ \ \ \\[0.15in]
& \le & \ O\left( {1\over {{\bf d}^{{{n\,+\,2}\over 2} \ + \ 2 }}}\right) \cdot
\max \, \left\{ {1\over {\lambda_1}}\,, \  {1\over {\lambda_2}} \right\}\ .\ \ \ \ \ \ \ \ \ \ \ \ \ \ \ \
 \end{eqnarray*}
Using the uniform bound on the operator $\,{\bf I}_o''\,({\bf z}_{\sigma})\,$, one has
$$|\,({\bf I}''_o \,({\bf z}_{\cal s})\, [\,w_{{\bf z}_{\cal s}}]\,
( \lambda_k \cdot D_{k_\ell}\, w_{{\bf z}_{\cal s}})\,| \ \le \ C \cdot \Vert\, w_{ {\bf z}_{\cal s}}\Vert_\btd \cdot \Vert\,
\lambda_k \cdot \btd\,w_{ {\bf z}_{\cal s}}\Vert_\btd \ \approx\ \ O\left( {1\over {{\bf d}^{{{n\,+\,2}\over 2} \ + \ 2 }}}
\right)\ .  \leqno (A.2.2)
$$

In the following we specify $\,n \ \ge \ 6\,.\,$ Argue as in the proof of Proposition A.5.27\, in the {\bf e}\,-\,Appendix of  \cite{III}\,,\, we obtain
 \begin{eqnarray*}
|\,\lambda \cdot {\bf I}\,| \ = \  -\,n\,(n - 2) \,\int_{\R^n} [\,\lambda \cdot D_{k_\ell} \,\,{{\bf z}_\sigma}\,]\!\!\!\!\!\!\!&\!\cdot\!\!\!\!\!&\!\!\!\!
  \left[ ({{\bf z}_\sigma} \, +\,\ w_{{\bf z}_\sigma})_+^{{n+2}\over {n - 2}} - \ {{\bf z}_\sigma}^{{n\,+\,2}\over {n \,-\, 2}} \,-\, \left({{n\,+\,2}\over {n \,-\, 2}} \right) {{\bf z}_\sigma}^{4\over {n - 2}}\cdot w_{{\bf z}_\sigma}  \right]\,,\ \ \ \ \ \\[0.15in]
(A.2.3) \ \ \ \  \ \cdot \cdot \cdot  \cdot \cdot \cdot  \cdot \cdot \cdot  \cdot \cdot \cdot \  \ \Longrightarrow \ \ |\,\lambda \cdot {\bf I}|
 & \le &  C \cdot \Vert\, w_{ {\bf z}_\sigma}\Vert_\btd^{{n\,+\,2}\over {n \,-\, 2}}
 \ \approx \ O\left( {1\over {{\bf d}^{{{n\,+\,2}\over 2} \cdot {{n\,+\,2}\over {n\,-\,2}} }}}\right).
 \end{eqnarray*}
  Recall that $\,\lambda \ = \ \sqrt{\lambda_1 \cdot \lambda_2\,}\,.\ $ Recall also (2.13)\,.\, Likewise
$$
|\,\lambda \cdot {\bf II}\,| \  \le   \  C \cdot \Vert\, w_{ {\bf z}_\sigma}\Vert_\btd^{{n\,+\,2}\over {n \,-\, 2}} \cdot
\Vert\,(\lambda\cdot \btd)\,w_{ {\bf z}_\sigma}\Vert_\btd \ \approx\ \ O\left( {1\over {{\bf d}^{{{n\,+\,2}\over 2} \ + \ 2 }}}
\right). \leqno (A.2.4)
$$

The following (highly coupled) terms are fine.
  \begin{eqnarray*}
 (A.2.5) \ \  \  \  \ |\,{\bf III}\,| & \le   & C \cdot   \left( \ \left[\,  \int_{\R^n} |\,{\tilde c}_n \cdot  H\,|^{{2n}\over {n + 2}} \cdot {\bf z}_\sigma^{{2n}\over {n - 2}} \right]^{{n + 2}\over 2n}\cdot \Vert  \,w_{ {\bf z}_\sigma} \Vert_\btd \right)      \cdot \left(\max \ \left\{ {1\over {\lambda_1}}\,, \ {1\over {\lambda_2}}  \right\} \right) \\[0.15in]
 & \ & \ \ \  \ \ + \ \,
 C \cdot    \left( \ \left[\,  \int_{\R^n} |\,{\tilde c}_n \cdot H\,|^{{2n}\over {n - 2}} \cdot {\bf z}_\sigma^{{2n}\over {n - 2}}
 \right]^{{n - 2}\over 2n}\cdot \Vert  \,w_{ {\bf z}_\sigma} \Vert_\btd^2 \right)
 \cdot \left(\max \ \left\{ {1\over {\lambda_1}}\,, \ {1\over {\lambda_2}}  \right\} \right)\ \ \ \ \ \ \ \  \ \ \ \\[0.2in]
 & \le &   C \cdot  {1\over { {\bf d}^{\,2 \cdot {{ n\,+\,2 }\over 2} \ -\,o\,(1\,)  } }}
 \cdot \left(\max \ \left\{ {1\over {\lambda_1}}\,, \ {1\over {\lambda_2}}  \right\} \right)\ .\\[0.3in]
 (A.2.6) \  \ \ \  \ \ |\,{\bf IV}\,| & \le &  C \cdot   \left( \ \left[\,  \int_{\R^n} |\,{\tilde c}_n \cdot \,H\,|^{{2n}\over {n + 2}}
 \cdot {\bf z}_\sigma^{{2n}\over {n - 2}} \right]^{{n + 2}\over 2n}\cdot \Vert
 \btd\,w_{ {\bf z}_\sigma} \Vert_\btd \right) \ + \ C   \cdot \Vert
 \btd\,w_{ {\bf z}_\sigma} \Vert_\btd   \cdot \Vert   \,w_{ {\bf z}_\sigma} \Vert_\btd^2 \ \ \ \ \ \ \  \ \ \   \\[0.2in]
 & \le &   C \cdot  {1\over { {\bf d}^{\,  {{ n\,+\,2 }\over 2} \ + \ 2  \ -\,o\,(1\,) } }}
  \cdot \left(\max \ \left\{ {1\over {\lambda_1}}\,, \ {1\over {\lambda_2}}  \right\} \right)\ .\\[0.3in]
 (A.2.7) \  \ \ \ \ \ \  |\,{\bf V}\,| & \le   & C \cdot {1\over { {\bf d}^{ {{n + 2}\over 2} \,-\,o\,(1)} }} \cdot
 \Vert \,\btd\,w_{ {\bf z}_\sigma} \Vert_\btd\   \le \    C \cdot  {1\over { {\bf d}^{\,  {{ n\,+\,2 }\over 2} \ + \ 2  \ -\,o\,(1\,) } }} \cdot \left(\max \ \left\{ {1\over {\lambda_1}}\,, \ {1\over {\lambda_2}}  \right\} \right)\ .
 \end{eqnarray*}
 Moreover

 $$
 {{n\,+\,2}\over 2} \,+\, 2 \ >  \ n \ - \ 2 \ \ \Longleftrightarrow \ \ 10 \ > \ n \  \ \Longrightarrow\ \ n \ =\ 6\,, \ 7\,, \ 8 \ \ \&\ \ 9\,, \leqno (A.2.8)
 $$

 $$
 {{n\,+\,2}\over 2} \cdot {{n\,+\,2  }\over {n\,-\,2   }} \ > \ n \ - \ 1  \ \ \Longleftrightarrow \ \ 10\,n  \ >  \ n^2
 \ \ \Longleftrightarrow \ \ n \ =\ 6\,, \ 7\,, \ 8 \ \ \&\ \ 9\,.\leqno (A.2.9)
 $$
Combining the estimates we obtain (4.4) and (4.5)\,,\, where the condition $\,10 \ > \ n \ \ge \ 6\,$ is used. \,.\qed

\newpage

{\bf \large \S\,A.3.} \ \ {\bf \large Proof of the Reduction Lemma 4.7.}\\[0.2in]
We largely follow \cite{2nd}\,.\\[0.2in]
 {\bf Lemma A.3.1 (Reduction Lemma)\,.} \ \  {\it In} $\,\R^n\,,\,$ $n \ \ge \ 3$\,,\,  {\it consider a homogeneous polynomial $\,{\cal Q}_{\ell}\,$}  {\it of even degree\,} $\,\ell \ \le \ n \,-\, 1$\,.\,  {\it  We have}
 $$
\int_{\R^n} \   {\cal Q}_{\ell}\, (y) \cdot \left( {1\over {1+ |\,y|^{\,2}}} \right)^{\!n} \,d y \ = \ {{J_n}\over {   \ell \cdot (\ell \,-\,2) \cdot \cdot \cdot \,2\cdot 1 }}\cdot  [\,\Delta_{\,o}^{\!(h_\ell)}\ Q_\ell\,]\ \ \ \ \ \ (h_\ell \ = \ \ell/\,2)\,, \leqno (A.3.2)
$$
 {\it where}
 $$
 J_n \ = \ \int_{R^n} y_1^2 \cdot \cdot \cdot y_{h_\ell}^2 \cdot \left( {1\over {1+ |\,y|^{\,2}}} \right)^{\!n} \,d y \ \ \ \ \ \ \ \  [\,y \ = \ (y_1\,, \cdot \cdot \cdot\,, \ y_\ell\,, \ \cdot \cdot \cdot\,, \ y_n)\,]\,.\leqno (A.3.3)
 $$
{\it Likewise, }
$$
\int_{\R^n} \   {\cal Q}_{\ell}\, (y) \cdot \left( {1\over {1+ |\,y|^{\,2}}} \right)^{\!n\,+\,1} \,d y \ = \ {{J_{n\,+\,1}}\over {   \ell \cdot (\ell \,-\,2) \cdot \cdot \cdot \,2\cdot 1 }}\cdot [\, \Delta_{\,o}^{\!(h_\ell)}\ Q_\ell\,]\ , \leqno (A.3.4)
$$
 {\it where}
 $$
 J_{n\,+\,1} \ = \ \int_{R^n} y_1^2 \cdot \cdot \cdot y_{h_\ell}^2 \cdot \left( {1\over {1+ |\,y|^{\,2}}} \right)^{\!n\,+\,1} \,d y \ \ \ \ (\,\Longrightarrow \ \  J_{n\,+\,1} \ < \ J_n\,)\,. \leqno (A.3.5)
 $$
{\it Here}
\begin{eqnarray*}
  \Delta_o^{\!(h_\ell)} \, {\cal Q}_\ell\, (  y\,) & = &  \Delta_o\ (\,\cdot \cdot \cdot\, [\,\Delta_o \,[\,\Delta_o  \   {\cal Q}_\ell\,( y\,)\,]\,]\,) )\ .
\\[0.1in]
& \ & \leftarrow \  \, h_\ell \ \,{\mbox{times}} \ \,\rightarrow
\end{eqnarray*}

 \vspace*{0.15in}

{\bf Proof.} \ \ We observe that, as $\,\ell \ \le\  n - 1\,,\,$   the integrals in (A.3.2) and (A.3.3) are  absolutely convergent. Keeping the notation $\,y = (y_{|_1}\,,\, \cdot \cdot \cdot\, \ y_{|_n} ) \in \R^n\,$,\, consider a typical term in $\,{\cal Q}_{\ell}\ $:
$$ \ \
y_{|_1}^{\alpha_1} \cdot  y_{|_2}^{\alpha_2}   \cdot \cdot \cdot \   y_{|_{n  }}^{\alpha_{n  }}\,, \ \ \ \ \ \ \ \ {\mbox{where}} \ \ \alpha_j \,\ge\, 0 \ \ \  {\mbox{and}}  \ \ \ \sum_{j = 1}^n \alpha_j  \ = \ \ell \ (\,\le \ n - 1\,)\,.   \leqno (A.3.6)
$$
If one of the indices (say, $\alpha_j$) is an odd natural number, via symmetry, we have
$$
\!\!\!\!\!\!\!\!\!\! \ \ \ \ \   \int_{\R^n} \  [\  y_{|_1}^{\alpha_1} \cdot  y_{|_2}^{\alpha_2}   \cdot \cdot\, y_{|_j}^{\alpha_j} \cdot \cdot \   y_{|_{n  }}^{\alpha_{n  }}] \cdot  \left( {1\over {1+ |\,y|^{\,2}}} \right)^{\!n} \,d y \ = \ 0\ \ \ \ \     (\alpha_j \ \ {\mbox{is\ \  odd}}\,)\,.  \leqno (A.3.7)
$$
Direct calculation also shows that in this situation
$$
 \ \ \ \ \ \ \ \ \ \ \ \ \ \ \ \ \ \ \ \  \ \ \   \ \ \   \ \ \  \ \ \Delta_{\,o}^{\!({h_\ell})}\, \left[\ y_{|_1}^{\alpha_1} \cdot  y_{|_2}^{\alpha_2}   \,\cdot \cdot\  y_{|_j}^{\alpha_j} \cdot \cdot \     y_{|_{n  }}^{\alpha_{n  }} \right] \ = \ 0 \ \ \ \ \ \ \  \ \ \ \ \ \ \ \  \ \ \ \ \ \ \ \ \ \ \  \ \ \ (\alpha_j \ \ {\mbox{is\ \  odd}}\,)\,.
$$
(\,Recall that $\,\Delta_o^{\!(h_\ell)}\, {\cal Q}_\ell\,$ is a number as \,$\ell$\, is even and $\,h_\ell \,=\, \ell/\,2$\,.\,)
Thus we are left with the case where any one index in (A.3.6) is an  even natural number  or zero. Let us introduce the following notion\,: a multi\,-\,index

\vspace*{-0.15in}

$$\ \ \  \ \ \ \ \ \
\alpha = (\alpha_1\,,\, \ \alpha_2\,,\, \ \cdot \cdot \cdot\,, \ \alpha_n) \ \ \ \ \ \ \ \ \  \ \ \ \ \ \ \bigg( \ |\, \alpha| \ = \ \sum_{j = 1}^n \,\alpha_j \ = \ \ell > 0 \, \bigg)
$$

\vspace*{-0.05in}

is {\it even} if each $\alpha_j\ \, (\,1 \,\le\, j \,\le\, n\,)\,$ is either  an even natural number or zero\,.\, With respect to this, the simplest case that can happen to the integral in (A.3.2) is
$$
J \, := \, \int_{\R^n}    \,[\,y_{|_1}^2 \,\cdot \cdot \cdot  \  y_{|_{h_\ell}}^2 \,]\cdot \left( {1\over {1+ r^{\,2}}} \right)^{\!n} \, d y \ \ \ \ \ \ \ (\,r \ = \ |\,y\,|\,)\,.\leqno (A.3.8)
$$
We seek to reduce other even multi\,-\,index cases to that in (A.3.8)\,.\,
As $\,\ell \,\le\, n - 1\,$,\, we  consider the following.
$$
y_{|_1}^{\,k + 2} \cdot  y_{|_2}^{\alpha_2}   \cdot \cdot \cdot \   y_{|_{n - 1}}^{\,\alpha_{n - 1}}\,, \ \ \ \ \ \ \ {\mbox{where}} \ \ \ \ \alpha \,= \ (k + 2\,,\, \ \alpha_2\,,\, \ \cdot \cdot \cdot\,, \ \alpha_{n-1}\,, \ 0) \ \ \ {\mbox{is \ \ even}}\,.
$$
Here $\,k \,\ge \,2\,$ is an even number\,.\,
By using  Fubini's theorem and integration by parts formula\,,\, we obtain  the following reduction formula.
\begin{eqnarray*}
(A.3.9)\ \ \ \ & \ & \ \ \ \int_{\R^n} \ \,  y_{|_1}^{\,k + 2} \cdot  \left[ \ \ y_{|_2}^{\alpha_2}   \cdot \cdot \cdot \   y_{|_{n - 1}}^{\alpha_{n - 1}}\, \right]   \cdot  \left( {1\over {1+ r^{\,2}}} \right)^{\!n} \, d y\\[0.15in]
 &\ &  \!\!\!\!\!  = \ (k + 1)    \int_{\R^n}   \ \,   y_{|_1}^{\,k  }  \cdot y_{|_n}^{\,2  }  \cdot \left[ \ \ y_{|_2}^{\alpha_2}   \cdot \cdot \cdot \   y_{|_{n - 1}}^{\alpha_{n - 1}}\, \right]  \cdot  \left( {1\over {1+ r^{\,2}}} \right)^{\!n} \, d y \ \ \ \     \mfor 2 \,\le\, k \,\le\, n - 3\,. \ \ \ \ \ \ \ \ \ \
\end{eqnarray*}
  See \S\,A.3 below\,.\bk
In view of (A.3.9), we introduce another notation.  For an integer $\,m\,\ge\, 0$\,,\, define
\begin{eqnarray*}
(A.3.10) \ \ \ \ \ \ \ \ \ \  m\,!_{-2} & = &  1 \ \ \ \ {\mbox{if}} \ \ m = 0\ \ \ {\mbox{or \  }} 2\,; \ \ \ \ \ m\,!_{-2} \ =  \ 0 \ \ \ \ {\mbox{if}} \ \ m \ \ {\mbox{is \ \ odd}}\,;
\\[0.1in]
m\,!_{-2} & =  & ( m - 1) \,( m - 3)\,(m - 5) \cdot \cdot \cdot \ \cdot \,3 \cdot 1 \ \ \ \ \ \ \ {\mbox{if}} \ \ m \,\ge\, 4 \ \ \, {\mbox{is \ \ even}}\,. \ \ \ \ \ \  \ \ \ \ \ \  \ \ \ \ \
\end{eqnarray*}
Via the vanishing formula (A.3.7) and the reduction formula (A.3.9)\,,\, we have
$$ \int_{\R^n} \  [\ y_1^{\,\alpha_1} \cdot \cdot \cdot  \,  y_n^{\,\alpha_n} \,]\cdot  \left( {1\over {1+ r^{\,2}}} \right)^{\!n} \, d y \ = \  (\alpha_1)\,!_{-2} \,\times \,\cdot \cdot \cdot \,\times (\alpha_n)\,!_{-2}\cdot J\ .\leqno (A.3.11)
$$
On the other side,  calculation shows that
$$
B:= \Delta_{\,o}^{\!( h_\ell)}\, \left[\ y_1^2 \cdot y_2^2 \ \cdot \cdot \cdot \ y_{h_\ell}^2\ \right]  \ =  \ \ell \, (\ell - 2) \,(\ell - 4) \,\cdot  \cdot \cdot \,2 \cdot 1\,. \leqno (A.3.12)
$$

\vspace*{0.15in}

{\bf Claim.} \ \ Let $\,\alpha_2\,, \ \cdot \cdot \cdot\,, \ \alpha_{n - 1}\,$   be even natural numbers or zero,\, and
$$
\ell \ =\ ( k + 2) + {\alpha_2}  \ + \ \cdot \cdot \cdot \ + \   \alpha_{n - 1}\,, \ \ \ \ {\mbox{where}} \ \ k \ \ge \ 2 \ \ {\mbox{is \ \ an \ \ even \ \ integer}}\,. \leqno (A.3.13)
$$
Then
\begin{eqnarray*}
(A.3.14) \ \  \ \ \ \ \  \Delta_{\,o}^{\!( h_\ell)}  \left\{\,   y_{|_1}^{\,k + 2}  \cdot \left[ \ \ y_{|_2}^{\alpha_2}   \cdot \cdot \cdot \   y_{|_{n - 1}}^{\alpha_{n - 1}}\, \right]   \, \right\} \, = \, (k + 1) \cdot \!\Delta_{\,o}^{\!( h_\ell)}\, \left\{\, y_{|_1}^{\,k  }  \cdot  y_{|_n}^{\,2  }  \left[ \ \ y_{|_2}^{\alpha_2}   \cdot \cdot \cdot \   y_{|_{n - 1}}^{\alpha_{n - 1}}\, \right]   \right\}. \ \ \ \ \
\end{eqnarray*}
Refer to \S\,A.3\,.\, \bk
Thus using (A.3.14) repeatedly, we are led to
$$
\Delta_o^{\!(h_\ell)} \ \left[\,  y_{|_1}^{\alpha_1}   \cdot \cdot \cdot \   y_{|_{n }}^{\alpha_{n }}\,  \right] \ =  \  (\alpha_1)\,!_{-2} \,\times \,\cdot \cdot \cdot \,\times (\alpha_n)\,!_{-2}\cdot B\,. \leqno (A.3.15)
$$
Using the linearity of the operations, together with (A.3.11) and (A.3.15), we obtain

\vspace*{-0.25in}

$$
\int_{\R^n} \   {\cal Q}_{\ell}\, (y) \cdot \left( {1\over {1+ |\,y|^{\,2}}} \right)^n\,d y \ =  \ {J\over B}\cdot  [\ \Delta_{\,o}^{\!(h_\ell)}\ {\cal Q}_{\ell}\,]\, . \leqno (A.3.16)
$$
In particular, we establish (A.3.2)\,.\,
The case for $\,n\,+\,1\,$ is a direct modification of the above argument.  {{\hfill {$\rlap{$\sqcap$}\sqcup$}}

 \vspace*{0.3in}

{\bf Verification of\ } (A.3.9) \,{\bf and\,} (A.3.14).\\[0.1in]
 Refer to (A.3.9) for the notation we use.
\begin{eqnarray*}
& \ & \\[-0.25in]
& \ & \int_{\R^n} \ \, (  {\mbox{\,terms \ \ without  \ }} y_{|_1} \ \,\& \ \,y_{|_n}  ) \cdot \ y_{|_1}^{\,k + 2} \cdot \left( {1\over {1+ r^{\,2}}} \right)^{\!n} \, d y \ \ \ (\,{\mbox{absolute \ \ convergence}})\\[0.15in]
& = & \int_{-\infty}^\infty \cdot  \cdot \cdot \ \int_{-\infty}^\infty \,(  {\mbox{\,terms \ \ without  \ }} y_{|_1} \ \,\& \ \, y_{|_n}  ) \ \times \  \ \ \ \ \ \ \ \ \ \ \ \ \ \ \ \    (\,{\mbox{Fubini's Theorem}}) \\[0.1in]
 & \ &  \!\!\!\!\!\!\!\!\!\times \left[ \ \int_{-\infty}^\infty\int_{-\infty}^\infty \,y_{|_1}^{\,k + 2} \cdot  \left( {1\over {1 + [\,y_{|_2}^{\,2} + \cdot \cdot \cdot + y_{|_{n-1}}^2 \,] + y_{|_1}^{\,2} + y_{|_n}^{\,2}}} \right)^{\!n}  dy_{|_1} \,dy_{|_n} \right]dy_{|_2} \cdot \cdot \cdot dy_{|_{n - 1}}\\[0.15in]
& = & \int_{-\infty}^\infty \cdot  \cdot \cdot \ \int_{-\infty}^\infty \,(  {\mbox{\,terms \ \ without  \ }} y_{|_1} \ \& \ y_{|_n}  ) \ \times \  \\[0.1in]
 & \ & \!\!\!\!\!\!\!\times\, \left[\  \int_0^\infty\int_0^{2 \,\pi} \,y_{|_1}^{\,k + 2}  \cdot   \left( {1\over {1 + [\,y_{|_2}^{\,2} + \cdot \cdot \cdot + y_{|_{n - 1}}^{\,2}\,] + \rho^{\,2}}} \right)^{\!n} \cdot \rho \  d\theta \, d \rho \right]dy_{|_2} \cdot \cdot \cdot dy_{|_{n-1}}\\
 & \ & \ \ \ (\,{\mbox{polar \ \ coordinates \ \ on \ \ }} \R^2\,,   \ \   \rho^2 \ = \ y_1^2 \,+\, y_n^2\,,\, \ \ y_1 \, = \ \rho\cdot \sin\,\theta\,, \ \ y_n \, = \ \rho\cdot \cos\,\theta )\\[0.15in]
& = & \int_{-\infty}^\infty \cdot  \cdot \ \int_{-\infty}^\infty \,(  {\mbox{\,terms \ \ without  \ }} y_{|_1} \ \,\&\, \ y_{|_n}  ) \ \times  \   \\[0.1in]
 & \ & \!\!\!\!\!\!\!\!\!\!\!\!   \times\, \left[ \ \int_0^\infty\int_0^{2 \,\pi} \,\rho^{\,k + 3}\, (\sin^{\,k + 2} \theta)  \cdot   \left( {1\over {1 + [\,y_{|_2}^{\,2} + \cdot \cdot \cdot + y_{|_{n-1}}^{\,2}\,] + \rho^{\,2}}} \right)^{\!n}  \cdot  d\theta \, d \rho\, \right]\,dy_{|_2} \cdot \cdot \cdot dy_{|_{n-1}}\\[0.15in]
 & = & \int_{-\infty}^\infty \cdot  \cdot \ \int_{-\infty}^\infty \,(  {\mbox{\,terms \ \ without  \ }} y_{|_1} \ \,\&\, \ y_{|_n}  ) \ \times  \   \\[0.1in]
& \ & \!\!\!\!\!\!\!\!\!\!\!\!   \times\,\left[  \!\int_0^\infty \!\!\rho^{\,k + 3}  \cdot   \left( {1\over {1 + [\,y_{|_2}^{\,2} + \cdot \cdot \cdot + y_{|_{n-1}}^{\,2}\,] + \rho^{\,2}}} \right)^{\!n}  \cdot  \left\{   \int_0^{2 \,\pi} \!\! \sin^{\,k \,+\, 2} \theta  \,d\theta \right\}  d \rho \, \right]\,dy_{|_2} \cdot \cdot\,   dy_{|_{n-1}}\ .
\end{eqnarray*}
Here ``\,(terms without $\,y_{|_1} \ \,\&\, \ y_{|_n}\,$)" is a polynomial on $y_{|_2}\,,\, \cdot \cdot \cdot\,, \ y_{|_{n-1}}$\,,\, having  sufficiently low degree so that the integral is absolutely convergent.
Likewise,
\begin{eqnarray*}
& \ & \int_{-\infty}^\infty \cdot  \cdot \ \int_{-\infty}^\infty \,(  {\mbox{\,terms \ \ without  \ }} y_{|_1} \ \,\&\, \ y_{|_n}  ) \ \times \   \\
 & \ &   \!\!\!\!\!\!\!\!\!\!\!\!\!\!\!\!\times \, \left[ \ \int_{-\infty}^\infty\int_{-\infty}^\infty \,y_{|_1}^{\,k}\,y_{|_n}^{\,2} \cdot  \left( {1\over {1 + [\,y_{|_2}^{\,2} + \cdot \cdot \cdot \,+\,  y_{|_{n-1}}^2\,] + y_{|_1}^{\,2} \,+\, y_{|_n}^{\,2}}} \right)^{\!n} \cdot   dy_{|_1} \,dy_{|_n} \right]dy_{|_2} \cdot \cdot \cdot dy_{|_{n-1}}\\[0.1in]
& = & \int_{-\infty}^\infty \cdot  \cdot \ \int_{-\infty}^\infty \,(  {\mbox{\,terms \ \ without  \ }} y_{|_1}\ \,\&\,\ y_{|_n}  ) \ \times  \   \\
 & \ & \!\!\!\!\!\!\!\!\!\!\!\!\!\!\!\!   \times   \left[\,  \int_0^\infty \!\!\!\rho^{\,k + 3}    \cdot  \left( \!{1\over {1 + [\,y_{|_2}^{\,2} \,+ \cdot \cdot   \,+ \, y_{|_{n-1}}^2] + \rho^{\,2}}} \!\right)^{\!n} \cdot \left\{ \int_0^{2 \,\pi} \!(\sin^k \theta)\,(\cos^{\,2} \theta) \,d\theta \right\} \, d \rho \right]dy_{|_2} \cdot \cdot  \, dy_{|_{n-1}}\,.
\end{eqnarray*}

A direct calculation using integration by parts shows that
$$
\int_0^{2 \,\pi} (\sin^{\,k + 2} \theta)\, d\theta \ = \ - \int_0^{2 \,\pi} (\sin^{\,k + 1} \theta)\, d\,[\,\cos\,\theta] \ = \  (k + 1) \int_0^{2 \,\pi} (\sin^k \theta)\,(\cos^{\,2} \theta) \,d\theta\,.
$$
Hence we deduce  (A.3.9). To show (A.3.14), let [\,refer to (A.3.14)\,]
$$
{\cal o} \ = \ {\alpha_2}  \ + \ \cdot \cdot \cdot \ + \   \alpha_{n - 1}\ . \leqno (A.3.17)
$$
We demonstrate how to use induction on $\,{\cal o}\,$ to prove the assertion. Recall that
$$
 \ell \,\in\, [\,0, \ n- 2) \ \ \ {\mbox{is \ \ even}}\,.
$$
{\bf (I)} \ \ When $\,{\cal o} \ = \ 0$\,,\, the term specified in (A.3.14) is a constant. We have
\begin{eqnarray*}
\Delta_o^{\!(h_\ell)} \  y_{|_1}^{\,k + 2}  & = & (k + 2) \, (k + 1) \cdot \cdot \cdot \, 3 \cdot 2 \cdot 1\,;
 \\[0.15in]
 \Delta_o  \ \left\{(k + 1) \, y_{|_1}^{\,k  } \,y_{|_n}^2  \right\} & = & (k + 1) \,k \, (k - 1) \, y_{|_1}^{\,k - 2  }\,y_{|_n}^2 \ + \ 2 \, (k + 1) \, y_{|_1}^{\,k  }\,,\\[0.15in]
 \Delta_o^{\!(2)}  \ \left\{(k + 1) \, y_{|_1}^{\,k  } \,y_{|_n}^2  \right\} & = & (k + 1) \,k \, (k - 1) \, (k - 2) \, (k - 3) \, y_{|_1}^{\,k - 4  }\,y_{|_n}^2 \ +\\[0.1in]
  & \ &  \ \ \ \ \ + \ 2 \times 2\, (k + 1) \,k \, (k - 1) \, y_{|_1}^{\,k - 2  }\,,\\
  & \cdot & \\
  & \cdot & \\
  & \cdot &\\
  \end{eqnarray*}

\vspace*{-0.5in}

  \begin{eqnarray*}
  \Delta_o^{ \, (h_\ell - 2)}  \ \left\{(k + 1) \, y_{|_1}^{\,k  } \,y_{|_n}^2  \right\} & = & (k + 1) \,k \, (k - 1) \, (k - 2) \, (k - 3) \, \cdot \cdot \cdot \, 3\cdot \, y_{|_1}^{\,2 }\,y_{|_n}^2 \ +\\[0.1in]
  & \ & \ \ \ \ \ \ \ \ \ + \ 2 \,\left[ \,{h_\ell} - 2 \right] (k + 1) \,k \, (k - 1) \cdot \cdot \cdot \, 5  \,\cdot y_{|_1}^{\,4  }\,,\\[0.1in]
  \Delta_o^{ \, (h_\ell - 1)}  \ \left\{(k + 1) \, y_{|_1}^{\,k  } \,y_{|_n}^2  \right\} & = & (k + 1) \,k \, (k - 1) \, (k - 2) \, (k - 3) \, \cdot \cdot \cdot \, 3\cdot 2 \cdot  \, [\,y_{|_1}^{\,2 }+ y_{|_n}^2] \ +\\[0.1in]
  & \ & \ \ \ \ \ \ \ \ \ \ + \ 2 \,\left[ \,{h_\ell} - 2\, \right] (k + 1) \,k \, (k - 1) \cdot \cdot \cdot \, 5 \cdot 4 \cdot 3\, y_{|_1}^{\,2 }\,,\\[0.1in]
  \Delta_o^{ \, h_\ell}  \ \left\{(k + 1) \, y_{|_1}^{\,k  } \,y_{|_n}^2  \right\} & = & (k + 1) \,k \, (k - 1) \, (k - 2) \, (k - 3) \, \cdot \cdot \cdot \, 3\cdot 2 \cdot  1 \cdot \, 4\ + \\[0.1in]
  & \ & \ \ \ \ \ + \  \left[\, \ell \,-\, 4 \right] (k + 1) \,k \, (k - 1) \cdot \cdot \cdot \, 5 \cdot 4 \cdot 3\,\cdot 2 \cdot 1\\[0.1in]
  & = & (k + 1) \,k \, (k - 1) \, (k - 2) \, (k - 3) \, \cdot \cdot \cdot \, 3\cdot 2 \cdot  1 \!\cdot\! [\,\ell - 4 + 4\,]\\[0.1in]
   & = & (k + 2) \,(k + 1) \,k \, (k - 1) \, (k - 2) \, (k - 3) \, \cdot \cdot \cdot \, 3\cdot 2 \cdot  1  \ \ \ \ \ \\[0.1in]
    & \ & \hspace*{1.55in}  (\,{\mbox{as}} \ \ \ell \,=\, k + 2  \ \ {\mbox{in \ this \ case}})\,.
\end{eqnarray*}
Hence the case $\,{\cal o}\, =\, 0\,$ is settled.

\vspace*{0.2in}

{\bf (II)} \ \ As an induction hypothesis, suppose that
$$
\Delta_o^{\! (h_\ell)} \ \left\{ y_{|_1}^{\,k + 2} \cdot [ \,\cdot \cdot \cdot \ {\mbox{degree}} \ = {\cal o} \ \cdot \cdot \cdot \,] \right\} = \Delta_o^{\!(h_\ell)} \ \left\{(k + 1) \, y_{|_1}^{\,k  } \,y_{|_n}^2 \times [ \,\cdot \cdot \cdot \ {\mbox{degree}} \ = \ {\cal o} \ \cdot \cdot \cdot \,]  \ \right\} \leqno (A.3.18)
$$
holds for $\,\ell = ( k + 2) \ + \ {\cal o}\,,\,$ where $\,k \, \ge \,2 \,$ (variable), but $\,{\cal o} \,>\, 0\,$ (fixed). We continue to use the notations above and {\it there is no\,  $\,y_{|_1}\,$ or \,$\,y_{|_n}\,$ inside the homogeneous polynomial denoted by\,} $\,[ \,\cdot \cdot \cdot \ {\mbox{degree}} \ = \,{\cal o}   \ \cdot \cdot \cdot \,]\,  $\,.\, Let us go on to show
\begin{eqnarray*}
(A.3.19)\ \ \ \ \ \ \ \ \ \ \ \  & \ & \Delta_o^{\! (h_\ell)} \ \left\{ y_{|_1}^{\,k + 2} \cdot [ \,\cdot \cdot \cdot \ {\mbox{degree}} \ =\, {\cal o} + 2 \ \cdot \cdot \cdot \,] \right\}\\[0.1in] &  = & \Delta_o^{\! (h_\ell)} \ \left\{(k + 1) \, y_{|_1}^{\,k  } \,y_{|_n}^2 \cdot [ \,\cdot \cdot \cdot \ {\mbox{degree}} \ = \,{\cal o} + 2 \ \cdot \cdot \cdot \,]  \ \right\}\ , \ \ \ \ \ \ \ \ \ \ \ \ \ \ \ \ \ \ \ \ \ \ \ \ \ \ \ \ \ \ \ \ \ \ \ \ \ \ \ \ \ \ \ \ \ \ \ \
\end{eqnarray*}
where $k \,\ge\, 2$ is even. Let us find the first Laplacians\,:
\begin{eqnarray*}
(A.3.20) \!\!& \ & \Delta_o  \, \left\{ y_{|_1}^{\,k + 2} \cdot [ \,\cdot \cdot \cdot \ {\mbox{degree}} \ = \  {\cal o} + 2 \ \cdot \cdot \cdot \,] \right\}\\[0.1in]
& = & (k + 2) \, (k + 1) \ y_{|_1}^{\,k  } \cdot [ \,\cdot \cdot \cdot \ {\mbox{degree}} \ =  \ {\cal o} + 2 \ \cdot \cdot \cdot \,] \\[0.1in]
 & \ & \ \ \ \ \ \ \ \ \ \ \ \ \ \ \ \ + \  y_{|_1}^{\,k + 2} \cdot \left\{ \Delta_o \ [ \,\cdot \cdot \cdot \ {\mbox{degree}} \ = \ {\cal o} + 2 \ \cdot \cdot \cdot \,] \right\}\\[0.1in]
& = & k \, (k + 1) \ y_{|_1}^{\,k  } \cdot [ \,\cdot \cdot \cdot \ {\mbox{degree}} \ =  \ {\cal o} + 2 \ \cdot \cdot \cdot \,] \\[0.1in]
& \ &  \ \ \  \ + \ 2 \, (k + 1) \ y_{|_1}^{\,k  } \cdot [ \,\cdot \cdot \cdot \ {\mbox{degree}} \ =  \ {\cal o} + 2 \ \cdot \cdot \cdot \,]\\[0.1in]
 & \ & \ \ \ \ \ \ \ \ \ \ \ \ \ \ \ \ + \  y_{|_1}^{\,k + 2} \cdot \left\{ \Delta_o \ [ \,\cdot \cdot \cdot \ {\mbox{degree}} \ =  \ {\cal o} + 2 \ \cdot \cdot \cdot \,] \right\}\\[0.1in]
 & \ & \ \ \ \ \ \ \ \ \ \  \ \ \ \ \ \ \ \ \ \  \ \ \ \ \ \ \  (\ \leftarrow \ \ \ \ \ \uparrow \ \ {\mbox{degree}} \ =  \ {\cal o} \ \ \ \ \ \rightarrow \ )\,;
 \\[0.22in]
& \ &  \Delta_o  \, \left\{(k + 1) \, y_{|_1}^{\,k  } \,y_{|_n}^2 \cdot [ \,\cdot \cdot \cdot \ {\mbox{degree}} \ =  \ {\cal o} + 2 \ \cdot \cdot \cdot \,]  \ \right\}\\[0.1in]
 &= & \  (k + 1)\,k\,(k - 1) \,y_{|_1}^{\,k  - 2} \,y_{|_n}^2 \cdot [ \,\cdot \cdot \cdot \ {\mbox{degree}} \ =  \ {\cal o} + 2 \ \cdot \cdot \cdot \,]\\[0.1in]
 & \ & \ \ \ \ \ \ \ \ \ \ \ \  + \ 2\,(k + 1) \,y_{|_1}^k \cdot [ \,\cdot \cdot \cdot \ {\mbox{degree}} \ = \ {\cal o} \ +\  2 \ \cdot \cdot \cdot \,]\\[0.1in]
 & \ & \ \ \ \ \ \ \ \ \ \ \ \ \ \ \ \  \ \ \ \ \ \  + \  (k + 1) \,y_{|_1}^{\,k  } \,y_{|_n}^2 \cdot \left\{ \Delta_o \ [ \,\cdot \cdot \cdot \ {\mbox{degree}} \ = \ {\cal o} + 2 \ \cdot \cdot \cdot \,] \right\}\\[0.1in]
 & \ & \ \ \ \ \ \ \ \ \ \ \ \ \ \ \  \ \ \ \ \ \ \ \ \ \  \ \ \ \ \ \ \  \ \ \ \ \ \ \ \ \   \ \ \  (\ \leftarrow \ \ \ \ \ \uparrow \ \ {\mbox{degree}} \ = \  {\cal o} \ \ \ \ \ \rightarrow \ )\,,\\[0.1in]
& \ & \ \ \ \ \ \ \ \ \ \ \ \  \ \ \ \ \ \   [\ {\mbox{observe \ \ that \ \ }}  (k + 2)\, (k + 1) - 2\, (k + 1) \ = \ k \, (k + 1)\,]\,. \ \ \ \ \ \ \ \ \ \ \ \ \ \ \ \ \ \ \ \ \ \ \ \ \ \ \ \ \ \ \ \ \
\end{eqnarray*}
Via the induction hypothesis, the last two terms in the respective expressions are equal. After simplification, to verify (A.3.19), it suffices to  show that
\begin{eqnarray*}
& \ & \Delta_o^{\!(\,h_\ell - 1)} \ \left\{ y_{|_1}^{\,k } \cdot [ \,\cdot \cdot \cdot \ {\mbox{degree}} \ = \   {\cal o} + 2 \ \cdot \cdot \cdot \,]\, \right\}\\[0.1in]
 &= & \Delta_o^{\!(\,h_\ell - 1)} \ \left\{(k - 1) \, y_{|_1}^{\,k - 2  } \ y_{|_n}^2 \cdot [ \,\cdot \cdot \cdot \ {\mbox{degree}} \ = \ {\cal o} + 2 \ \cdot \cdot \cdot \,]  \ \right\}\ .
\end{eqnarray*}
Applying $\,\Delta_o\,$ on the  terms
$$
\left\{ y_{|_1}^{\,k } \cdot [ \,\cdot \cdot \cdot \ {\mbox{degree}} \ = \   {\cal o} + 2 \ \cdot \cdot \cdot \,]\, \right\} \ \   {\mbox{and}} \ \  \left\{(k - 1) \, y_{|_1}^{\,k - 2  } \ y_{|_n}^2 \cdot [ \,\cdot \cdot \cdot \ {\mbox{degree}} \ = \   {\cal o} + 2 \, \cdot \cdot \cdot \,]  \, \right\},
$$
using  similar calculation and cancelation as in (A.3.20), we come down gradually to verify
\begin{eqnarray*}
(A.3.21) \ \ \ \ \ \ & \ & \Delta_o^{\! \left(h_\ell \,-\, {{k - 2}\over 2} \right)} \ \left\{ y_{|_1}^{\,4 } \cdot [ \,\cdot \cdot \cdot \ {\mbox{degree}} \ = \   {\cal o} + 2 \ \cdot \cdot \cdot \,]\, \right\}\\[0.1in]
 &= &\Delta_o^{\! \left(h_\ell \,-\, {{k - 2}\over 2} \right)} \   \left\{3 \, y_{|_1}^{\,2 } \,y_{|_n}^2 \cdot [ \,\cdot \cdot \cdot \ {\mbox{degree}} \ =  \ {\cal o} + 2 \ \cdot \cdot \cdot \,]  \ \right\}\,.\ \ \ \ \ \ \ \ \ \ \ \ \ \ \ \ \ \ \ \ \ \ \ \ \ \ \ \ \ \ \ \ \ \ \ \  \ \ \ \ \ \ \ \ \ \ \ \
\end{eqnarray*}
Note that
$$
\ell = (k + 2) + ( {\cal o} + 2) \ \ \Longrightarrow \ \ {\ell\over 2} \,-\, {{k - 2}\over 2} \ = \ {{4 + ( {\cal o} + 2)  }\over 2}\ .
$$
Apply the Laplacian on the two terms inside the brackets in (A.3.21) and obtain
\begin{eqnarray*}
& \ & \Delta_o  \ \left\{ y_{|_1}^{\,4} \cdot [ \,\cdot \cdot \cdot \ {\mbox{degree}} \ =  \ {\cal o} + 2 \ \cdot \cdot \cdot \,] \right\}\
= \ 4 \cdot 3 \ y_{|_1}^{\,2 } \cdot [ \,\cdot \cdot \cdot \ {\mbox{degree}} \ =  {\cal o} + 2 \ \cdot \cdot \cdot \,]\ +  \\[0.1in]
 & \ & \ \ \ \ \ \ \ \ \ \ \ \ \ \ \ \ + \  y_{|_1}^{\,4} \cdot \left\{ \ \Delta_o \ [ \,\cdot \cdot \cdot \ {\mbox{degree}} \ =  \ {\cal o} + 2 \ \cdot \cdot \cdot \,] \,\right\}\,;\\[0.1in]
 & \ & \ \ \ \ \ \ \ \ \ \  \ \ \ \ \ \ \ \ \ \  \ \ \ \ \ \,    (\ \leftarrow \ \ \ \ \ \uparrow \ \ {\mbox{degree}} \ =  \ {\cal o} \ \ \ \ \ \rightarrow \ )
 \\[0.25in]
& \ &  \Delta_o  \ \left\{3 \, y_{|_1}^{\,2  } \ y_{|_n}^2 \cdot [ \,\cdot \cdot \cdot \ {\mbox{degree}} \ = \ {\cal o} + 2 \ \cdot \cdot \cdot \,]  \ \right\}\\[0.1in]
 &= & \  3 \cdot 2 \cdot 1 \cdot [\ y_{|_1}^2 + y_{|_n}^2\,] \cdot [ \,\cdot \cdot \cdot \ {\mbox{degree}} \ = \ {\cal o} + 2 \ \cdot \cdot \cdot \,]\ + \\[0.1in]
 & \ & \ \ \ \ \  \ \ \ \ \ \ \ \ \ \ + \  3 \,y_{|_1}^{\,2  }\, y_{|_n}^2 \cdot \left\{ \ \Delta_o \ [ \,\cdot \cdot \cdot \ {\mbox{degree}} \ =  \ {\cal o} + 2 \ \cdot \cdot \cdot \,] \, \right\}\\[0.1in]
 & \ & \ \ \ \ \ \ \ \ \ \ \ \ \ \ \  \ \ \ \ \ \ \ \ \ \  \ \, \ \ \ \  (\ \leftarrow \ \ \ \ \ \uparrow \ \ {\mbox{degree}} \ =  {\cal o} \ \ \ \ \ \rightarrow \ )\,.
\end{eqnarray*}
Again we apply the induction hypothesis to cancel the last  term in each expression above.  Apply the Laplacian again and obtain

\newpage

\begin{eqnarray*}
(A.3.22)& \ & \!\!\!\!\!\!\!\Delta_o \ \left\{   4 \cdot 3 \ y_{|_1}^{\,2 } \cdot [ \,\cdot \cdot \cdot \ {\mbox{degree}} \ = \   {\cal o} + 2 \ \cdot \cdot \cdot \,] \ \right\} \\[0.1in]
 &\ & \!\!\!\!\!\!\!\!\!\!\!\!\!\!\!\!\!\!\!\!\!\!\!\!\!\!\!\!=\ 4 \cdot 3 \cdot 2 \cdot  [ \,\cdot \cdot \cdot \ {\mbox{degree}} \ =  {\cal o} + 2 \ \cdot \cdot \cdot \,] \ + \ 4 \cdot 3 \cdot y_{|_1}^2 \cdot \left\{ \Delta_o \ [ \,\cdot \cdot \cdot \ {\mbox{degree}} \ = \ {\cal o}\ +\ 2 \ \cdot \cdot \cdot \,] \, \right\}\ ;\\[0.25in]
(A.3.23)& \ &  \Delta_o \ \left\{   3 \cdot 2 \ [\, y_{|_1}^{\,2 } + y_{|_n}^2\,] \cdot [ \,\cdot \cdot \cdot \ {\mbox{degree}} \ =  {\cal o} + 2 \ \cdot \cdot \cdot \,] \,\right\}\\[0.1in]
 &  = &     3 \cdot 2 \cdot 4 \cdot  [ \,\cdot \cdot \cdot \ {\mbox{degree}} \ =  \ {\cal o}\ + \ 2 \ \cdot \cdot \,\cdot \,] \\[0.1in]
   & \ & \ \ \ \ \ \ \ \ + \  3 \cdot 2 \cdot [\, y_{|_1}^2 \ + \  y_{|_n}^2\,] \cdot \left\{ \Delta_o \ [ \,\cdot \cdot \cdot \ {\mbox{degree}} \ = \  {\cal o} \ +\  2 \ \cdot \cdot \cdot \,] \, \right\}\,. \ \ \ \ \ \ \ \ \ \ \ \ \ \ \ \ \ \ \ \ \ \ \ \ \ \ \ \ \ \ \ \ \ \ \
\end{eqnarray*}
\begin{eqnarray*}
{\mbox{As}} \  \ & \ & 4 \cdot 3 \cdot  y_{|_1}^2\,  \cdot \left\{ \Delta_o \ [ \,\cdot \cdot \cdot \ {\mbox{degree}} \ = \ {\cal o} + 2 \ \cdot \cdot \cdot \,] \ \right\} \\[0.1in]
& = & 2\cdot 3 \cdot  y_{|_1}^2\,  \cdot \left\{ \Delta_o \ [ \,\cdot \cdot \cdot \ {\mbox{degree}} \ = \ {\cal o} \ +\ 2 \ \cdot \cdot \cdot \,] \ \right\} \ + \\[0.1in]
& \ & \ \ \ \ \ \ + \  2\cdot 3 \cdot  y_{|_1}^2\,  \cdot \left\{ \Delta_o \ [ \,\cdot \cdot \cdot \ {\mbox{degree}} \ =  \ {\cal o} + 2 \ \cdot \cdot \cdot \,] \ \right\}\ ,\\[0.25in]
& \ &   3 \cdot 2 \cdot [\, y_{|_1}^2 + y_{|_n}^2\,] \cdot \left\{ \Delta_o \ [ \,\cdot \cdot \cdot \ {\mbox{degree}} \ =  \  {\cal o} + 2 \ \cdot \cdot \cdot \,] \ \right\}\\[0.1in]
 & = &   3 \cdot 2 \cdot  y_{|_1}^2   \cdot \left\{ \Delta_o \ [ \,\cdot \cdot \cdot \ {\mbox{degree}} \ =  \ {\cal o} \ + \ 2 \ \cdot \cdot \cdot \,] \ \right\} \\[0.1in]
   & \ & \ \ \ \ + \ 3 \cdot 2 \cdot  y_{|_n}^2\,  \cdot \left\{ \Delta_o \ [ \,\cdot \cdot \cdot \ {\mbox{degree}} \ = \  {\cal o} + 2 \ \cdot \cdot \cdot \,] \ \right\}\ ,\ \ \ \ \ \ \ \ \ \ \ \ \ \ \ \ \ \ \ \ \ \ \ \ \ \ \ \\[0.1in]
{\mbox{and}} \,& \ & \!\!\!\!\!\!\Delta_o^{\left({{2 \,+\, ( \,{\cal o} \,+\, 2)}\over 2} \right) } \ \left[\  y_{|_1}^2   \cdot \left\{ \Delta_o \ [ \,\cdot \cdot \cdot \ {\mbox{degree}} \ =  \ {\cal o} + 2 \ \cdot \cdot \cdot \,] \ \right\} \right] \\[0.1in]
 &\ &  \!\!\!\!\!\!\!\!\!\!\!\!\!\!\!\!\!\!\!\!\!\!\!\!\!=\ \Delta_o^{\left({{2 \,+\, ( \,{\cal o} \,+\, 2)}\over 2} \right) } \ \left[\  y_{|_n}^2   \cdot \left\{ \Delta_o \ [ \,\cdot \cdot \cdot \ {\mbox{degree}} \ =\   {\cal o} \ + \ 2 \ \cdot \cdot \cdot \,] \ \right\} \right] \   (\leftarrow \   {\mbox{equal \   to \   a \   number}})\,, \ \ \ \ \ \ \ \ \ \ \ \ \ \ \ \ \ \ \ \ \ \ \ \
\end{eqnarray*}
we apply the remaining order of Laplacian on (A.3.22) and (A.3.23), yielding the same numbers. Hence we verify (A.3.20), and so (A.3.21)\,.\, This  completes the induction step.\qed

\newpage

\newpage

{\large \bf \S\,A.\,4. \ \ Interaction terms (4.14)\,.}

\vspace*{0.2in}

We begin with
\begin{eqnarray*}
 & \ & {\partial\over {\partial \xi_{1_j}}} \,\int_{\R^n} ({\tilde c}_n \cdot H) \!\cdot\! \left( \,V_1 \ + \ V_2\, \right)^{{{2n}\over {n\,-\,2}}   }\ = \ {{2n}\over {n\,-\,2}} \cdot \int_{\R^n} ({\tilde c}_n \cdot H) \!\cdot\! \left( \,V_1 \ + \ V_2\, \right)^{{{n\,+\,2}\over {n\,-\,2}}   }
\cdot {{\partial V_1}\over {\partial \xi_{1_j}}}\\[0.15in]
& = & {{2n}\over {n\,-\,2}} \cdot \left[ \ \int_{\R^n} ({\tilde c}_n \cdot H) \!\cdot\!  V_1^{{{n\,+\,2}\over {n\,-\,2}}   }
\cdot {{\partial V_1}\over {\partial \xi_{1_j}}}   \  \right] \cdot [\, 1 \ + \ o\,(1)\,]\\[0.15in] \ \ \   & \ &  \!\!\!\!\!\!\!\!\!\! [\, {\mbox{similar \ \ to     \ \ Weak  \ \ Interaction \ \  Lemma \ \ 2.16\,;}}  \ \  \uparrow \ \ \ \  \ \ \  \ \ \ \downarrow \ \ {\mbox{cf. \ \ }} (A.1.30)  \,]\\[0.15in]
& = &   {{n - 2}\over 2} \cdot {{2n}\over {n\,-\,2}} \cdot \left[\,  \int_{\R^n} ({\tilde c}_n \cdot H) \!\cdot\!  V_1^{{{n\,+\,2}\over {n\,-\,2}}   }
\!\cdot\!  \left\{  -\, \lambda_1^{{n -\, 2}\over 2} \cdot {{ 2\,(\xi_{1_j} \,-\, y_j)}\over {(\,\lambda^2_1
\,+\, |\, y \,-\, \xi_1|^{\,2})^{{n }\over 2} }} \,\right\} \  \right] \cdot [\, 1 \, + \, o\,(1)\,]\\[0.15in]
& = &  -\,2\,n \cdot \left[ \  \int_{\R^n} ({\tilde c}_n \cdot H) \!\cdot\!
V_1^{{{n\,+\,2}\over {n\,-\,2}}   }
\cdot  V_1 \cdot \left\{ \   {{  (\xi_{1_j} - y_j)}\over {(\,\lambda^2_1
\,+\, |\, y \,-\, \xi_{1}|^{\,2})  }} \,\right\} \ dy \  \right] \cdot [\, 1 \ + \ o\,(1)\,]\\[0.15in]
& = &  -\,2\,n \cdot \left[ \,  \int_{\R^n} ({\tilde c}_n \cdot H) \!\cdot\!
 \left(  {{\lambda_1}\over {\lambda_1^2 \, + \, |\,y|^2}} \right)^{\!\!n}  \cdot \left\{ \   {{  (\xi_{1_j} - y_j)}\over {(\,\lambda^2_1
\,+\, |\, y \,|^{\,2} )  }} \,\right\} dy\,  \right] \!\cdot\! [\, 1 \, + \, o\,(1)\,] \ \ [\,{\mbox{cf. \ \ (A.10)}}\,]\\[0.15in]
& = &  -\,2\, n \cdot \xi_{1_j} \cdot \left[\   \int_{\R^n} ({\tilde c}_n \cdot H) \!\cdot\!
\left(  {{\lambda_1}\over {\lambda_1^2 \ + \ |\,y|^2}} \right)^{\!\!n}
  \cdot \left\{ \   {{  1}\over {(\,\lambda^2_1
\,+\, |\, y \,|^{\,2})  }} \,\right\} \ dy \right] \cdot [\, 1 \ + \ o\,(1)\,] \\[0.15in]
 & \ & \hspace*{4.7in}(\,{\mbox{via \ \ symmetry}}\,)\\[0.15in]
& = &  -\,2\, n \cdot {{ \xi_{1_j} }\over {\lambda_1}} \cdot \left[\   \int_{\R^n} ({\tilde c}_n \cdot H) \!\cdot\!
\left(  {{\lambda_1}\over {\lambda_1^2 \ + \ |\,y|^2}} \right)^{\!\!n}
  \cdot \left\{ \   {{  \lambda_1}\over {(\,\lambda^2_1
\,+\, |\, y \, |^{\,2})  }} \,\right\} \ dy \right] \cdot [\, 1 \ + \ o\,(1)\,] \ .
\end{eqnarray*}
It follows that
\begin{eqnarray*}
 & \ & -\,{{n\,-\,2 }\over {2n}} \cdot \left( \lambda_1 \cdot {\partial\over {\partial \xi_{1_j}}} \right)\,\int_{\R^n}  ({\tilde c}_n \cdot H) \!\cdot\! \left( \,V_1 \ + \ V_2\, \right)^{{{2n}\over {n\,-\,2}}   }\\[0.2in]
 & = &
(n\,-\,2) \cdot \lambda_1 \cdot {{ \xi_{1_j} }\over {\lambda_1}} \cdot  \left[\   \int_{\R^n} ({\tilde c}_n \cdot H) \!\cdot\!
 \left(  {{\lambda_1}\over {\lambda_1^2 \ + \ |\,y|^2}} \right)^{\!\!n\,+\,1}   dy \right] \cdot [\, 1 \ + \ o\,(1)\,]  \\[0.2in]
& =  &  (n\,-\,2) \cdot \lambda_1 \cdot  {{ \xi_{1_j} }\over {\lambda_1}} \cdot {{\lambda_1^\ell }\over { \lambda_1 }} \cdot \left[\   \int_{\R^n} ({\tilde c}_n \cdot H)\,(Y) \!\cdot\!
 \left(  {1\over {1\ + \ |\,Y|^2}} \right)^{\!\!n\,+\,1}   dY \right] \cdot [\, 1 \ + \ o\,(1)\,]  \\[0.2in]
& =  &  {\bar C}^+_2\,(n\,,\,\ell) \cdot  {\varpi}_1 \cdot   \lambda_1 \cdot {{ \xi_{1_j} }\over {\lambda_1}} \cdot {{\lambda_1^\ell }\over { \lambda_1 }} \cdot [\, 1 \ + \ o\,(1)\,]\   \ \ \ \ \ \  \ \ \   \ \ \     [\,{\mbox{similar \ \ to \ \ (4.6)}}\,; \ \ \ Y \ = \ \lambda_1^{-1} \cdot y\,].
\end{eqnarray*}
Here we use (1.8) and (1.10)\,.

\end{document}